\documentclass[reqno]{amsart}
\usepackage{tkz-euclide}
\usepackage{hyperref}
\usepackage{tikz-3dplot}
\usepackage{amsmath}  
\usepackage{amsfonts}
\usepackage{amssymb} 
\usepackage{mathrsfs} 
\usepackage{graphicx}  
\usepackage{bm}
\usepackage{tikz} 
\usetikzlibrary{calc} 
\usepackage{xcolor}  
\usepackage{a4wide}

\usetikzlibrary{arrows.meta} 
\usepackage{amsthm} 

\usepackage{float}
\usepackage{enumitem}
\usepackage[english]{babel}
\usepackage[font=footnotesize, labelfont=bf]{ caption}
\usepackage{ esint }
\usepackage{stackengine,scalerel,graphicx}
\stackMath

\definecolor{trueblue}{rgb}{0.0, 0.45, 0.81}

\newcommand{\EEE}{\color{black}}

\newcommand{\eps}{\varepsilon} 
\newcommand{\e}{\varepsilon}

\newcommand{\eg}{{\it e.g.}, }
\newcommand{\ie}{;{ \it i.e.}, }

\theoremstyle{plain}
\begingroup
\newtheorem{theorem}{Theorem}[section]
\newtheorem{lemma}[theorem]{Lemma}

\newtheorem{remark}[theorem]{Remark}
\newtheorem{proposition}[theorem]{Proposition}
\newtheorem{corollary}[theorem]{Corollary}
\endgroup

\theoremstyle{definition}
\begingroup

\endgroup

\renewcommand{\tilde}{\widetilde}

\DeclareMathOperator{\diam}{diam}

\renewcommand{\d}{ \mathrm{d}}

\numberwithin{equation}{section}
\newcommand{\N}{\mathbb{N}}
\newcommand{\Z}{\mathbb{Z}}
\newcommand{\F}{\mathcal{F}}

\newcommand{\R}{\mathbb{R}}

\renewcommand{\S}{\mathbb{S}}
\renewcommand{\L}{\mathcal{L}}

\newcommand{\dx}{\,\mathrm{d}x}
\newcommand{\dy}{\,\mathrm{d}y}
\newcommand{\ds}{\displaystyle}
\newcommand{\PP}{\mathbb{P}}
\let\O=\Omega
\let\o=\omega

\usepackage[normalem]{ulem}

\begin{document}

\title[Homogenisation of nonlinear Dirichlet problems in randomly perforated domains]{Homogenisation of nonlinear Dirichlet problems in randomly perforated domains under minimal assumptions on the size of perforations}

\author[Lucia Scardia]{Lucia Scardia} 
\address[Lucia Scardia]{Department of Mathematics, Heriot-Watt University, Edinburgh EH14 4AS, Scotland}
\email{L.Scardia@hw.ac.uk
}

\author{Konstantinos Zemas}
\address[Konstantinos Zemas]{Applied Mathematics M\"unster, University of M\"unster\\
	Einsteinstrasse 62, 48149 M\"unster, Germany}
\email{konstantinos.zemas@uni-muenster.de}

\author{Caterina Ida Zeppieri}
\address[Caterina Ida Zeppieri]{Applied Mathematics M\"unster, University of M\"unster\\
	Einsteinstrasse 62, 48149 M\"unster, Germany}
\email{caterina.zeppieri@uni-muenster.de}

\begin{abstract}
In this paper we study the convergence of integral functionals with $q$-growth in a randomly perforated domain of $\mathbb R^n$, with $1<q<n$.
Under the assumption that the perforations are small balls whose centres and radii are generated by a \emph{stationary short-range marked point process}, we obtain in the critical-scaling limit an averaged analogue of the nonlinear capacitary term obtained by Ansini and Braides in the deterministic periodic case \cite{Ansini-Braides}. In analogy to the random setting introduced by Giunti, H\"ofer, and Vel\'azquez \cite{Giunti-Hofer-Velasquez} to study the Poisson equation, we only require that the random radii have finite $(n-q)$-moment. This assumption on the one hand ensures that the expectation of the nonlinear $q$-capacity of the spherical holes is finite, and hence that the limit problem is well defined. On the other hand, it does not exclude the presence of balls with large radii, that can cluster up. We show however that the critical rescaling of the perforations is sufficient to ensure that no percolating-like structures appear in the limit.
\end{abstract} 

\maketitle

\section{Introduction}\label{introduction}
In this paper we study the limit behaviour as $\eps\to 0^+$ of sequences of \emph{nonlinear} Dirichlet boundary value problems of the type
\begin{equation}\label{intro:Pe}
\begin{cases}
\mathcal A[u]=0 & \text{in}\quad D_\eps
\cr
u=0 & \text{on} \quad \partial D_\eps\,, 
\end{cases}  
\end{equation}
where $\mathcal A$ is an elliptic differential operator and $D_\e:=D\setminus H_\eps$ is obtained by removing from an open, bounded, Lipschitz set $D\subset \R^n$ a collection $H_\e$ of small, spherical inclusions. 
Here we assume that $H_\e=H^\omega_\e$ is a \emph{random} set, namely, that the centres of the spherical holes are generated according to a stationary point process in $\R^n$ and that the associated radii are (suitably scaled) unbounded random variables with short-range correlations. 
Under minimal assumptions both on the nonlinearity and on the set of perforations, we prove that, \emph{almost surely}, solutions to \eqref{intro:Pe} converge weakly in a suitable Sobolev space to the solution of a limit problem
 \begin{equation}\label{intro:NL-PDE}
\begin{cases}
\mathcal A_0[u]=0 & \text{in}\quad D
\cr
u=0 & \text{on} \quad \partial D\,.
\end{cases}  
\end{equation}
In \eqref{intro:NL-PDE}, the nonlinear homogenised operator $\mathcal A_0$ depends both on $\mathcal A$ and on the geometry of the perforated domain, through some sort of limit averaged nonlinear capacity density of $H_\e^\omega$. 
In particular our result extends both the stochastic homogenisation result of Giunti, H\"ofer, and Vel\'azquez \cite{Giunti-Hofer-Velasquez} for the Poisson equation to the nonlinear setting, and the classical result of Ansini and Braides \cite{Ansini-Braides} to the random setting.  

\medskip
     
The study of homogenisation problems in perforated domains has a long history with seminal contributions of Marchenko and Khruslov \cite{Marchenko_Khruslov, Marchenko_Khruslov2}, Cioranescu and Murat \cite{Cionarescu_strange}, and Papanicolaou and Varadhan \cite{Papa} (see also \cite{Casado-Diaz-Garroni-2, DalMaso-Defranceschi, DalMaso-Garroni}). 
In a \emph{periodic} setting a typical $H_\e$ is chosen as
\begin{equation}\label{intro:periodic-holes}
H_\e = \bigcup_{i\in \Z^n} \overline{B}_{\eps^\alpha \rho}(\eps i)\,, 
\end{equation}
for some $\rho>0$ and $\alpha>1$. In \eqref{intro:periodic-holes} the parameter $\eps$ represents the characteristic distance between the centres of the spherical holes, while $\eps^\alpha \ll \e$ is proportional to the size of their (common) radius.  
If, moreover, the \textit{linear} case of the Poisson equation is considered, the boundary value problem \eqref{intro:Pe} becomes  
\begin{equation}\label{intro:Poi-e}
\begin{cases}
\mathcal -\Delta u=\psi & \text{in}\quad D_\eps
\cr
u=0 & \text{on} \quad \partial D_\eps\,, 
\end{cases}  
\end{equation}
where $\psi \in W^{-1,2}(D)$. 
In this linear, periodic framework, Cioranescu and Murat  \cite{Cionarescu_strange} showed the existence of a critical scaling for the perforation radius such that the sequence of solutions $(u_\eps)$ to \eqref{intro:Poi-e} converges weakly to the solution of a limit Dirichlet problem. Namely, assuming that $n>2$, for $\alpha=n/(n-2)$ and $\e>0$, the unique solution $u_\eps \in W^{1,2}_0(D_\eps)$ to \eqref{intro:Poi-e} converges weakly in $W^{1,2}(D)$, as $\e\to 0^+$, to the unique solution $u_0\in W^{1,2}_0(D)$ of
\begin{equation}\label{intro:Poi-lim}
\begin{cases}
\mathcal -\Delta u+\mu_0 u=\psi & \text{in}\quad D
\cr
u=0 & \text{on} \quad \partial D.
\end{cases}
\end{equation}
In \eqref{intro:Poi-lim} the zero-order term $\mu_0 u$ is reminiscent of the homogeneous Dirichlet boundary conditions prescribed on the boundary of the spherical holes in \eqref{intro:Poi-e}, and $\mu_0$ is a positive constant of geometric nature which represents the limit \emph{capacity density} generated by the set $H_\e$. Namely, we have
\begin{equation*}
\mu_0= \lim_{\e \to 0^+} \textrm{Cap}(H_\e\cap Q,\R^n)\,,
\end{equation*}
where $Q$ is a unit cube in $\R^n$ and $\textrm{Cap}(H_\e\cap Q,\R^n)$ denotes the elliptic, or $2$-capacity of $H_\e \cap Q$ in $\R^n$\ie
\[
\mathrm{Cap}(H_\e\cap Q,\R^n):= \inf\bigg\{\int_{\R^n}|\nabla v|^2\,\mathrm{d}x\colon  v\in W_0^{1,2}(\R^n),\ v\equiv 1 \ \text{on } H_\e\cap Q\bigg\}\,.
\]
In view of \eqref{intro:periodic-holes}, by the subadditivity of the capacity and the fact that $\e^\alpha \ll \e$ it is easy to see that
\begin{equation}\label{intro:mu}
\mu_0 =\lim_{\e \to 0^+} \textrm{Cap}\left(\bigcup _{i\in \Z^n} B_{\eps^\alpha \rho}(\eps i) \cap Q,\R^n\right)=\textrm{Cap}(B_\rho(0),\R^n)\,,
\end{equation}
where the last equality follows from the fact that the number of holes in $Q$ is of order $\e^{-n}$, and $\textrm{Cap}(B_{\e^\alpha \rho}(\e i),\R^n)= \e^n \,\textrm{Cap}(B_\rho(0),\R^n)$ for every $i\in \Z^n$, if $\alpha=n/(n-2)$. Moreover, in this case the constant $\mu_0$ is explicit and given by 
\[
\mu_0=(n-2) \mathcal H^{n-1}(\mathbb S^{n-1})\rho^{n-2}\,. 
 \]
Choices for the scaling of the perforation radius different from the critical value $\alpha=n/(n-2)$ give trivial convergence results. More precisely,  
in the case of \emph{tiny holes}, corresponding to the choice $\alpha>n/(n-2)$, it is immediate to see that $\mu_0=0$, so that the limit problem \eqref{intro:Poi-lim} reduces to 
\begin{equation*}
\begin{cases}
\mathcal -\Delta u=\psi & \text{in}\quad D
\cr
u=0 & \text{on} \quad \partial D\,.
\end{cases}  
\end{equation*}
For \emph{large holes}, corresponding to choosing $\alpha<n/(n-2)$, the sequence of solutions to \eqref{intro:Poi-e} converges strongly to zero in $W^{1,2}(D)$ (see \cite[Lemma 3.4.1]{Allaire}).  

In the last three decades the result of Cioranescu and Murat \cite{Cionarescu_strange} has been extended in a number of directions, ranging from the case of general \emph{nonlinear} elliptic operators  in periodically perforated domains \cite{Ansini-Braides, DalMaso-Defranceschi} to the case where a \emph{random} distribution of holes is also allowed \cite{Caffarelli_Mellet, Casado-Diaz, Casado-Diaz2, Focardi, Giunti, Giunti-Hofer-Velasquez}, just to mention a few examples.   

\begin{figure}
\begin{center}
\begin{tikzpicture}
\node[inner sep=0pt] (ran) at (0,0)
    {\includegraphics[width=.45\textwidth]{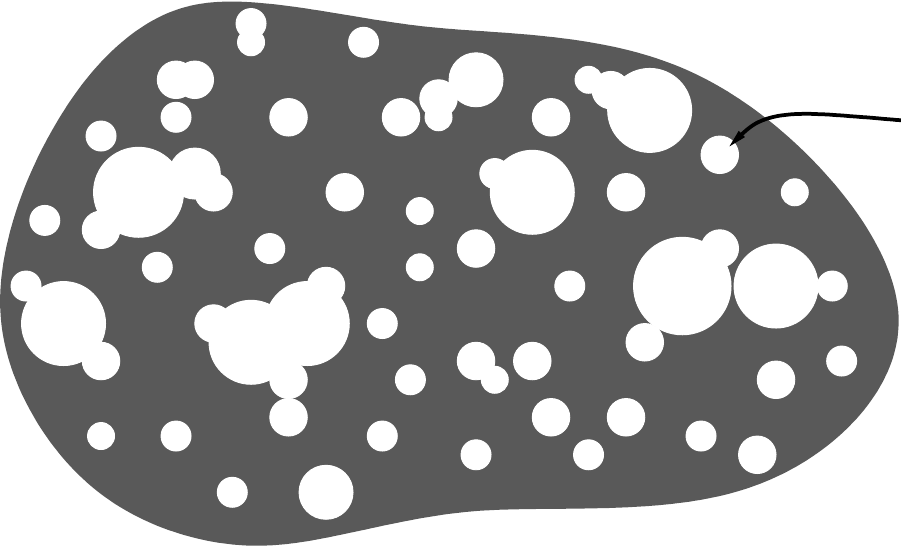}};
\draw(3,-1.6)node[right]{${D^\omega_\eps}$};
\draw(3.4,1.2)node[right]{${B}_{\e^\alpha \rho_i^\omega} (\e x_i^\omega)$};
\end{tikzpicture}
\end{center}
\caption{Realisation of the random domain $D_\e^\omega$.}\label{fig:realisation}
\end{figure}

As far as a nonlinear variant of \cite{Cionarescu_strange} is concerned, in \cite{Ansini-Braides} Ansini and Braides proved a \emph{nonlinear vectorial} version of the Cioranescu and Murat result when the Dirichlet boundary value problem \eqref{intro:Pe} is of variational nature\ie when $\mathcal A [u]=0$ in $D_\e$ is the Euler-Lagrange system associated to an integral functional defined on $D_\e$.  In \cite{Ansini-Braides} the corresponding limit problem is obtained by resorting to a direct $\Gamma$-convergence approach instead of the more classical PDE approach. The variational methods used in \cite{Ansini-Braides} allow for  minimal assumptions on the integral functionals and hence on the nonlinearity; these assumptions are the same considered in this paper and will be discussed below.  On the other hand, as far as the geometry of the perforated domain is concerned, in \cite{Giunti-Hofer-Velasquez} Giunti, H\"ofer, and Vel\'azquez proved a stochastic counterpart of \cite{Cionarescu_strange} for the linear problem \eqref{intro:Poi-e}, where $H_\e=H_\e^\omega$ is a \emph{random set} given by the union of small balls with random centres and radii, for which clusters occur with overwhelming probability (see Figure \ref{fig:realisation} for an illustration). In fact, the assumptions on $H_\e^\omega$ formulated in \cite{Giunti-Hofer-Velasquez} are shown to be the minimal ones in order to have homogenisation.  

The aim of this work is to combine the two general frameworks described above to devise \emph{minimal assumptions} both on the nonlinearity and on the random set of the spherical perforations $H_\e^\omega$, for which, almost surely, the corresponding Dirichlet problems \eqref{intro:Pe} admit a homogenised limit of the type \eqref{intro:NL-PDE}.

For the sake of the exposition, to illustrate our main result we consider here a prototypical random geometry for the set $H_\e^\omega$, while we refer the reader to Section \ref{Assumptions_on_m.p.p.} for the more general probabilistic framework considered in the paper.

In what follows $(\Omega, \mathcal T, \mathbb P)$ denotes a given underlying probability space. We consider a marked point process $(\Phi, \mathcal R)$ on $\R^n\times \R_+$ where  $\Phi:=(x_i^\omega)_i$ is a Poisson point process of constant intensity $\lambda>0$\ie the average number of points of the process per unit volume satisfies $\langle N(Q)\rangle =\lambda$, where $N(Q)=\# \big(\Phi \cap Q\big)$. For the marks we assume that  $\mathcal R=(\rho_i^\omega)_{x^\o_i\in \Phi}$, with $\rho_i$ identical and independently distributed \emph{unbounded} random variables.

The random set of perforations associated to $(\Phi, \mathcal R)$ is defined as
\begin{equation}\label{intro:He-sto} 
H_\e^\omega:= \bigcup_{x_i^\omega \in \Phi} \overline{B}_{\e^\alpha \rho_i^\omega} (\e x_i^\omega)\,,
\end{equation}
where $\alpha=n/(n-2)$ is the critical scaling for problems with quadratic growth. The analogue of \eqref{intro:mu} follows from the strong law of large numbers, which guarantees that 
\emph{almost surely}
\begin{equation}\label{intro:av-cap}
\lim_{\e \to 0^+}\textrm{Cap}(H_\e^\omega\cap Q,\R^n)= \lim_{\e \to 0^+} \e^n \hspace{-0.3cm}\sum_{x_i^\omega \in \Phi \cap (\e^{-1} Q)} \textrm{Cap}(B_{\rho^\omega_i}(\e x_i^\omega),\R^n)= 
\lambda \langle \textrm{Cap}(B_{\rho}(0),\R^n)\rangle\,.
\end{equation}
Moreover, an explicit calculation gives 
\begin{equation}\label{intro:tilde-mu}
\tilde \mu_0= \lambda \langle \textrm{Cap}(B_{\rho}(0),\R^n)\rangle= (n-2) \mathcal H^{n-1}(\mathbb S^{n-1})\lambda \langle \rho^{n-2}\rangle\,,
\end{equation}
where clearly $\tilde \mu_0$ reduces to $\mu_0$ if $\Phi=\Z^n$ and $\rho_i^\omega=\rho$, constant and deterministic.  

We observe that in view of \eqref{intro:tilde-mu}, the minimal condition for $\tilde \mu_0$ to be well-defined is that the random variables $\rho_i$ satisfy the stochastic integrability condition 
\begin{equation}\label{intro:moment}
\langle \rho_i^{n-2} \rangle<+\infty\,. 
\end{equation}
We note, however, that \eqref{intro:moment} does not prevent the balls generating $H^\o_\e$ from overlapping. Indeed, it is easy to check that the expected number of holes that may overlap (that is, for which $\e^{n/(n-2)} \rho_i^\omega > \e$, with $\lambda \e$ being the typical distance between two points in $\Phi$) is of the order $\e^{-n+2}$, while the expected total number of holes is of order $\e^{-n}$. In \cite{Giunti-Hofer-Velasquez} Giunti, H\"ofer, and Vel\'azquez proved that, even though with probability one there are regions where the holes cluster, the moment condition \eqref{intro:moment} is indeed sufficient to ensure that almost surely these regions have a vanishing capacity, as $\e \to 0^+$.  Moreover, this moment condition allows to extend the Cioranescu and Murat construction of the ``oscillating test functions'' to this random setting and to prove that the stochastic analogue of \cite{Cionarescu_strange} holds true, almost surely, up to replacing in \eqref{intro:Poi-lim} $\mu_0$ with $\tilde \mu_0$.         
  
 \medskip 

In this paper we extend the result by Giunti, H\"ofer, and Vel\'azquez to the nonlinear vectorial setting. In the same spirit as in \cite{Ansini-Braides}, we work with functionals rather than with the associated Euler-Lagrange systems. Therefore, we consider $1<q<n$, $m\in \N$, and a Borel-measurable function $f \colon \R^{m\times n} \to \R$, with $f(0)=0$, satisfying a $q$-growth and coercivity condition\ie   
\begin{equation*}
c_1(|\xi|^q-1)\leq f(\xi)\leq c_2(|\xi|^q+1)\, \quad \forall\, \xi\in \R^{m\times n}
\end{equation*} 
with $0<c_1<c_2$. Then, we introduce the vectorial, random functionals defined on $W_0^{1,q}(D;\R^m)$ as
\begin{equation}\label{intro:Fe}
\F_\eps^\omega(u):= \int_{D\setminus H^\omega_\e}f(\nabla u)\,\mathrm{d}x,\quad   \text{ if }\ u=0 \; \text{ in }\; H_\eps^\omega\cap D\,,
\end{equation} 
and $+\infty$ otherwise. In \eqref{intro:Fe}, $H^\omega_\e$ is as in \eqref{intro:He-sto} with $\alpha$ being the critical scaling in the case of $q$-growth, namely $\alpha=n/(n-q)$.
Moreover, for every $\xi\in \R^{m\times n}$ we set
\[
g(\xi)=\lim_{\e \to 0^+} \eps^{nq/(n-q)}
Qf(\eps^{-n/(n-q)}
\xi)\,,
\]
where $Qf$ denotes the quasiconvex envelope of $f$. Upon passing to a subsequence the function $g$ is always well-defined (see Section \ref{main_result} for more details), moreover the limit in $\e$ becomes redundant when $f$ (and hence $Qf$) is positively homogeneous of degree $q$.   

Eventually, for every $z\in \R^m$ we define the random variables 
\begin{equation*}
\hspace{-0.8em}\varphi_{\rho_i^\omega}(z):=\inf\bigg\{\int_{\R^n}g(\nabla\zeta)\,\mathrm{d}x\colon \zeta-z\in W_0^{1,q}(\R^n;\R^m),\ \zeta\equiv 0 \ \text{on } \overline{B}_{\rho^\omega_i}(0)\bigg\}\,.
\end{equation*}
We observe that for $\omega \in \Omega$ fixed, $\varphi_\rho(z)$ can be interpreted as a \emph{nonlinear $g$-capacity density} of the ball $B_{\rho}(0)$ in $\R^n$; furthermore, for $f(\xi)=|\xi|^q$ we have
\[
\varphi_\rho(z)=\mathrm{Cap}_q(B_{\rho}(0),\R^n) |z|^q \,,
\]
where $\textrm{Cap}_q(B_{\rho}(0),\R^n)$ denotes the classical $q$-capacity of $B_\rho(0)$ relative to $\R^n$ (cf.\ \eqref{relative_scalar_q_capacity}). In general, however, $\varphi_\rho$ is not positively $q$-homogeneous as a function of $z$ (as observed by Casado-Diaz and Garroni in \cite{Casado-Diaz-Garroni-2}), hence in the definition of $\varphi_\rho(z)$ (and in the case $m=1$) we cannot reduce to the case $z=1$ just by unscaling.

Thanks to the growth conditions of order $q$ satisfied by $f$, it is easy to show that 
\begin{equation}\label{intro:n-q-estimate}
C_1|z|^q\rho^{n-q}\leq \varphi_{\rho}(z)\leq C_2|z|^q\rho^{n-q}\,, 
\end{equation}  
for some $0<C_1<C_2$ only depending on $n, q, c_1$, and $c_2$ (see Lemma \ref{growth_of_phi_rho}). Hence, in analogy to \eqref{intro:tilde-mu}-\eqref{intro:moment}, in this nonlinear framework the homogenised problem is well-defined up to assuming the stochastic integrability condition
\begin{equation}\label{intro:stoch-int}
\langle \rho^{n-q} \rangle <+\infty
\end{equation}
for the perforations radii.

The main result of this paper is Theorem \ref{main_thm}, which establishes an \emph{almost sure} $\Gamma$-convergence result for the random functionals $\mathcal F^\omega_\e$, under the assumption \eqref{intro:n-q-estimate}. Namely, Theorem \ref{main_thm} states that if $(\Phi, \mathcal R)$ is a marked point process as above whose marks additionally satisfy the moment condition \eqref{intro:stoch-int}, then, almost surely, the functionals $\mathcal F_\e^\omega$ $\Gamma$-converge, as $\eps\to 0^+$, to the \emph{deterministic} functional given by 
\[
\mathcal F_0(u):=\int_{D}Qf(\nabla u)\,\mathrm{d}x+\lambda\,\int_D\langle \varphi_\rho (u)\rangle\,\mathrm{d}x, \quad \text{if} \quad u\in W_0^{1,q}(D;\R^m)\,.
\]
In particular, if $f$ is convex and differentiable, by the fundamental property of $\Gamma$-convergence we deduce that for $\mathbb P$-a.e.\ $\omega \in \Omega$ the unique solution $u^\omega_{\eps} \in W^{1,q}_0(D;\R^m)$ to 
\[
\begin{cases}
-\mathrm{div}Df(\nabla u^\omega_{\eps})=\psi & \text{in } D^\omega_{\eps} \cr
u^\omega_{\eps}=0  & \text{on } \partial D^\omega_{\eps}
\end{cases}
\]
converges weakly in $W^{1,q}(D;\R^m)$ to the unique solution $u \in W_0^{1,q}(D;\R^m)$ of the boundary value problem 
\begin{equation*}
\begin{cases}
-\mathrm{div}Df(\nabla u)+\lambda \varphi'(u)=\psi & \text{in } D \cr
u=0 & \text{on } \partial D\,,
\end{cases}
\end{equation*}
where $\varphi(z):=\langle \varphi_{\rho}(z)\rangle$ and $\psi \in W^{-1,\frac{q}{q-1}}(D; \R^m)$. Thus, in particular, for $f(\xi)=|\xi|^q$ we obtain a convergence result for the $q$-Laplace equation in randomly perforated domains of general geometry.

To prove the $\Gamma$-convergence result in Theorem \ref{main_thm} we follow a proof methodology which is close in spirit to that of Ansini and Braides \cite{Ansini-Braides}. This approach is purely variational and is based on a so-called ``joining lemma on perforated domains'' which allows to replace any sequence with equibounded energy $\F_\eps^\omega$ with a sequence which is constant in a spherical layer surrounding each perforation, without essentially changing the energy. In \cite{Ansini-Braides} this construction  is then pivotal both in the proof of the lower-bound and in the upper-bound inequalities. Indeed, when proving the lower-bound inequality the joining lemma allows to estimate separately the energy close to and far from the perforations. Moreover, using the modified sequence it is possible to recover the nonlinear capacitary term as the limit of some suitable discrete energy densities. On the other hand, in the proof of the upper-bound inequality the joining lemma enters in the construction of the recovery sequence.  

Crucially, in \cite{Ansini-Braides} to prove the joining lemma and thus the $\Gamma$-convergence result, it is of fundamental importance that the perforations are well-separated from one another. In the periodic setting this is a straightforward consequence of the regular arrangement of the holes on an $\e$-scale together with the $\e^\alpha \ll \e$ scaling for their (constant) radius.  Whereas, as already observed, in our framework for almost every realisation we have to deal with the presence of large radii or centres very close to each other and thus of clustering holes. This technical issue is tackled similarly as in \cite{Giunti-Hofer-Velasquez} by showing that, almost surely, the set of perforations $H_\e^\omega$ can be partitioned into two sets: a set of ``good'' holes $H^\omega_{\e,g}$ and a set of ``bad'' holes $H^\omega_{\e,b}$. In the set of good holes we identify a subset of  ``small'' balls which are ``$\e$-separated'' from one another, in which then a stochastic variant of the joining lemma holds. There are the only balls that contribute to the nonlinear capacitary term in $\mathcal F_0$  (cf.\ Proposition \ref{main_prop_for_approximating_the_capacity}). The set of bad holes $H^\omega_{\e,b}$ contains, among others, all those balls which overlap with probability one. We then show that $H_{\e,b}^\omega$ can be enclosed into a ``safety layer'' $D_{\e,b}^\omega$ which is still well-separated from $H_{\e,g}^\omega$, and such that the nonlinear $g$-capacity of $H^\omega_{\e,b}$ relative to $D^\omega_{\e,b}$ vanishes as $\e\to 0^+$. In other words, the set of bad holes is well-separated from $H^\omega_{\e,g}$ and is asymptotically negligible. 

Similarly, in the construction of the recovery sequences for the $\Gamma$-limit the only energy contribution that is relevant for the capacitary term  is the one carried by the balls in $H^\omega_{\e,g}$. However, since a recovery sequence needs to be admissible for $\F_\e^\omega$, it has to vanish also in the bad balls and hence in particular in the clusters. This constraint makes for a rather delicate proof of the upper bound which also relies on a corrector-like construction in the bad holes for the $q$-Laplace equation (cf.\ Lemma \ref{construction_of_correctors}).             

\medskip

This paper is organized as follows. In Section \ref{setting} we recall the basics of marked point processes $(\Phi, \mathcal R)$, and we list the assumptions that the process generating the holes $H_\e^\omega$ needs to satisfy. These assumptions are quite mild: our analysis is valid for rather general stationary point processes $\Phi$ whose associated marks in $\mathcal R$ need not be independent, as long as their correlation-range is suitably controlled (cf.\ \eqref{2_correlation_functions} and \eqref{short_range_integrable_correlations}).  In Section \ref{main_result} we state the main result of this paper, Theorem \ref{main_thm}.  To prove it, we need a number of technical results of both deterministic and probabilistic nature. The deterministic preliminaries are collected in Section \ref{deterministic_preliminaries} which contains, among other results, a variant of the joining lemma which is relevant in our case (cf. Lemma \ref{joining_lemma_on_randomly_perforated_domains}). Section \ref{probabilistic_preliminaries} is instead entirely devoted to some stochastic auxiliary results. In this section we prove, in particular, a version of the strong law of large numbers for correlated marked point processes in the nonlinear setting (cf.\ Proposition \ref{probabilistic_lemma_3}). In Section \ref{discrete_capacity_approximation} we build upon sections \ref{deterministic_preliminaries} and \ref{probabilistic_preliminaries} to prove a discrete approximation result for the nonlinear capacitary term (cf.\ Proposition \ref{main_prop_for_approximating_the_capacity}). Section \ref{s:main-thm-proof} is devoted to the proof of the $\Gamma$-convergence result Theorem \ref{main_thm} (cf.\ Proposition \ref{newprop_inf} and Proposition \ref{Gamma_limsup_ineq}).       

\section{Notation and setup}\label{setting} 

In this section we collect some useful notation and we introduce the probabilistic setup. 

\subsection{Notation} We denote with $\mathcal{B}(\R^n)$ the $\sigma$-algebra of Borel subsets of $\R^n$. For every $A\subset \R^n$ and $x\in \R^n$, we denote by $\tau_xA$ the shift of the set $A$ by $x$\ie $\tau_xA:=x+A$. The diameter of $A$ is denoted with $\mathrm{diam}(A)$ and the characteristic function of $A$ with $\chi_A$. Given $A, B\subset \R^n$ we write $A \subset \subset B$ iff $\overline{A} \subset B$. By $\#A$ we denote the cardinality of a discrete set $A$. For $n,k\in \mathbb{N}$ we denote by $\mathcal{L}^n$ and $\mathcal{H}^{k}$ the $n$-dimensional Lebesgue measure and the $k$-dimensional Hausdorff measure, respectively. 

Given $\rho>0$ and $x\in \R^n$, we denote with $B_\rho(x)$ the open ball centred at $x$ with radius $\rho$. (We also use the notation $B^n_\rho(x)$ to clarify the dimension, if needed.) We denote with $Q_{\rho}(x)$ the half-open cube centred at $x$ with side-length $\rho>0$, namely
$$Q_{\rho}(x):=x+\rho[-\tfrac{1}{2},\tfrac{1}{2})^{n}\,,$$ 
and we omit the subscript when $\rho=1$, so that $Q(x)=Q_1(x)$. The unit sphere in $\R^n$ is denoted with $\S^{n-1}$. 
Moreover, we use the notation $\beta_n:=\L^n(B_1(0))$ for the volume of the unit ball in $\R^n$. 

For every $a,b\in \R$ we use the standard notation $a\wedge b:=\min\{a,b\}$ and $a\vee b:=\max\{a,b\}$. From time to time we use the notation $\varlimsup$ and $\varliminf$ to indicate the $\limsup$ and $\liminf$ respectively. The symbols $\sim_{M_1,M_2,\dots}, \lesssim_{M_1,M_2,\dots}$ indicate that the corresponding equality, inequality is valid up to a (positive multiplicative) constant that depends only on the parameters $M_1, M_2,\dots$ and the space dimension, but is allowed to vary from line to line. 

Let $(\Omega,\mathcal{T},\mathbb{P})$ denote an underlying given probability space; the expected value of a random variable $X:\Omega\to \mathbb R$ with respect to the probability measure $\mathbb{P}$ is denoted by $\langle X\rangle$\ie
\begin{equation*}
\langle X\rangle :=\int_{\Omega} X(\omega)\,\mathrm{d}\mathbb{P}(\omega)=:\int_{\R} x\,\mathrm{d}(X_{*}\mathbb{P})(x)\,,
\end{equation*}
where $X_{*}\mathbb{P}$ is the push-forward measure of $\mathbb{P}$ onto $\R$, or the probability distribution of $X$, defined via $(X_{*}\mathbb{P})(B):=\mathbb{P}(X^{-1}(B))$, for every $B\in \mathcal B(\R)$. 

\subsection{Marked point processes}  We refer the reader to \cite[Chapter 9, Definitions 9.1.I - 9.1.IV]{Daley-Vere-Jones} and \cite[Section 3.5]{Schneider-Weil} for a systematic treatment of marked point processes.  

Throughout the paper $(\Phi,\mathcal{R})$ denotes a marked point-process (in short, m.p.p.) with state space $\Phi \subset \R^n$ and mark space $\mathcal{R}\subset \R_+$. Then, an element $(x,\rho)$ of $(\Phi,\mathcal{R})$ is referred to as a point-mark pair and is given by random variables $x:\Omega\to \R^n$ and $\rho:\Omega\to \R_+$. For a fixed realisation $\omega\in \Omega$, the set $\Phi^\omega:=\Phi(\omega) \subset \R^n$ is a countable collection of points in $\R^n$ and $\mathcal{R}^\omega:=\mathcal{R}(\omega) \subset \R_+$ are the associated marks. By the definition of point processes the set $\Phi^\omega \cap B$ is \emph{finite} for every bounded $B \in \mathcal{B}(\R^n)$.

Let $P^{(\Phi,\mathcal R)}$ denote the probability distribution of the m.p.p.\ $(\Phi,\mathcal R)$, which is a probability measure on the product space $(\R^n\times \R_+,\mathcal{B}(\R^n)\otimes \mathcal{B}(\R_+))$ and let  $P^{\Phi}$ denote its marginal obtained by projecting $P^{(\Phi,\mathcal R)}$ onto $(\R^n,\mathcal{B}(\R^n))$. By \cite[Proposition A.1.5. III]{Daley-Vere-Jones1} there exists a family of probability kernels $\mathcal K: \R^n\times  \mathcal{B}(\R_+)\to [0,1]$ (namely a family of probability measures on $\R_+$ for every $x\in \R^n$, where $\mathcal K(\cdot,B)$ is $\mathcal{B}(\R^n)$-measurable for every fixed $B\in \mathcal{B}(\R_+)$) such that 
$$
P^{(\Phi,\mathcal R)}(A\times B)=\int_{A}\mathcal K(x,B)\, \mathrm{d}P^{\Phi} (x) \quad \text{for all } A\in \mathcal{B}(\R^n) \, \text{and } B\in \mathcal{B}(\R_+)\,.
$$
Then $P_{x}^\mathcal{R}(\cdot):=\mathcal K(x,\cdot)$ is called the regular conditional probability distribution of the marks or \emph{1-point mark distribution}. In practice, it is obtained by disintegrating $P^{(\Phi,\mathcal R)}$ with respect to $P^{\Phi}$. Its density with respect to the one-dimensional Lebesgue measure is denoted with $f_1(x,\rho)$, namely 
\begin{equation}\label{c:f-uno}
P^\mathcal{R}_x(B)=\int_{B} f_1(x,\rho)\, \d\rho \quad \text{for all }\,  B\in \mathcal{B}(\R_+)\,.
\end{equation}
Similarly, we denote with $P_{x_1,x_2}^\mathcal{R}$ the \emph{2-point mark distribution}, namely the regular conditional joint probability distribution of two marks in $\mathcal{R}$, given that $\Phi$ has two points at the locations $x_1,x_2\in \R^n$. We denote its density with respect to the two-dimensional Lebesgue measure with $f_2((x_1,\rho_1),(x_2,\rho_2)\big)$, namely 
\begin{equation}\label{c:f-due}
P^\mathcal{R}_{x_1,x_2}(B\times C)=\int_{B}\int_C f_2\big((x_1,\rho),(x_2,\rho') \big)\d\rho\,\d\rho' \quad 
\quad \text{for all }\,  B, C\in \mathcal{B}(\R_+)\,.
\end{equation}

\subsection{Assumptions on the m.p.p. $(\Phi,\mathcal{R})$}\label{Assumptions_on_m.p.p.} 
Below we list the assumptions we require on the m.p.p.\ $(\Phi,\mathcal{R})$; these are in the same spirit as the ones formulated in \cite[Section 2]{Giunti-Hofer-Velasquez}. 

\begin{enumerate}[label=(H$\arabic*$)]
\item\label{H1} The point process $\Phi$ is \textit{stationary}\ie for every $x\in \R^n$ the processes $\tau_x\Phi:=\{x+x_i\}_{x_i\in \Phi}$ and $\Phi$ have the same probability distribution.

\smallskip

\item\label{H2} The point process $\Phi$ is \emph{locally square integrable}\ie
there exists $\lambda>0$ such that for every unitary cube $Q\subset \R^n$,
\begin{equation}\label{second_moment_condition}
\big\langle (N(Q))^2\big\rangle\leq \lambda^2\,,
\end{equation}
where $N(B): =\# \Phi(B)$, with $\Phi(B):=\Phi \cap B$, for every $B\in \mathcal B(\R^n)$. 
(We note that the stationarity of $\Phi$ ensures that the bound in \eqref{second_moment_condition} is independent of the centre of the cube $Q$.) 

\smallskip

\item\label{H3} The point process $\Phi$ satisfies a \textit{strong mixing condition.} For $A\in \mathcal B(\R^n)$, let $\mathcal{T}(A)$ denote the smallest $\sigma$-algebra with respect to which the random variables $N(B)$ are $\mathbb{P}$-measurable for every Borel subset $B\subset A$. We assume that there exist constants $C>0$ and $\gamma>n$ with the following property. For every $A\in \mathcal B(\R^n)$, every $x\in \R^n$ with $|x|>\mathrm{diam}(A)$, and for every random variables $Z_1,Z_2$ measurable with respect to $\mathcal{T}(A), \mathcal{T}(\tau_xA)$ respectively, there holds 
\begin{equation}\label{strong_mixing}
|\langle Z_1Z_2\rangle-\langle Z_1\rangle\langle Z_2\rangle|\leq \frac{C}{1+(|x|-\mathrm{diam}(A))^\gamma}\langle Z_1^2\rangle^{1/2}\langle Z_2^2\rangle^{1/2}\,.
\end{equation}
We observe that \eqref{strong_mixing} in particular ensures the ergodicity of $\Phi$, cf. \cite[Paragraph 12.3]{Daley-Vere-Jones}.

\smallskip

\item\label{H4} Let $f_1$ and $f_2$ be as in \eqref{c:f-uno} and \eqref{c:f-due}, respectively. 
In view of the stationarity of $\Phi$, the density $f_1$ is independent of $x$\ie for every $x\in \R^n$ we have 
$f_1(x,\rho)=h(\rho)$, for some $h\in L^1(\R_+;\R_+)$ with $\int_{\R_+}h(\rho)\,\mathrm{d}\rho=1\,.$

We assume that $h$ satisfies the following integrability condition:
\begin{equation}\label{finite_average_capacity}
\int_{0}^{+\infty}\rho^{n-q}h(\rho)\,\mathrm{d}\rho<+\infty\,,
\end{equation}
which is equivalent to asking $\langle\rho^{n-q}\rangle<+\infty$. Moreover, for the density $f_2$ we assume that
\begin{equation}\label{2_correlation_functions}
f_2\big((x_1,\rho_1),(x_2,\rho_2)\big)=h(\rho_1)h(\rho_2)+K(|x_1-x_2|,\rho_1,\rho_2)\,,
\end{equation} 
for some function $K$ satisfying 
\[
\int_{\R_+\times \R_+}K(r,\rho_1,\rho_2)\,\mathrm{d}\rho_1\,\mathrm{d}\rho_2=0 \quad \text{for every $r\geq 0$}\,,
\] 
and 
\begin{equation}\label{short_range_integrable_correlations}
|K(r,\rho_1,\rho_2)|\leq \frac{C}{(1+r^\gamma)(1+\rho_1^s)(1+\rho_2^s)}\,,
\end{equation}
for some $C>0$ and $\ s>n-q$, and where $\gamma>n$ is the constant in \ref{H3}. 
\end{enumerate}
The interested reader is referred to \cite[Subsection 2.1]{Giunti-Hofer-Velasquez} for some explicit examples of m.p.p. $(\Phi,\mathcal{R})$ satisfying \ref{H1}--\ref{H4}.

\begin{remark}[Independent marking]
\normalfont
Under our assumptions the marks have the same distribution, but they are not independent. If the m.p.p. is in addition independently marked, then the expression of $f_2$ simplifies to 
$$
f_2\big((x_1,\rho_1),(x_2,\rho_2)\big)=h(\rho_1)h(\rho_2)\,.
$$
The additional (location-dependent) term in \eqref{2_correlation_functions} introduces a \textit{short-range correlation} between the marks conditioned on the point positions, thus giving a measure of the lack of independence of the marks. 
We note that if, conditional to the point process $\Phi$, all the marks are independent but not necessarily identically distributed, then the density of the 2-point mark distribution still factorises as 
$$
f_2\big((x_1,\rho_1),(x_2,\rho_2)\big) = f_1(x_1,\rho_1)\,f_1(x_2,\rho_2)\,,
$$
but is location-dependent.
\end{remark}

\begin{remark}
\normalfont As in \cite{Giunti-Hofer-Velasquez}, the assumptions \ref{H1}--\ref{H4} guarantee the validity of the \textit{strong law of large numbers-type results} stated in Subsection \ref{probabilistic_preliminaries} (see Lemmata \ref{probabilistic_lemma_1}--\ref{probabilistic_lemma_3} therein), which will play an important role in what follows.  
\end{remark}

\section{Statement of the main result}\label{main_result}

In this section we state the main result of the paper. To this end, we need to introduce some additional notation. 

Let $n\in \N$ and $1<q<n$, and let $(\Phi, \mathcal R)$ be a m.p.p.\ satisfying \ref{H1}--\ref{H4}. For $\omega \in \Omega$ fixed, we consider a countable family of points $(x^\omega_i)_i\subset \Phi^\omega$ and the corresponding marks $(\rho_i^\omega)_i\subset \mathcal{R}^\omega$. For fixed $\eps>0$, we associate to $(x^\omega_i, \rho_i^\omega)_i$  the family of open balls $\big({B}_{\alpha_\eps
\rho^\omega_i}(\eps x_i^\omega)\big)_i$\,, where 
\begin{equation}\label{eps_scaling}
\alpha_\eps:=\eps^{n/(n-q)}\,,
\end{equation}  
and where we have used the shorthand notation $x^\omega_i:=x_i(\omega)$ and $\rho^\omega_i:=\rho_i(\omega)$. 

Let $D\subset \R^n$ be an open, bounded, Lipschitz set, star-shaped with respect to the origin. 
The set of random spherical perforations in $D$ is given by
\begin{equation}\label{random_holes}
H_\eps^\omega:=\bigcup_{x^\omega_i\in \Phi^\omega\cap({\eps}^{-1}D)} \overline{B}_{\alpha_\eps	
	\rho^\omega_i}(\eps x_i^\omega)\,;
\end{equation}
note that the sets $({\eps}^{-1}D)_{\eps>0}$ are nested as $\eps\to 0^+$, since $D$ is star-shaped with respect to the origin. We finally define the randomly perforated domain as
\begin{equation*}
D_\eps^\omega:=D\setminus H_\eps^\omega\,.
\end{equation*}
Now, let $m\in \N$, and let $f:\R^{m\times n}\to 
\R$ be a Borel-measurable function of $q$-growth\ie there exist two constants $0<c_1<c_2$, such that
\begin{equation}\label{growth_of_f}
c_1(|\xi|^q-1)\leq f(\xi)\leq c_2(|\xi|^q+1) \ \forall\, \xi\in \R^{m\times n} \,. 
\end{equation} 
Without loss of generality we assume that $f(0)=0$. 

Finally, we introduce the nonlinear vectorial (random) functionals $\F_\eps^\omega:L^{1}(D;\R^m)\longrightarrow  
\R\cup\{+\infty\}$ defined as
\begin{equation}\label{random_functionals}
\F_\eps^\omega(u):=\begin{cases} \ds \int_{D_\eps^\omega}f(\nabla u)\,\mathrm{d}x & \text{if } u\in W_0^{1,q}(D;\R^m) \, \text {and } u=0 \ \text{in } H_\eps^\omega\cap D\,,
\cr
+\infty\, & \text{otherwise in } L^1(D;\R^m)\,.
\end{cases}
\end{equation} 
We recall that the $\alpha_\eps$-scaling \eqref{eps_scaling} for the radii of the perforations is the critical one in the case of energies with $q$-growth, under Dirichlet boundary conditions (cf.\ \cite{Ansini-Braides}).

The aim of this paper is to study the limit behaviour of the functionals $\F_\eps^\omega$ as $\eps\to 0^+$ (see Theorem \ref{main_thm}).
 
\begin{remark}\label{on_boundary_conditions}
\normalfont
 The homogeneous trace condition imposed on $\partial D$ in \eqref{random_functionals} allows us to work under the additional assumption $H_\eps^\omega\subset\subset D$, for $\eps>0$ small enough. In fact, if $D'\supset \supset D$ is open, bounded, Lipschitz, and star-shaped with respect to the origin, and is such that $H_\eps^\omega\subset\subset D'$, then, thanks to the condition $u\in W_0^{1,q}(D;\R^m)$, we can readily extend $u$ to $D'$ by setting $u\equiv0$ in $D'\setminus D$, thus getting $u=0$ on $H_\eps^\omega \cap D'$. We will then assume with no loss of generality that we are always in this situation.
\end{remark}

We now define the nonlinear capacitary term which appears in the limit functional of $\F_\eps^\omega$.

Let $(\eps_j) \searrow 0$; we define the functions $g_j:\R^{m\times n}\to \R$ as
\begin{equation}\label{j_capacity_densities}
g_j(\xi):=\alpha_{\eps_j}^q
Qf({\alpha_{\eps_j}^{-1}
}\xi) \quad \forall \, \xi\in \R^{m\times n}\,,   
\end{equation}
where $\alpha_\eps$ is as in \eqref{eps_scaling} and $Qf$ denotes the \emph{quasiconvex envelope} of $f$\ie
\begin{equation*}
Qf(\xi):=\inf\bigg\{\int_{(0,1)^n}f(\xi+\nabla\psi)\,\mathrm{d}x: \ \psi\in W_0^{1,q}((0,1)^n;\R^m)\bigg\}\,.
\end{equation*}
We note that $g_j$ is quasiconvex for every $j\in \N$; moreover, by \eqref{growth_of_f}, 
\begin{equation}\label{growth_of_g_j}
c_1(|\xi|^q-\alpha_{\eps_j}^q 
)\leq g_j(\xi)\leq c_2(|\xi|^q+\alpha_{\eps_j}^q
)
\end{equation} 
for every $\xi \in \R^{m\times n}$. Invoking, \eg \cite[Remark 4.13]{Braides-De Franceschi}, we then get that $(g_j)$ are locally equi-Lipschitz continuous, namely
\begin{equation}\label{equilipschitzianity}
|g_j(\xi_1)-g_j(\xi_2)|\leq L(\alpha_{\eps_j}^{q-1}
+|\xi_1|^{q-1}+|\xi_2|^{q-1})|\xi_1-\xi_2|\,,
\end{equation} 
for every $\xi_1,\xi_2 \in \R^{m\times n}$, every $j\in \N$, and for some constant $L:=L(c_1,c_2,q)>0$. Consequently, up to a subsequence (not relabelled), for every $\xi \in \R^{m\times n}$ there exists the limit
\begin{equation}\label{limit_capacity_density}
g(\xi):=\lim_{j\to +\infty} g_j(\xi)\,.
\end{equation}
The function $g$ is quasiconvex, and in view of \eqref{growth_of_g_j} and \eqref{equilipschitzianity} it satisfies the growth conditions
\begin{equation}\label{growth_of_g}
c_1|\xi|^q\leq g(\xi)\leq c_2|\xi|^q \ \ \forall \, \xi\in \R^{m\times n}\,,
\end{equation} 
as well as the bound
\begin{equation}\label{lipschitzianity_of_g}
|g(\xi_1)-g(\xi_2)|\leq L(|\xi_1|^{q-1}+|\xi_2|^{q-1})|\xi_1-\xi_2|\ \ \forall \, \xi_1, \xi_2 \in\R^{m\times n}\,. 
\end{equation}
Given two open sets $A\subset \subset B\subset \R^n$, with $A$ bounded, we define the $q$-capacity of $A$ relative to $B$ as 
\begin{equation}\label{relative_scalar_q_capacity}
\mathrm{Cap}_q(A,B):= \inf\bigg\{\int_{\R^n}|\nabla v|^q\,\mathrm{d}x\colon  v\in W_0^{1,q}(B;\R),\ v\equiv 1 \ \text{on } \overline{A}\bigg\}\,. 
\end{equation}
In the definition above $\mathrm{Cap}_q(A,B)$ depends on $A$ only via its closure $\overline A$, which is compact since $A$ is bounded. Hence \eqref{relative_scalar_q_capacity} agrees with the classical definition of $q$-capacity for compact sets, see \eg  \cite{Federer-Ziemer}. Note that, if $A\subset A'\subset\subset B\subset B'$, then 
\begin{equation}\label{Capq-increasing}
\mathrm{Cap}_q(A,B')\leq \mathrm{Cap}_q(A,B) \leq \mathrm{Cap}_q(A',B)\,.
\end{equation}
Moreover, $\mathrm{Cap}_q$ is countably subadditive with respect to the first entry (see, e.g. \cite{Federer-Ziemer}, Section 3), namely 
\begin{equation}\label{Capq-subadditive}
\text{if } A\subset \bigcup_{i=1}^\infty A_i, 
\quad \text{then } \,
\mathrm{Cap}_q(A,B) \leq \sum_{i=1}^\infty \mathrm{Cap}_q(A_i,B)\,.
\end{equation}
Given $A\subset \subset B\subset \R^n$ and $z\in \R^m$, we define the $g$-\textit{capacity} of $A$ relative to $B$ at $z$ as
\begin{equation*}
\mathrm{Cap}_g(A,B;z):=\inf\bigg\{\int_{\R^n}g(\nabla\zeta)\,\mathrm{d}x\colon \zeta-z\in W_0^{1,q}(B;\R^m),\ \zeta\equiv 0 \ \text{on } \overline{A}\bigg\}\,. 
\end{equation*}
As in the deterministic setting \cite{Ansini-Braides}, the $g$-capacity of the ball (with respect to $\R^n$) appears in the definition of the limit functional. For this reason  
we consider 
\begin{equation}\label{nonlinear_capacity_density}
\hspace{-0.8em}\varphi_{\rho}(z):=\mathrm{Cap}_g({B_\rho}(0),\R^n;z)=\inf\bigg\{\int_{\R^n}g(\nabla\zeta)\,\mathrm{d}x\colon \zeta-z\in W_0^{1,q}(\R^n;\R^m),\ \zeta\equiv 0 \ \text{on } \overline{B}_\rho(0)\bigg\}\,,
\end{equation}
which is well-defined for every $z\in \R^m$ and every $\rho>0$. Throughout the paper the function $\varphi_{\rho}$ is referred to as the \textit{nonlinear $g$-capacity density}.

We note that if $\rho \in \mathcal{R}$, namely $\rho:\Omega\to \R_+$ is a random variable, then $\omega\mapsto \varphi_{\rho^\omega}(z)$ is also a random variable. Indeed, for every $z\in \R^{m}$ the function $\rho\mapsto \varphi_{\rho}(z)$ is continuous (see Lemma \ref{growth_of_phi_rho} below and the estimates in \eqref{continuity_of_phi_rho_capacity}), hence the composite function $\omega\mapsto \varphi_{\rho^\omega}(z)$ is $\mathcal T$-measurable. We can then define the \textit{average $g$-capacity density} at $z\in \R^m$ as the expected value of $\varphi_{\rho^\omega}(z)$. That is, we set
\begin{equation}\label{expected_g_capacity}
\varphi(z):=\langle\varphi_{\rho}(z)\rangle:=\int_{0}^{+\infty}\varphi_\rho(z)h(\rho)\,\mathrm{d}\rho\,,
\end{equation}
where $h\in L^1(\R_+;\R_+)$ satisfies \eqref{finite_average_capacity}.

\medskip

We are now in a position to state the main result of this paper. 

\begin{theorem}\label{main_thm}
Let $(\Phi, \mathcal R)$ be a m.p.p.\ satisfying \ref{H1}--\ref{H4} and let $(\F_\eps^\o)_{\eps>0}$ be the functionals defined in \eqref{random_functionals}, for $\omega\in \Omega$.
Then there exists a sequence $(\eps_j) \searrow 0$ and a set $\Omega' \in \mathcal T$ with $\mathbb P(\Omega')=1$, such that for every $\omega\in \Omega'$
\begin{equation*}
\mathcal{F}_{\eps_j}^\omega {\overset{\Gamma}{\longrightarrow}} \, \mathcal{F}_0 \ \text{with respect to the strong } L^{1}(D;\R^m)\text{-topology}\,,
\end{equation*}
where $\F_0:L^{1}(D;\R^m)\longrightarrow \R\cup \{+\infty\}$ is the deterministic functional defined as
\begin{equation}\label{deterministic_functionals}
\F_0(u):=\begin{cases} \ds \int_{D}Qf(\nabla u)\,\mathrm{d}x+\langle N(Q)\rangle\int_D\varphi (u)\,\mathrm{d}x & \text{if } u\in W_0^{1,q}(D;\R^m)\,,
\cr
+ \infty & \text{otherwise in } L^{1}(D;\R^m)\,,
\end{cases}
\end{equation} 
with $N(Q)$ and $\varphi$ as in \eqref{second_moment_condition} and \eqref{expected_g_capacity}, respectively.
\end{theorem}

In general the $\Gamma$-convergence result in Theorem \ref{main_thm} holds true only up to subsequences. In fact the limit density $g$ appearing in the definition of $\varphi$ clearly depends on the choice of the subsequence (cf.\ \eqref{j_capacity_densities} and  \eqref{limit_capacity_density}). This phenomenon is typical of the nonlinear setting and is also observed in \cite{Ansini-Braides} (see Remark 2.7 therein).

\begin{remark} [Convergence of the Euler-Lagrange equations] \normalfont 
Let $(\e_j)$ be the vanishing sequence and $\Omega' \in \mathcal T$ be the set of probability one whose existence is established by Theorem \ref{main_thm}. 
Let $\psi\in W^{-1,\frac{q}{q-1}}(D; \R^m)$ be fixed; we consider the functionals defined for $u\in L^1(D; \R^m)$ as
\begin{equation}\label{c:forcing-term}
\mathcal F^\omega_{\eps_j}(u)+ \int_{D}\psi:u\,\mathrm{d}x\,,
\end{equation}
where $\int_{D}\psi:u\,\mathrm{d}x$ denotes the duality pairing between $W^{1,q}_0(D;\R^m)$ and  $W^{-1,\frac{q}{q-1}}(D;\R^m)$. By continuity it is immediate to check that for every $\omega \in \Omega'$ the functionals in \eqref{c:forcing-term} $\Gamma$-converge, with respect to the strong $L^{1}(D;\R^m)$-topology, to the deterministic functional defined for $u\in L^1(D; \R^m)$ as
\begin{equation*}
\mathcal F_0(u)+\int_D \psi:u\,\mathrm{d}x\,.
\end{equation*}
Let now $f$ be convex and differentiable; by the $L^1(D;\R^m)$ equi-coerciveness of \eqref{c:forcing-term} and the fundamental property of $\Gamma$-convergence we can deduce a convergence result for the corresponding Euler-Lagrange equations. That is, for every $\omega\in \Omega'$, the sequence $(u^\omega_{\eps_j})$, where $u^\omega_{\eps_j} \in W^{1,q}_0(D;\R^m)$ is the unique solution to 
\[
\begin{cases}
-\mathrm{div}Df(\nabla u^\omega_{\eps_j})=\psi & \text{in } D^\omega_{\eps_j} \cr
u^\omega_{\eps_j}=0  & \text{on } \partial D^\omega_{\eps_j}\,,
\end{cases}
\]
converges weakly in $W^{1,q}(D;\R^m)$ to the unique solution $u_0 \in W_0^{1,q}(D;\R^m)$ of the following deterministic problem 
\begin{equation*}
\begin{cases}
-\mathrm{div}Df(\nabla u_0)+\langle N(Q)\rangle \varphi'(u_0)=\psi & \text{in } D \cr
u_0=0 & \text{on } \partial D\,.
\end{cases}
\end{equation*}
Finally, we observe that when $f$ is $q$-homogeneous, the convergence in \eqref{limit_capacity_density} defining $g$ holds true for the whole sequence. Hence the function $g$ (and, consequently, $\varphi$ in \eqref{expected_g_capacity}) is independent of the subsequence. Therefore both the convergence of the functionals in Theorem \ref{main_thm} and the convergence of the optimality conditions above hold true for the whole sequence. Our result then provides an extension of \cite[Theorem 2.1]{Giunti-Hofer-Velasquez} to the nonlinear $q$-homogeneous vectorial setting.

\end{remark}

\section{Deterministic building blocks}\label{deterministic_preliminaries}

In this section we collect some (deterministic) technical results which will be used in the proof of Theorem \ref{main_thm}. %

We start by establishing some properties of the nonlinear capacity density $\varphi_{\rho}$ defined in \eqref{nonlinear_capacity_density}. Moreover, we introduce some auxiliary capacity densities whose role will become apparent in the next sections. 
Finally, we state and prove a so-called \textit{joining lemma}, Lemma \ref{joining_lemma_on_randomly_perforated_domains}, in the same spirit of \cite[Lemma 3.1]{Ansini-Braides}.

\begin{lemma}\label{growth_of_phi_rho} There exist two constants $C_1:=C_1(n,q,c_1), C_2:=C_2(n,q,c_2)>0$ such that 
\begin{equation}\label{bounds_on_phi_rho_capacity}
C_1|z|^q\rho^{n-q}\leq \varphi_{\rho}(z)\leq C_2|z|^q\rho^{n-q}\,, 
\end{equation}
for every $z\in \R^m$ and every $\rho>0$, where $c_1, c_2>0$ are the constants in \eqref{growth_of_f}. Moreover, for every $0<\rho_1<\rho_2$ and every $z\in \R^m$ we have that 
\begin{equation}\label{continuity_of_phi_rho_capacity}
\varphi_{\rho_{1}}(z)\leq \varphi_{\rho_{2}}(z) \leq \varphi_{\rho_{1}}(z)\Big(\frac{\rho_2}{\rho_1}\Big)^n\left(1+\frac{L}{c_1\rho_2^q}\Big(\rho_1^{q-1}+\rho_2^{q-1}\Big)(\rho_2-\rho_1)\right),
\end{equation}
where $L$ is the constant in \eqref{lipschitzianity_of_g}.

\end{lemma}
\begin{proof}
Let $z\in \R^m$ and $\rho>0$ be fixed; let
\begin{equation}\label{model_capacity}
\widehat \varphi_{\rho}(z):= \inf\Big\{\int_{\R^n}|\nabla\zeta|^q\,\mathrm{d}x\colon \zeta-z\in W_0^{1,q}(\R^n;\R^m),\ \zeta\equiv 0 \ \text{on } \overline{B}_\rho(0)\Big\}\,
\end{equation}
be the nonlinear capacity density corresponding to the $q$-Dirichlet energy; namely $\widehat \varphi_{\rho}$ is the density defined as \eqref{nonlinear_capacity_density} in the model case $f(\xi)=g(\xi)=|\xi|^q$. The unique solution $\widehat\zeta_{z,\rho}$ of \eqref{model_capacity} satisfies the following $q$-Laplace boundary value problem
\begin{equation*}
\begin{cases}
\mathrm{div}(|\nabla \widehat\zeta_{z,\rho}|^{q-2}\nabla \widehat\zeta_{z,\rho})=0 \quad \text{in}\; \R^{n}\setminus \overline{B}_{\rho}(0)\,,
\cr\cr
\widehat\zeta_{z,\rho}-z \in  W_0^{1,q}(\R^n;\R^m), \quad  \widehat\zeta_{z,\rho}|_{\overline{B}_{\rho}(0)}\equiv 0\,. 
\end{cases}
\end{equation*}
This can be computed explicitly and is given by the radially symmetric function
\begin{equation*}
\widehat\zeta_{z,\rho}(x)=z\big(-(\rho/|x|)^{(n-q)/(q-1)}+1\big)\chi_{\R^n\setminus \overline{B}_\rho(0)}.
\end{equation*}
A direct calculation yields
\begin{equation}\label{explicit_formula_for_model_capacity}
\widehat\varphi_\rho(z)=C_{n,q}|z|^q\rho^{n-q}\,,
\end{equation}
where
\[C_{n,q}:=\Big(\frac{n-q}{q-1}\Big)^{q-1}\mathcal{H}^{n-1}(\S^{n-1})\,.
\]
Therefore, gathering the upper bound in \eqref{growth_of_g}, \eqref{nonlinear_capacity_density}, and \eqref{explicit_formula_for_model_capacity}, we get
\begin{equation*}
\varphi_{\rho}(z)\leq \int_{\R^n}g(\nabla\widehat \zeta_{z,\rho})\,\mathrm{d}x\leq c_2\widehat\varphi_\rho(z)= c_2C_{n,q}|z|^q\rho^{n-q}\,,
\end{equation*}
which gives the second inequality in \eqref{bounds_on_phi_rho_capacity} with $C_2:=c_2C_{n,q}$. 
	
Similarly, for every $\zeta$ such that $\zeta-z\in W_0^{1,q}(\R^n;\R^m)$ and  $\zeta\equiv 0 \ \text{on } \overline{B_\rho}(0)$, again using  \eqref{growth_of_g}, \eqref{nonlinear_capacity_density}, and \eqref{explicit_formula_for_model_capacity}, we obtain
\begin{equation*}
\int_{\R^n}g(\nabla \zeta)\,\mathrm{d}x\geq c_1\int_{\R^n}|\nabla \zeta|^q\,\mathrm{d}x\geq c_1 \widehat\varphi_\rho(z)= c_1C_{n,q}|z|^q\rho^{n-q}\,.
\end{equation*}
Taking the infimum in $\zeta$ we get the first inequality in \eqref{bounds_on_phi_rho_capacity} with $C_1:=c_1C_{n,q}$. This concludes the proof of \eqref{bounds_on_phi_rho_capacity}.

For \eqref{continuity_of_phi_rho_capacity}, note that the first inequality follows by monotonicity, since for  $0<\rho_1<\rho_2 $ every competitor $\zeta$ for the minimisation problem defining $\varphi_{\rho_{2}}(z)$ is also a competitor for the minimisation problem defining $\varphi_{\rho_{1}}(z)$. To prove the second inequality, let $\zeta_1$ be the minimiser of the problem defining $\varphi_{\rho_{1}}(z)$. Set $\zeta_2(x):=\zeta_1\left(\rho_1x/\rho_2\right)$; then $\zeta_2$ is a competitor for the minimisation problem defining $\varphi_{\rho_2}(z)$. Moreover, by \eqref{lipschitzianity_of_g} and \eqref{growth_of_g} we get
\begin{align*}
g(\nabla \zeta_2 (x)) &\leq g(\nabla \zeta_1(\rho_1x/\rho_2)) + L\big(|\nabla \zeta_2(x)|^{q-1}+|\nabla \zeta_1(\rho_1x/\rho_2)|^{q-1}\big)|\nabla\zeta_2(x)-\nabla \zeta_1(\rho_1x/\rho_2)|\\
&\leq g(\nabla \zeta_1(\rho_1x/\rho_2)) + L\left(1+\Big(\frac{\rho_1}{\rho_2}\Big)^{q-1}\right)\left(1-\frac{\rho_1}{\rho_2}\right)|\nabla \zeta_1(\rho_1x/\rho_2)|^q\\
&\leq g(\nabla \zeta_1(\rho_1x/\rho_2)) \left(1+\frac{L}{c_1\rho_2^q}\Big(\rho_1^{q-1}+\rho_2^{q-1}\Big)(\rho_2-\rho_1)\right).
\end{align*}
By integrating over $\R^n$, a change of variables and the definition of $\zeta_1$ give
$$
\int_{\R^n}g(\nabla\zeta_2)\,\mathrm{d}x \leq \varphi_{\rho_{1}}(z)\Big(\frac{\rho_2}{\rho_1}\Big)^n\left(1+\frac{L}{c_1\rho_2^q}\Big(\rho_1^{q-1}+\rho_2^{q-1}\Big)(\rho_2-\rho_1)\right).
$$
Eventually, since $\zeta_2$ is a competitor for the minimization problem defining $\varphi_{\rho_{2}}(z)$, we obtain 
$$
\varphi_{\rho_{2}}(z) \leq \varphi_{\rho_{1}}(z)\Big(\frac{\rho_2}{\rho_1}\Big)^n\left(1+\frac{L}{c_1\rho_2^q}\Big(\rho_1^{q-1}+\rho_2^{q-1}\Big)(\rho_2-\rho_1)\right)
$$
and thus \eqref{continuity_of_phi_rho_capacity}.
\end{proof}

\smallskip

\begin{remark}
{\rm As an immediate corollary of Lemma \ref{growth_of_phi_rho} (cf.\ the estimates in  \eqref{bounds_on_phi_rho_capacity}) we have that $\varphi_\rho(u)\in L^1(D)$ whenever $u\in L^{q}(D;\R^m)$ and $0<\rho<+\infty$. Moreover, again by Lemma \ref{growth_of_phi_rho} and by \eqref{finite_average_capacity} and \eqref{expected_g_capacity} we have also that $\varphi(u)\in L^1(D)$ whenever $u\in L^{q}(D;\R^m)$, where $\varphi$ is defined in \eqref{expected_g_capacity}.}
\end{remark}

Let $\theta\in (0,1)$ and $\rho>0$ be fixed. Let $(\eps_j)\searrow 0$; we set
\begin{equation}\label{K_j_radii}
 K_j:=\frac{\eps_j}{\alpha_{\eps_j}}=\eps_j^{-q/n-q}\,,
\end{equation}
where $\alpha_{\e_j}$ is defined as in \eqref{eps_scaling}, and assume that $j\in \N$ is large enough to guarantee that
\begin{equation}\label{large_K_j}
\theta K_j \geq 2\rho\,.
\end{equation}
Moreover, for $z\in \R^m$ we define the class of functions  
\begin{equation*}
X^j_{\theta,\rho,z}:=\Big\{\zeta\colon \zeta-z\in W_0^{1,q}(B_{\theta K_{j}}(0);\R^m),\ \zeta\equiv 0 \ \text{on } \overline{B_\rho}(0)\Big\}\,,
\end{equation*}
and the auxiliary capacity densities 
\begin{equation}\label{K_j_rho_truncation_of_capacity}
\varphi^j_{\theta,\rho}(z):=\inf\bigg\{\int_{B_{\theta  K_j}(0)}g_j(\nabla\zeta)\,\mathrm{d}x\colon \zeta\in X^j_{\theta,\rho,z}\bigg\}\,.
\end{equation}
We observe that the function $\rho\mapsto \varphi^j_{\theta,\rho}(z)$ is increasing. Indeed, if $0<\rho_1\leq \rho_2$, then every competitor for the minimisation problem defining $\varphi^j_{\theta,\rho_2}(z)$ is also a competitor for the minimisation problem defining $\varphi^j_{\theta,\rho_1}(z)$.

\medskip

The next lemma is the analogue of Lemma \ref{growth_of_phi_rho} for the functions $\varphi_{\theta,\rho}^j$.

\begin{lemma}\label{bounds_for_the_truncated_capacities} 
Let $(\eps_j) \searrow 0$, $\theta\in (0,1)$  and $\rho>0$ be fixed, and let $K_j$ be defined as in \eqref{K_j_radii}. Assume that \eqref{large_K_j} is satisfied. Let $\varphi^j_{\theta,\rho}$ be as in \eqref{K_j_rho_truncation_of_capacity}; then
\begin{itemize}
\item[(i)] there exist constants $C_1,\dots,C_4>0$ depending only on $n,q,c_1,c_2$, such that for every $z\in \R^m$
\begin{equation}\label{lower_bound_for_the_truncated_capacities}
\varphi^j_{\theta,\rho}(z)\geq C_1|z|^q\big(\rho^{(q-n)/(q-1)}-(\theta K_j)^{(q-n)/(q-1)}\big)^{1-q}-C_2\theta^n\,
\end{equation}
and
\begin{equation}\label{upper_bound_for_the_truncated_capacities} 
\varphi^j_{\theta,\rho}(z)\leq C_3|z|^q\big(\rho^{(q-n)/(q-1)}-(\theta  K_j)^{(q-n)/(q-1)}\big)^{1-q}+C_4\theta^n.
\end{equation}

\item[(ii)]  Let $M>0$, and assume in addition that $\rho\in (0,M]$, and that $j\in \N$ is large enough so that
\begin{equation}\label{larger_j}
\theta  K_{j}\geq 2 M\,. 
\end{equation}
Then, there exists a constant $C_M>0$, with $C_M\to +\infty$ as $M\to +\infty$, so that 
\begin{equation}\label{Lipschitz_estimate_for_phi_j_theta_rho}
|\varphi^j_{\theta,\rho}(z)-\varphi^j_{\theta,\rho}(w)|\leq C_M\big(\theta^{n(q-1)/q}+\alpha_{\eps_j}^{q-1}+|z|^{q-1}+|w|^{q-1}\big)|z-w|\,, 
\end{equation}
for every $z,w\in \R^m$\,.
\end{itemize} 
\end{lemma}
\begin{proof} 
We start proving $(i)$. To this end, for every $z\in \R^m$ we set 
\begin{equation}\label{K_model_capacity}
\widehat \varphi^j_{\theta,\rho}(z):= \inf\bigg\{\int_{B_{\theta K_j}(0)}|\nabla\zeta|^q\,\mathrm{d}x\colon  \zeta\in X^j_{\theta,\rho,z}\bigg\}\,.
\end{equation} 
Similarly as in the proof of \eqref{bounds_on_phi_rho_capacity}, the unique solution $\widehat\zeta^j_{\theta,\rho,z}$ to \eqref{K_model_capacity} can be computed explicitly. Moreover, a direct calculation gives
\begin{equation}\label{explicit_K_model_capacity}
\widehat\varphi^j_{\theta,\rho}(z)=c_{n,q}|z|^q\big(\rho^{(q-n)/(q-1)}-(\theta K_j)^{(q-n)/(q-1)}\big)^{1-q}\,,
\end{equation}
where the constant $c_{n,q}>0$ can be computed explicitly. 

To prove \eqref{lower_bound_for_the_truncated_capacities}, let $ \zeta\in X^j_{\theta,\rho,z}$; note that $\zeta$ is a competitor for the minimisation problems defining $\varphi^j_{\theta,\rho}(z)$ and $\widehat \varphi^j_{\theta,\rho}(z)$. By combining \eqref{growth_of_g_j} and \eqref{explicit_K_model_capacity}, also recalling \eqref{eps_scaling}, \eqref{K_j_radii}, we find  
\begin{align*}
\int_{B_{\theta K_{j}}(0)}g_j(\nabla \zeta)\,\mathrm{d}x&\geq c_1\bigg(\int_{B_{\theta K_{j}}(0)}|\nabla \zeta|^q\,\mathrm{d}x-\beta_n(\theta K_{j})^n\alpha_{\eps_j}^q
\bigg)
\geq c_1\big(\widehat \varphi^j_{\theta,\rho}(z)-\beta_n\theta^n\big)\\
&= c_1 \big(c_{n,q}|z|^q\big(\rho^{(q-n)/(q-1)}-(\theta K_j)^{(q-n)/(q-1)}\big)^{1-q}-\beta_n\theta^n\big)\,.
\end{align*}
Then, by taking the infimum over $\zeta\in X^j_{\theta,\rho,z}$ we obtain \eqref{lower_bound_for_the_truncated_capacities} with $C_1:=c_1c_{n,q}>0$ and $C_2:=c_1\beta_n>0$.
 	
To prove  \eqref{upper_bound_for_the_truncated_capacities} we use the fact that $\widehat\zeta^j_{\theta,\rho,z}$ is a competitor for the minimisation problem defining $\varphi^j_{\theta,\rho}(z)$, which together with \eqref{growth_of_g_j} and  \eqref{explicit_K_model_capacity} gives
\begin{align*}
\varphi^j_{\theta,\rho}(z)&\leq \int_{B_{\theta K_{j}}(0)} g_j(\nabla \widehat\zeta^j_{\theta,\rho,z})\,\mathrm{d}x
\leq c_2\bigg(\int_{B_{\theta K_{j}}(0)}|\nabla \widehat\zeta^j_{\theta,\rho,z}|^q\,\mathrm{d}x+\beta_n(\theta K_{j})^n \alpha_{\eps_j}^q\bigg)\\
&= c_2\big(\widehat\varphi^j_{\theta,\rho}(z)+\beta_n \theta^n\big)
= c_2 \big(c_{n,q}|z|^q\big(\rho^{(q-n)/(q-1)}-(\theta K_{j})^{(q-n)/(q-1)}\big)^{1-q}+\beta_n\theta^n
\big)\,, 
\end{align*}
which proves \eqref{upper_bound_for_the_truncated_capacities} with $C_3:=c_2c_{n,q}>0$ and $C_4:=c_2\beta_n>0$. 

We now prove $(ii)$. By the quasiconvexity and coercivity of $g_j$ (cf. \eqref{j_capacity_densities} and \eqref{growth_of_g_j}), the Direct Methods of the Calculus of Variations ensure the existence of a function $\zeta^j_{\theta,\rho,z}\in X^j_{\theta,\rho,z}$ such that 
\begin{equation}\label{attainment_of_minimizer_for_phi_K_j_rho}
\varphi^j_{\theta,\rho}(z)=\int_{B_{\theta K_{j}}(0)}g_j(\nabla \zeta^j_{\theta,\rho,z})\,\mathrm{d}x\,. 
\end{equation}
We now modify $\zeta^j_{\theta,\rho,z}$ to construct a competitor for the minimization problem defining $\varphi^j_{\theta,\rho}(w)$. To this end, we consider a radial cut-off function $\eta_{M}\in C_c^{\infty}(B_{2M}(0))$ satisfying
\begin{equation}\label{cut_off_in_B_K}
 0\leq \eta_{M}\leq 1\,, \ \ \eta_{M}|_{\overline{B}_{\rho}(0)}\equiv 1\,, \ \ \|\nabla \eta_{M}\|_{L^\infty}\leq \frac{1}{2M-\rho}\leq\frac{1}{M} \,, 
\end{equation}
where in the last inequality we used that $\rho\in (0, M]$.  We now define the function $\zeta^j_{\theta,\rho,w}$ as  
\begin{equation*} 
\zeta^j_{\theta,\rho,w}:=\zeta^j_{\theta,\rho,z}+(1-\eta_{M})(w-z)\,.
\end{equation*}
In view of \eqref{cut_off_in_B_K} and \eqref{larger_j}, for $j\in \N$ large enough we have that $\eta_M|_{\partial B_{\theta K_{j}}}\equiv 0$, and therefore $\zeta^j_{\theta,\rho,w}\in X^j_{\theta,\rho,w}$, thus, it is an admissible competitor for the minimisation problem defining $\varphi^j_{\theta,\rho}(w)$. Note that
\begin{equation}\label{grad_zw}
\nabla \zeta^j_{\theta,\rho,w} -\nabla \zeta^j_{\theta,\rho,z} =\nabla \eta_{M}\otimes(z-w)\in C^\infty_c(B_{2M}(0);\R^{m\times n})\,.
\end{equation}
By \eqref{K_j_rho_truncation_of_capacity}, \eqref{attainment_of_minimizer_for_phi_K_j_rho},  \eqref{equilipschitzianity}, \eqref{growth_of_g_j}, \eqref{larger_j} \eqref{cut_off_in_B_K} and \eqref{grad_zw}, we obtain 
\begin{align}\label{auxiliary_lipschitz_estimate}
 \varphi^j_{\theta,\rho}(w)-\varphi^j_{\theta,\rho}(z)
&\leq \int_{B_{\theta K_{j}}(0)}\big|g_j(\nabla \zeta^j_{\theta,\rho,w})-g_j(\nabla \zeta^j_{\theta,\rho,z})\big|\,\mathrm{d}x \nonumber\\
&\leq L \int_{B_{2M}(0)}\Big(\alpha_{\eps_j}^{q-1}+|\nabla \zeta^j_{\theta,\rho,w}|^{q-1}+|\nabla \zeta^j_{\theta,\rho,z}|^{q-1}\Big)|\nabla  \zeta^j_{\theta,\rho,w}-\nabla \zeta^j_{\theta,\rho,z}|\,\mathrm{d}x \nonumber\\
&\leq \frac{L}{M}\bigg(\beta_n(2M)^n\alpha_{\eps_j}^{q-1}+\int_{B_{2M}(0)}\Big(|\nabla \zeta_{\theta,\rho,w}^j|^{q-1}+|\nabla\zeta^j_{\theta,\rho,z}|^{q-1}\Big)\,\mathrm{d}x\bigg)|w-z|\nonumber\\
&\leq \frac{C}{M}\bigg(M^n\alpha_{\eps_j}^{q-1}
+\frac{M^{n}}{M^{q-1}}|z-w|^{q-1}+M^{n/q}\Big(\int_{B_{2M}(0)}|\nabla\zeta^j_{\theta,\rho,z}|^{q}\Big)^{\tfrac{q-1}{q}}\bigg)|w-z|\nonumber\\
&\leq \frac{CM^n}{M}\left(\alpha_{\eps_j}^{q-1}+|z|^{q-1}+|w|^{q-1}+\Big(\varphi^j_{\theta,\rho}(z)+M^n\alpha_{\eps_j}^{q}
\Big)^{\tfrac{q-1}{q}}\right)|w-z|\,,
\end{align}
where $C$ is a positive constant depending only on $c_1,c_2,n,q, L$. Now, by \eqref{upper_bound_for_the_truncated_capacities} 
\begin{align}\label{c:c-star}
\varphi^j_{\theta,\rho}(z) & \leq C_3|z|^q\big(\rho^{(q-n)/(q-1)}-(\theta  K_j)^{(q-n)/(q-1)}\big)^{1-q}+C_4\theta^n\nonumber \\
& \leq C_3|z|^q\big(1-2^{(q-n)/(q-1)}\big)^{1-q}\rho^{n-q} +C_4\theta^n\,,
\end{align}
where the last inequality follows from \eqref{larger_j}, for $j\in \N$ large enough (depending on $\theta,M$). 
Finally, \eqref{auxiliary_lipschitz_estimate} and \eqref{c:c-star} imply that 
\begin{align*}
\varphi^j_{\theta,\rho}(w)-\varphi^j_{\theta,\rho}(z)&\leq
 CM^{n-1}\left(\alpha_{\eps_j}^{q-1}+|z|^{q-1}+|w|^{q-1}+\Big(|z|^q\rho^{n-q}+\theta^n+M^n\alpha_{\eps_j}^{q}\Big)^{\tfrac{q-1}{q}}\right)|w-z|\,\\
 & \leq CM^{2n-1} \left(\theta^{{n(q-1)/q}}+\alpha_{\eps_j}^{q-1}+|z|^{q-1}+|w|^{q-1}\right)|w-z|\,.
\end{align*}
Then, interchanging the role of $z$ and $w$ in the previous argument we obtain \eqref{Lipschitz_estimate_for_phi_j_theta_rho}, and this concludes the proof.
\end{proof}

Note that, as detailed in the remark below, sequences of minimisers of $(\varphi^j_{\theta, \rho})_j$ are pre-compact, modulo a straightforward extension.

\begin{remark}[Compactness after extension]\label{compactness_after_modification}
{\rm
Let $\theta\in (0,1)$, $M>0$, $\rho\in(0,M]$, $z\in \R^m$ be fixed, and let $j\in \N$ be so large that \eqref{larger_j} holds. Let $(\zeta^j)\subset X^j_{\theta,\rho,z}$ be such that
\begin{equation}\label{equibounded_j_energy}
\sup_{j \in \N} \int_{B_{\theta K_{j}}(0)}g_j(\nabla \zeta^j)\,\mathrm{d}x=: C_{\theta,\rho,z}<+\infty\,,	
\end{equation}
where $g_j$ is as in \eqref{j_capacity_densities}. We set 
\begin{equation*}
\tilde \zeta^j:=\begin{cases}
\zeta^j \ \ \text{in } B_{\theta K_{j}}(0)\,,\\
z \ \   \text{in } \R^n\setminus B_{\theta K_{j}}(0)\,.
\end{cases}	
\end{equation*}
Clearly $(\tilde \zeta^j)\subset W^{1,q}_{\mathrm{loc}}(\R^n;\R^m)$, and $\tilde\zeta^j-z\in W_0^{1,q}(\R^n;\R^m)$. Moreover, in view of \eqref{growth_of_g_j}, \eqref{eps_scaling}, \eqref{K_j_radii}, and \eqref{equibounded_j_energy} we have
\begin{align}\nonumber
\sup_{j \in \N}\int_{\R^n}|\nabla (\tilde \zeta^j-z)|^q\,\mathrm{d}x&=
\sup_{j \in \N}\int_{\R^n}|\nabla \tilde \zeta^j|^q\,\mathrm{d}x= \sup_{j \in \N}\int_{B_{\theta K_{j}(0)}}|\nabla \zeta^j|^q\,\mathrm{d}x\\\nonumber
&\lesssim \sup_{j \in \N}\Big(\int_{B_{\theta K_{j}(0)}}g_j(\nabla \zeta^j)\,\mathrm{d}x+(\theta K_{j})^n\alpha_{\eps_j}^q\Big)\\\label{tilde_zeta_j_equibounded_in_W_1_q}
&\lesssim C_{\theta,\rho,z}+\theta^n.
\end{align}
Then, by the Sobolev Embedding Theorem we deduce the existence of $c_{n,q}>0$ such that
\begin{equation}\label{Sobolev_inequality}
\sup_{j \in \N} \int_{\R^n}|\tilde \zeta^j-z|^{q^*}\,\mathrm{d}x\leq c_{n,q} \sup_{j \in \N}\int_{\R^n}|\nabla \tilde \zeta^j|^q\,\mathrm{d}x\lesssim C_{\theta,\rho,z}+\theta^n,
\end{equation}
where $q^*:=nq/n-q$ denotes the conjugate Sobolev exponent of $q$. 
The estimates \eqref{tilde_zeta_j_equibounded_in_W_1_q} and \eqref{Sobolev_inequality} then guarantee that 
\begin{equation*}
\tilde \zeta^j-z\rightharpoonup \tilde \zeta \; \text{ weakly in }\; W^{1,q}_{\mathrm{loc}}(\R^n;\R^m)\ \text{as } j\to +\infty\,, 
\end{equation*}
for some $\tilde \zeta\in W^{1,q}_{\mathrm{loc}}(\R^n;\R^m)$. Equivalently, setting $\zeta:=\tilde \zeta+z$, we have that $\zeta-z\in W_0^{1,q}(\R^n;\R^m)$, $\zeta|_{\overline{B}_\rho(0)}\equiv 0$, and 
\[
\tilde\zeta^j\rightharpoonup \zeta\ \text{weakly in } W^{1,q}_{\mathrm{loc}}(\R^n;\R^m)\ \text{as } j\to +\infty\,.
\]
}
\end{remark}

The following result is an immediate consequence of Lemma \ref{bounds_for_the_truncated_capacities}.

\begin{corollary}\label{local_uniform_convergence_for_capacity_densities}
Let $M>0$ and $\rho\in (0,M]$ be fixed. Let $\theta\in (0,1)$ and $j \in \N$ satisfy \eqref{large_K_j} and let $\varphi^j_{\theta,\rho}$ be as in \eqref{K_j_rho_truncation_of_capacity}. Then we have the following.
\begin{itemize}
\item [(i)]  There exists a subsequence (not relabelled) and a measurable function $\varphi_{\theta,\rho}\colon \R^m \to \R$ such that
\begin{equation}\label{local_uniform_convergence_phi_j}
\varphi^j_{\theta,\rho}\longrightarrow \varphi_{\theta,\rho} \; \text{ in }\; L^\infty_{\rm{loc}}(\R^m),\; \text{as } j\to+\infty\,. 
\end{equation}
Moreover, the function $\varphi_{\theta,\rho}$ satisfies
\begin{equation}\label{phi_theta_rho__bounds} 
C_1|z|^q\rho^{n-q}-C_2\theta^n\leq \varphi_{\theta,\rho}(z)\leq C_3|z|^q\rho^{n-q}+C_4\theta^n\, 
\quad \, \forall \, z\in \R^m\,,
\end{equation}
and 
\begin{equation}\label{phi_theta_rho_Lipschitz_estimate}
|\varphi_{\theta,\rho}(z)-\varphi_{\theta,\rho}(w)|\leq C_M\big(\theta^{{n(q-1)/q}}+|z|^{q-1}+|w|^{q-1}\big)|z-w|\, \quad \forall \,z,w\in \R^m\,,
\end{equation} 
where $C_1, \ldots, C_4>0$ and $C_M>0$ are as in Lemma \ref{bounds_for_the_truncated_capacities}.  
Additionally, for $\theta\in (0,1)$ and $z\in \R^m$ fixed, the function  $\rho \mapsto \varphi_{\theta,\rho}(z)$ is increasing.

\smallskip

\item[(ii)] There exists a subsequence (not relabelled) and a measurable function $\tilde\varphi_{\rho}\colon \R^m \to \R$ such that
\begin{equation}\label{theta_local_uniform_convergence_phi_j}
\varphi_{\theta,\rho}\longrightarrow \tilde \varphi_{\rho} \; \text{ in }\; L^\infty_{\rm{loc}}(\R^m),\; \text{ as } \theta\to 0^+\,. 
\end{equation}
Moreover, the function $\tilde \varphi_{\rho}$ satisfies
\begin{equation}\label{tilde_phi_rho_bounds} 
C_1|z|^q\rho^{n-q}\leq \tilde\varphi_{\rho}(z)\leq C_3|z|^q\rho^{n-q}\, \quad \, \forall \, z\in \R^m\,,
\end{equation}
and 
\begin{equation}\label{tilde_phi_rho_Lipschitz_estimate}
|\tilde \varphi_{\rho}(z)-\tilde \varphi_{\rho}(w)|\leq C_M\big(|z|^{q-1}+|w|^{q-1}\big)|z-w|\, \quad \forall \,z,w\in \R^m\,,
\end{equation} 
where $C_1, C_3>0$ and $C_M>0$ are as in Lemma \ref{bounds_for_the_truncated_capacities}.  
Additionally, for $z\in \R^m$ fixed the function  $\rho \mapsto \tilde \varphi_{\rho}(z)$ is increasing. 
\end{itemize}
\end{corollary}

Our next goal is to prove that the abstract limit density $\tilde\varphi_\rho$ given by Corollary \ref{local_uniform_convergence_for_capacity_densities} ($ii$) coincides with $\varphi_{\rho}$ defined in \eqref{nonlinear_capacity_density}. To do so, we follow a similar approach as in \cite[Section 7]{Ansini-Babadjian-Zeppieri}.

\begin{proposition}\label{identification_of_capacity_density} Let $M>0$, and $\rho\in(0,M]$ be fixed. Let $\varphi_{\rho}$ be as in \eqref{nonlinear_capacity_density} and let  $\tilde\varphi_\rho$ be given by Corollary \ref{local_uniform_convergence_for_capacity_densities} $(ii)$. 
Then, $\tilde  \varphi_\rho\equiv\varphi_\rho$.
Therefore, in particular, $\varphi_{\rho}$ satisfies the estimates \eqref{tilde_phi_rho_bounds} and \eqref{tilde_phi_rho_Lipschitz_estimate}.
\end{proposition}

\begin{proof}
We divide the proof into two main steps.

\medskip

\emph{Step 1: $\tilde \varphi_\rho\geq \varphi_{\rho}$.} Let $z\in \R^m$, $\theta\in (0,1)$ and $K\geq 2M$ be fixed. Moreover, let $(\eps_j) \searrow 0$, and take $j\in \N$ so large that $\theta K_{j}> K$ holds, where $K_j$ is defined in \eqref{K_j_radii}.

In view of the quasiconvexity and coercivity of $g_j$ (cf.\ \eqref{j_capacity_densities} and \eqref{growth_of_g_j}), there exists $\zeta^j\in X^j_{\theta,\rho,z}$ so that 
\begin{equation*}
\varphi^j_{\theta,\rho}(z)=\int_{B_{\theta K_{j}(0)}}g_j(\nabla \zeta^j)\,\mathrm{d}x\,.	
\end{equation*}
Hence, up to a subsequence,  \eqref{local_uniform_convergence_phi_j} ensures that
$$
\sup_{j\in \N}\int_{B_{\theta K_{j}}(0)}g_j(\nabla \zeta^j)\,\mathrm{d}x\leq C_{\theta,\rho,z}\,,
$$
for some constant $C_{\theta,\rho,z}>0$. As observed in Remark \ref{compactness_after_modification} we can extend $(\zeta^j)$ to obtain a new sequence $(\tilde \zeta^j)$ and a function $\zeta \in W^{1,q}_{\rm{loc}}(\R^n;\R^m)$ such that $\tilde\zeta^j\rightharpoonup \zeta$ weakly in $W^{1,q}_{\mathrm{loc}}(\R^n;\R^m)$ as $j\to +\infty$.
Then, in particular using \eqref{growth_of_g_j}, \eqref{K_j_radii}, \eqref{eps_scaling}, and since
\[\int_{B_{\theta K_j}(0)\setminus B_K(0)}g_j(\nabla \tilde \zeta_j)\geq -c_1\alpha_{\eps_j}^q\L^n(B_{\theta K_j}(0)\setminus B_{K}(0))\geq -\beta_nc_1\theta^n\,,\]
we get
\begin{equation}\label{1_one_side_ineq}
\varphi_{\theta,\rho}(z)=\lim_{j\to+\infty}\varphi_{\theta,\rho}^j(z)=\lim_{j\to+\infty}\int_{B_{\theta K_{j}}(0)}g_j(\nabla \zeta^j)\,\mathrm{d}x\geq \varliminf_{j\to+\infty}\int_{B_{K}(0)}g_j(\nabla \tilde \zeta^j)\,\mathrm{d}x-\beta_nc_1\theta^n\,.
\end{equation}
Now consider the auxiliary integral functionals defined as 
\[
\mathcal{G}^j_K(\zeta):=\begin{cases}\ds\int_{B_K(0)} g_j(\nabla \zeta)\,\mathrm{d}x & \text{if } \zeta\in W^{1,q}(B_K(0);\R^m)\,,
\cr
+\infty &  \text{otherwise in } L^1(B_K(0);\R^m)\,,
\end{cases}
\]
and 
\[
\mathcal{G}_K(\zeta):=\begin{cases}\ds\int_{B_K(0)} g(\nabla \zeta)\,\mathrm{d}x & \text{if } \zeta\in W^{1,q}(B_K(0);\R^m)\,,
\cr
+\infty &  \text{otherwise  in } L^1(B_K(0);\R^m)\,,
\end{cases}
\]
where $g$ is as in \eqref{limit_capacity_density}.
Then, we can invoke \cite[Proposition 12.8]{Braides-De Franceschi} to deduce that the functionals $(\mathcal G^j_K)_j$ $\Gamma$-converge  to $\mathcal{G}_K$ with respect to the strong $L^{1}(B_K(0);\R^m)$-topology and the weak-$W^{1,q}(B_K(0);\R^m)$-topology. 
Since $\tilde \zeta^j\rightharpoonup \zeta$ weakly in $W^{1,q}(B_K(0);\R^m)$, in particular we can deduce that 
\begin{equation}\label{2_one_side_ineq}
\underset{j\to+\infty}{\liminf}\int_{B_{K}(0)}g_j(\nabla \tilde \zeta^j)\,\mathrm{d}x\geq \int_{B_K(0)}g(\nabla \zeta)\,\mathrm{d}x\,.
\end{equation}
Therefore, gathering \eqref{1_one_side_ineq} and \eqref{2_one_side_ineq} gives 
$$\varphi_{\theta,\rho}(z)\geq \int_{B_K(0)}g(\nabla \zeta)\,\mathrm{d}x-\beta_nc_1\theta^n\,.
$$
 Passing to the limit as $K\to +\infty$ and using the Dominated Convergence Theorem, we obtain
$$
\varphi_{\theta,\rho}(z)\geq \int_{\R^n}g(\nabla \zeta)\,\mathrm{d}x-\beta_nc_1\theta^n\geq \varphi_{\rho}(z)-\beta_nc_1\theta^n\,,
$$
where we also used that $\zeta$ is a competitor for the minimisation problem  defining $\varphi_{\rho}(z)$. Eventually, passing to the limit as $\theta\to 0^+$ and recalling \eqref{theta_local_uniform_convergence_phi_j}, we get that $\tilde \varphi_\rho(z)\geq \varphi_{\rho}(z)$.

\medskip

\emph{Step 2: $\tilde \varphi_\rho\leq \varphi_{\rho}$.} Let $z\in \R^m$ be fixed; by the quasiconvexity and coercivity of $g$ (cf.\ \eqref{limit_capacity_density} and \eqref{growth_of_g}), there exists $\zeta$ with $\zeta-z\in W_0^{1,q}(\R^n;\R^m)$, $\zeta|_{\overline{B}_\rho(0)}\equiv 0$, such that
\begin{equation}\label{reverse_attainment}
\varphi_{\rho}(z)=\int_{\R^n}g(\nabla \zeta)\,\mathrm{d}x\,.
\end{equation}
Let $\theta\in (0,1)$ and $K\geq 2M$ be fixed, and let $(\eps_j) \searrow 0$; we take $j\in \N$ so large that $\theta K_{j}> K$ holds. Consider a smooth radial cut-off function $\eta_K \in C_c^{\infty}(B_{K}(0))$ satisfying 
\begin{equation*}\label{eta_k_cut_off}
0\leq \eta_{K}\leq 1\,, \ \ \eta_{K}|_{\overline{B}_{K/2}(0)}\equiv 1\,, \ \ \|\nabla \eta_{K}\|_{L^\infty}\leq \frac{2}{K}\,, 	
\end{equation*}	
and define 
\begin{equation}\label{zeta_k_modification}
\zeta_{K}:=\eta_K\zeta+(1-\eta_K)z\,. 	
\end{equation}	
Then $\zeta_K$ is an admissible test function for the following minimisation problem
\begin{equation*}
\tilde \varphi^j_{K,\rho}(z):=\min\left\{\int_{B_K(0)}g_j(\nabla \tilde \zeta)\,\mathrm{d}x\colon \tilde \zeta-z\in W^{1,q}_0(B_K(0);\R^m)\,,\ \tilde \zeta|_{\overline{B}_{\rho}(0)}\equiv 0\right\}\,. 
\end{equation*}
By the $\Gamma$-convergence of $(\mathcal{G}_K^j)_j$ to $\mathcal{G}_K$, and by \cite[Proposition 11.7]{Braides-De Franceschi}, there exists a recovery sequence for $\mathcal{G}_K^j$ converging to $\zeta_K$. More precisely, there is $(\zeta_K^j)_j$, with $\zeta_K^j \rightharpoonup \zeta_K$ weakly in $W^{1,q}(B_K(0);\R^m)$,
\begin{equation*}
\zeta_K^j-z\in W_0^{1,q}(B_K(0);\R^m)\,,\ \ \zeta_K^j|_{\overline{B}_{\rho}(0)}\equiv 0\,,
\end{equation*}  
and for which the energies converge, namely
\begin{equation}\label{second_property_recovery_sequence}
\lim_{j\to +\infty}\int_{B_K(0)}g_j(\nabla \zeta_K^j)\,\mathrm{d}x=\int_{B_K(0)} g(\nabla \zeta_K)\,\mathrm{d}x\,.
\end{equation}
We extend $\zeta_K^j$ to $B_{\theta K_{j}}(0)$ by setting 
\begin{equation}\label{tilde_zeta_K_j}
\tilde \zeta_K^j:=\begin{cases}
\zeta_K^j & \text{in } B_K(0)
\cr
z &\text{in } B_{\theta K_{j}}(0)\setminus B_K(0)\,.
\end{cases}
\end{equation}
Then $\tilde \zeta_K^j\in X^j_{\theta,\rho,z}$, so $\tilde \zeta_K^j$ is an admissible competitor for the minimisation problem defining $\varphi^j_{\theta,\rho}(z)$ in \eqref{K_j_rho_truncation_of_capacity}. Therefore, by \eqref{tilde_zeta_K_j} we have
\begin{equation*}
\varphi^j_{\theta,\rho}(z)\leq \int_{B_{\theta K_{j}}(0)} g_j(\nabla \tilde \zeta_K^j)\,\mathrm{d}x= \int_{B_{K}(0)} g_j(\nabla\zeta_K^j)\,\mathrm{d}x\,.
\end{equation*}
Taking the limit as $j\to +\infty$ (possibly passing to a subsequence), and using \eqref{local_uniform_convergence_phi_j}, \eqref{second_property_recovery_sequence}, \eqref{growth_of_g} and \eqref{zeta_k_modification}, we obtain
\begin{align}\nonumber
\varphi_{\theta,\rho}(z)&\leq \int_{B_K(0)}g(\nabla \zeta_K)\,\mathrm{d}x=\int_{B_K(0)\setminus B_{K/2}(0)}g(\nabla \zeta_K)\,\mathrm{d}x+\int_{B_{K/2}(0)}g(\nabla \zeta_K)\,\mathrm{d}x\\\label{1_reverse_bound}
&\leq c_2\int_{B_K(0)\setminus B_{K/2}(0)}|\nabla \zeta_K|^q\,\mathrm{d}x+\int_{B_{K/2}(0)}g(\nabla \zeta)\, \mathrm{d}x\,.
\end{align}
We claim that 
\begin{equation}\label{remainder_goes_to_0_wrt_K}
\varlimsup_{K\to +\infty}\int_{B_K(0)\setminus B_{K/2}(0)}|\nabla \zeta_K|^q\,\mathrm{d}x=0\,.
\end{equation}
By \eqref{zeta_k_modification}, and since by the Sobolev Embedding Theorem $\zeta-z\in L^{q^*}(\R^n;\R^m)$, we have
\begin{align*}	
&\int_{B_K(0)\setminus B_{K/2}(0)}|\nabla \zeta_K|^q\,\mathrm{d}x= \int_{B_K(0)\setminus B_{K/2}(0)}\big|\nabla \eta_K\otimes (\zeta-z)+\eta_K\nabla \zeta\big|^q\,\mathrm{d}x\\
&\qquad\lesssim \|\nabla \eta_K\|^q_{L^\infty}\int_{B_K(0)\setminus B_{K/2}(0)}\big|\zeta-z\big|^q \,\mathrm{d}x+\| \eta_K\|^q_{L^\infty}\int_{B_K(0)\setminus B_{K/2}(0)}\big|\nabla \zeta|^q\,\mathrm{d}x\\
&\qquad\lesssim  K^{-q}K^{n(1-\frac{q}{q^*})}\Big(\int_{B_K(0)\setminus B_{K/2}(0)}|\zeta-z|^{q^*}\mathrm{d}x\Big)^{\frac{q}{q^*}}+\int_{B_K(0)\setminus B_{K/2}(0)}|\nabla \zeta|^q\,\mathrm{d}x\\[1pt]
&\qquad\lesssim  \Big(\int_{\R^n\setminus B_{K/2}(0)}|\zeta-z|^{q^*}\,\mathrm{d}x\Big)^{\frac{q}{q^*}}+\int_{\R^n\setminus B_{K/2}(0)}|\nabla \zeta|^q\,\mathrm{d}x\,,
\end{align*}
where we have used that $K^{-q}K^{n(1-\frac{q}{q^*})}=1$. Hence, passing to the limsup as $K\to +\infty$ yields \eqref{remainder_goes_to_0_wrt_K}. Taking the limsup as $K\to +\infty$ in \eqref{1_reverse_bound}, by \eqref{reverse_attainment}, \eqref{remainder_goes_to_0_wrt_K}, and the Dominated Convergence Theorem we get
\begin{equation*}
\varphi_{\theta,\rho}(z)\leq \varphi_\rho(z)\,.
\end{equation*}
Eventually, by taking the limit as $\theta\to 0^+$ in the above inequality, and by  \eqref{theta_local_uniform_convergence_phi_j}, we have that 
$\tilde \varphi_\rho(z)\leq \varphi_{\rho}(z)$.
\end{proof}
The following result, Lemma \ref{joining_lemma_on_randomly_perforated_domains}, is an adaptation to the non-periodic setting of the so-called ``Joining Lemma'' (cf.\  \cite[Lemma 3.1]{Ansini-Braides}). This is a technical tool which allows to modify sequences of functions near ``good perforations" without increasing the energy too much, and will be crucial in the proof of Theorem \ref{main_thm}. 
First, we need to define the class of the good perforations, consisting  of balls which are well-separated from one another and not too large, in the sense specified below.

\medskip

Let $\eps>0$ and  $M>0$ be fixed, and let $D\subset \R^n$ be an open, bounded, Lipschitz set, star-shaped with respect to the origin. We define $\mathscr{G}_{\eps,M}$ as the collection of points in $\R^n$  satisfying the following two properties:
\begin{enumerate}[label=(\alph*)]
\item \label{c:dis}  $|x_i-x_j|\geq 2/M$ for every $i\neq j$\,;

\smallskip

\item \label{c:card} $\bigcup_{x_i\in \mathscr{G}_{\eps,M}} \overline{B}_{\eps/M}(\eps x_i)\subset D$.
\end{enumerate}
We note that by \ref{c:dis} and \ref{c:card} the family of balls $({B}_{\eps/M}(\eps x_i))_{x_i\in \mathscr{G}_{\eps,M}}$ consists of pairwise disjoint subsets of $D$; therefore we immediately get that 
\begin{equation}\label{former:i}
(\beta_n M^{-n}\varepsilon^n) \# \mathscr{G}_{\eps,M} \leq \mathcal{L}^n(D)\,.
\end{equation}
For $\theta\in(0,1)$ fixed, we refer to $(B_{\theta\eps/M}(\eps x_i))_{x_i\in \mathscr{G}_{\eps,M}}$ as the family of ``truncated good perforations'' in $D$. By \ref{c:card} we have that $B_{\theta\eps/M}(\eps x_i) \subset \subset D$, for every
$x_i\in \mathscr{G}_{\eps,M}$.

For $x_i\in \mathscr{G}_{\eps,M}$ and $l\in \N$ we define the annulus
\begin{equation}\label{def: annuli}
C_{\eps,\theta, M}^{l}(\e x_i):=\Big\{x\in \R^n\colon 2^{-(l+1)}\theta\eps/M<|x-\eps x_{i}|< 2^{-l}\theta\eps/M\Big\} \subseteq B_{\theta\eps/M}(\eps x_i)\,.
\end{equation}
If $(\eps_j)\searrow 0$ we adopt the shorthand notation $\mathscr{G}_{j,M}:=\mathscr{G}_{\eps_j,M}$\,.

\medskip

We are now ready to state and prove the following variant of the Joining Lemma. 

\begin{lemma}[Joining Lemma]\label{joining_lemma_on_randomly_perforated_domains}
Let $(\eps_j) \searrow 0$, $M>0$, $\theta\in(0,1)$ and $k\in \N$ be fixed. 
Let $ (u_j)\subset W_0^{1,q}(D;\R^m)$ be such that 
\begin{equation}\label{u_j_to_u_weakly}
u_j\rightharpoonup u \ \text{ weakly in } W^{1,q}(D;\R^m)\,,
\end{equation}
for some $u \in W_0^{1,q}(D;\R^m)$, and let $\mathscr{G}_{j,M}$ be a collection of points in $\R^n$ satisfying \ref{c:dis}-\ref{c:card}. 
Then for every $x_{ j,i}\in \mathscr{G}_{j,M}$ there exists $k_{j,i}\in\{0,\dots,k-1\}$ and a corresponding annulus $C_{\eps_j,\theta, M}^{k_{j,i}}(\e_j x_{j,i})$ (defined as in \eqref{def: annuli} with $\e$, $l$, and $x_i$ replaced by $\e_j, k_{j,i}$, and $x_{j,i}$, respectively), such that 
we can construct a sequence $(w_j)\subset W_0^{1,q}(D;\R^m)$ satisfying the following properties:
\begin{enumerate}[label=(\roman*)]
\item \label{properties_of_the_modification_w_j_M-1} $w_j\equiv u_j \ \text{ in } D\setminus \bigcup_{x_{j,i} \in \mathscr{G}_{j,M}} C_{\eps_j,\theta, M}^{k_{j,i}}(\e_j x_{j,i})$\,;
\smallskip
\item \label{properties_of_the_modification_w_j_M-2} $w_j\equiv \bar u_{j,i} \ \text{ on } \partial B_{ \bar \sigma_{j,i}}(\eps_j x_{j,i})$, where 
\begin{equation}\label{intermediate_values}
\bar u_{j,i}:=\fint_{C_{\eps_j,\theta, M}^{k_{j,i}}(\e_j x_{j,i})} u_j\,\mathrm{d}x\,, \qquad  
\bar \sigma_{j,i}:=\frac{3}{4}2^{-k_{j,i}}\frac{\theta\eps_j}{M}\,;
\end{equation}
\item \label{w_j_M_to_u_weakly} $w_j\rightharpoonup u \ \text{ weakly in } W^{1,q}(D;\R^m)\,;$
\smallskip
\item \label{c:jl_estimate}$\displaystyle\Big|\int_{D} f(\nabla w_j)\,\mathrm{d}x-\int_{D}f(\nabla u_j)\,\mathrm{d}x\Big|\leq \frac{C}{k}\,$, for some $C>0$ depending on $c_2$, $n$, $m$, $q$, $D$, and $\sup_{j\in \N}\|\nabla u_j\|_{L^q(D;\R^{m\times n})}$.
\end{enumerate}

If, additionally, the sequence $(|\nabla u_j|^q)$ is equi-integrable, then also $(|\nabla w_j|^q)$ is equi-integrable, and one can take $k_{j,i} =0$ for all $x_{j,i} \in \mathscr{G}_{j,M}$, up to replacing \ref{c:jl_estimate} with the following estimate
\begin{equation}\label{modified_energy_estimate_equiiintegrability}
\Big|\int_{D} f(\nabla w_j)\,\mathrm{d}x-\int_{D}f(\nabla u_j)\,\mathrm{d}x\Big|\leq C_k\theta^n+\frac{C}{k}\,,
\end{equation}
where $C_k>0$ can blow up as $k\to +\infty$.
\end{lemma}
\begin{proof}
The proof is an adaptation of that of \cite[Lemma 3.1]{Ansini-Braides}, and we present it in detail for the convenience of the readers.

For every $j\in \N$, $x_{j,i} \in \mathscr{G}_{j,M}$, and $l\in \{0,\dots,k-1\}$ we define the shorthand
\[
C_{j,i}^{l}:=C_{\eps_j,\theta, M}^{l}(\e_j x_{j,i})
\]
and we denote with $\bar u_{j,i}^l, \bar \sigma_{j,i}^l$ the quantities defined as in \eqref{intermediate_values}, with $l$ replacing $k_{j,i}$. Let $\psi_{j,i}^{l}\in C_c^\infty(C_{j,i}^{l})$ be a cut-off function satisfying
\begin{equation}\label{cut_off_properties}
\psi_{j,i}^{l}\big|_{\partial B_{\bar \sigma_{j,i}^l}}\equiv 1\,, \quad 0\leq \psi_{j,i}^{l}\leq 1\,, \quad \|\nabla \psi_{j,i}^{l}\|_{L^\infty(C_{j,i}^{l})}\leq \frac{c}{\bar \sigma_{j,i}^l}\,, 
\end{equation}
with $c>0$, and set
\begin{equation}\label{cut_off_modification}
w_{j,i}^{l}:=u_j+\psi_{j,i}^{l}(\bar u_{j,i}^{l}-u_j)\,.
\end{equation}
Note that $w_{j,i}^{l}\equiv u_j$ outside the annulus $C_{j,i}^{l}$ and $w_{j,i}^{l}=\bar u_{j,i}^{l}$ on $\partial B_{\bar \sigma_{j,i}^{l}}$. By \eqref{growth_of_f} and \eqref{cut_off_properties} we have 
\begin{align}\nonumber
\int_{C_{j,i}^{l}} f(\nabla w_{j,i}^{l})\,\mathrm{d}x&\leq c_2\int_{C_{j,i}^{l}}\Big(\big|\nabla \psi_{j,i}^{l}\otimes(\bar u_{j,i}^{l}-u_j)+(1-\psi_{j,i}^{l})\nabla u_j\big|^q+1\Big)\,\mathrm{d}x\\\label{first_energy_estimate_on_an_annulus}
&\lesssim (\bar \sigma_{j,i}^{l})^{-q}\int_{C_{j,i}^{l}}|u_j-\bar u_{j,i}^{l}|^q\,\mathrm{d}x+\int_{C_{j,i}^{l}}\big(1+|\nabla u_j|^q\big)\,\mathrm{d}x
\lesssim \int_{C_{j,i}^{l}}\big(1+\big|\nabla u_{j}\big|^q\big)\,\mathrm{d}x\,,
\end{align}
where to conclude we have used the Poincar{\' e} inequality in the annulus $C_{j,i}^{l}$. Since the sets $(C_{j,i}^{l})_{l=0,\dots,k-1}$ are pairwise disjoint and
$$\bigcup_{l=0}^{k-1} C_{j,i}^{l}\subset B_{\theta\eps_j/M}(\eps_jx_{j,i})\setminus {\overline B}_{2^{-k}\theta\eps_j/M}(\eps_jx_{j,i})\subset B_{\theta\eps_j/M}(\eps_jx_{j,i})\,,
$$
we obtain for every $j\in \N$ and $x_{j,i} \in \mathscr{G}_{j,M}$, 
\begin{equation*}
\sum_{l=0}^{k-1}\int_{C_{j,i}^{l}}\big(1+\big|\nabla u_j\big|^q\big)\,\mathrm{d}x\leq \int_{B_{\theta\eps_j/M}(\eps_jx_{j,i})}\big(1+\big|\nabla u_{j}\big|^q\big)\,\mathrm{d}x\,.
\end{equation*}
In particular, for fixed $j\in \N$ and $x_{j,i} \in \mathscr{G}_{j,M}$, there exists $k_{j,i}\in \{0,\dots,k-1\}$ so that
\begin{equation*}
\int_{C_{j,i}^{k_{j,i}}}\big(1+\big|\nabla u_j\big|^q\big)\,\mathrm{d}x\leq \frac{1}{k} \int_{B_{\theta \eps_j/M}(\eps_jx_{j,i})}\big(1+\big|\nabla u_{j}\big|^q\big)\,\mathrm{d}x\,.
\end{equation*}
Setting $C_{j,i}:=C_{j,i}^{k_{j,i}}, \ \bar u_{j,i}:=\bar u_{j,i}^{k_{j,i}}, \ \bar\sigma_{j,i}:=\bar \sigma_{j,i}^{k_{j,i}}$ and $w_{j,i}:=w_{j,i}^{k_{j,i}}$, we define the sequence $(w_j)$ as   
\begin{equation}\label{modified_sequence}
w_j:=\begin{cases}
u_j & \text{in } D\setminus \underset{x_{j,i}\in \mathscr{G}_{j,M}}{\bigcup} C_{j,i}\,, \cr
w_{j,i} & \text{in } C_{j,i}\,.
\end{cases}
\end{equation}
From the definition and by \eqref{cut_off_properties} and  \eqref{cut_off_modification} we have that  $w_j$ satisfies properties \ref{properties_of_the_modification_w_j_M-1}-\ref{properties_of_the_modification_w_j_M-2}. Moreover by \eqref{first_energy_estimate_on_an_annulus} - \eqref{modified_sequence}, and the fact that the balls $(B_{\theta\eps_j/M}(\eps_jx_{j,i}))_{x_{j,i}\in \mathscr{G}_{j,M}}$ are pairwise disjoint subsets of $D$, we have 
\begin{align*}
\Big|\int_D f(\nabla w_j)\,\mathrm{d}x-\int_D f(\nabla u_{j})\,\mathrm{d}x\Big|&\leq \sum_{x_{j,i}\in \mathscr{G}_{j,M}}\int_{C_{j,i}}|f(\nabla w_{j,i})-f(\nabla u_j)|\,\mathrm{d}x\,\\
&\lesssim \frac{1}{k}\sum_{x_{j,i}\in \mathscr{G}_{j,M}}\int_{B_{\theta\eps_j/M}(\eps_j x_{j,i})}\big(1+|\nabla u_j|^q\big)\,\mathrm{d}x\\
&\lesssim \frac{1}{k}\int_{D}\big(1+|\nabla u_j|^q\big)\,\mathrm{d}x\,.
\end{align*}
Then, \ref{c:jl_estimate} follows immediately by \eqref{u_j_to_u_weakly}.
	
It only remains to check that $w_j$ satisfies \ref{w_j_M_to_u_weakly}. To do so, we start by showing that $w_j\rightarrow u$ strongly in $L^q(D;\R^m)$, as $j\to +\infty$. Indeed, by \eqref{cut_off_properties}, \eqref{cut_off_modification}, \eqref{modified_sequence}, and by using the Poincar{\' e} inequality in the annuli $C_{j,i}$, we get
\begin{align}\label{conv:wjtou}
\int_D|w_j-u|^q\,\mathrm{d}x&= \int_{D\setminus \underset{x_{j,i}\in \mathscr{G}_{j,M}}{\bigcup}C_{j,i}}\big|u_j-u\big|^q\,\mathrm{d}x+ \sum_{x_{j,i}\in \mathscr{G}_{j,M}}\int_{C_{j,i}}\big|u_j+\psi_{j,i}^{k_{j,i}}(\bar u_{j,i}-u_j)-u\big|^q\,\mathrm{d}x\nonumber\\
&\lesssim  \int_D|u_j-u|^q\,\mathrm{d}x+\sum_{x_{j,i}\in \mathscr{G}_{j,M}}\int_{C_{j,i}} \big|u_j-\bar u_{j,i}\big|^q\,\mathrm{d}x\nonumber\\
&\lesssim \int_{D}\big|u_j-u\big|^q\,\mathrm{d}x+(\theta \eps_j /M)^q\sup_{j\in \N}\int_D |\nabla u_j|^q\,\mathrm{d}x\,,
\end{align}
which is infinitesimal as $j\to +\infty$ by \eqref{u_j_to_u_weakly}, and since $(\eps_j)\searrow 0$. By combining \eqref{modified_sequence}, \eqref{growth_of_f}, and \eqref{first_energy_estimate_on_an_annulus}
we also obtain
\begin{align*}
\int_{D}|\nabla w_j|^q\,\mathrm{d}x & = \int_{D\setminus \underset{x_{j,i}\in \mathscr{G}_{j,M}}{\bigcup}C_{j,i}}|\nabla u_j|^q\,\mathrm{d}x+ \sum_{x_{j,i}\in \mathscr{G}_{j,M}}\int_{C_{j,i}}|\nabla w_{j,i}|^q\,\mathrm{d}x\\
&\lesssim \int_D|\nabla u_j|^q\,\mathrm{d}x+\sum_{x_{j,i}\in \mathscr{G}_{j,M}}\int_{C_{j,i}} (f(\nabla w_{j,i})+1)\,\mathrm{d}x\\
& \lesssim \int_D|\nabla u_j|^q\,\mathrm{d}x+\sum_{x_{j,i}\in \mathscr{G}_{j,M}}\int_{B_{\theta\eps_j/M}(\eps_jx_{j,i})}(1+|\nabla u_j|^q)\,\mathrm{d}x\lesssim \int_{D}(1+|\nabla u_j|^q)\,\mathrm{d}x\,,
\end{align*}
which by \eqref{u_j_to_u_weakly} and \eqref{conv:wjtou} yields the desired convergence \ref{w_j_M_to_u_weakly}. 	

Suppose now that the sequence $(|\nabla u_j|^q)$ is equi-integrable. In this case, for each $x_{j,i}\in \mathscr{G}_{j,M}$ we set 
\begin{equation*}
C_{j,i}:=\{\theta \eps_j/(2M)<|x-\eps_jx_{j,i}|<\theta \eps_j/M\}\,,\quad  
\bar u_{j,i}:=\fint_{C_{j,i}} u_j\,\mathrm{d}x\,, 
\quad \bar \sigma_{j,i}:=\frac{3}{4}\frac{\theta\eps_j}{M}\,,
\end{equation*}
and 
\begin{equation*}
w_{j,i}:=u_j+ \psi_{j,i}(\bar u_{j,i}-u_j)\,,
\end{equation*}
where $\psi_{j,i}\in C_c^\infty(C_{j,i})$ is a cut-off function such that
\begin{equation*}
 0\leq \psi_{j,i}\leq 1\,, \ \ \psi_{j,i}\big|_{\partial B_{\bar \sigma_{j,i}}}\equiv 1\,, \ \ \|\nabla \psi_{j,i}\|_{L^\infty(C_{j,i})}\leq \frac{c}{\bar \sigma_{j,i}}\,, 
\end{equation*}
with $c>0$. Similarly to \eqref{first_energy_estimate_on_an_annulus}, also in this case we have  
\begin{equation*}
\int_{C_{j,i}}f(\nabla w_{j,i})\,\mathrm{d}x\lesssim \int_{C_{j,i}}\big(1+|\nabla u_j|^q\big)\,\mathrm{d}x\,.
\end{equation*}
Setting
\begin{equation*}
w_j:=\begin{cases}
u_j\quad \  \text{in } D\setminus \underset{x_{j,i}\in \mathscr{G}_{j,M}}{\bigcup} C_{j,i}\,,\\
w_{j,i}\ \ \ \ \text{in each } C_{j,i}\,,
\end{cases}
\end{equation*} 
one can easily check that \ref{properties_of_the_modification_w_j_M-1},\ref{properties_of_the_modification_w_j_M-2} and \ref{w_j_M_to_u_weakly} are satisfied. 
To prove \eqref{modified_energy_estimate_equiiintegrability}  note that for every $T>0$ large enough,  
\begin{align*}
\Big|\int_D f(\nabla w_j)\,\mathrm{d}x-\int_D f(\nabla u_{j})\,\mathrm{d}x\Big|&\lesssim \sum_{x_{j,i}\in \mathscr{G}_{j,M}}\int_{B_{\theta \eps_j/M}(\eps_jx_{j,i})}(1+|\nabla u_j|^q)\,\mathrm{d}x\,\\
&\lesssim (1+T^q)(\theta \eps_j/M)^n\#\mathscr{G}_{j,M}+\int_{D\cap \{|\nabla u_j|>T\}}|\nabla u_j|^q\,\mathrm{d}x\,.
\end{align*} 
By \eqref{former:i} we get 
\begin{align*}
\Big|\int_D f(\nabla w_j)\,\mathrm{d}x-\int_D f(\nabla u_{j})\,\mathrm{d}x\Big|\lesssim (1+T^q)\theta^n+\int_{D\cap \{|\nabla u_j|>T\}}|\nabla u_j|^q\,\mathrm{d}x\,.
\end{align*}
For $k\in \N$, taking first $T=T_k>0$ large enough, using the equi-integrability assumption on $(|\nabla u_j|^q)$ we have 
$$
\sup_{j\in \N}\int_{D\cap \{|\nabla u_j|>T_k\}}|\nabla u_j|^q\,\mathrm{d}x\leq \frac{1}{k}
$$ 
and hence \eqref{modified_energy_estimate_equiiintegrability}. 
\end{proof}

\section {Stochastic building blocks} \label{probabilistic_preliminaries}
In this section we collect some probabilistic results that we use in the proof of Theorem \ref{main_thm}. 
Preliminarily we recall that for a bounded set $E\subset \R^n$ we have set
\[
\Phi(E):=\Phi\cap E, \quad N(E):=\# \Phi(E)\,. 
\]
Moreover for $\e>0$ we also define the $\e$-dependent random variables
\begin{align}\label{notation_for_r.v.}
\Phi_\eps(E):=\Phi({\eps}^{-1}E), \quad N_\eps(E):=\# \Phi_\eps(E)
\end{align}
and for $\delta>0$, we introduce the \textit{thinning process} $\Phi^{\delta}$ of $\delta$-isolated centres, defined by
\begin{equation}\label{notation_for_thinning}
\Phi^{\delta,\omega}:=\Big\{x^\omega\in \Phi^\omega\colon \underset{y^\omega\in \Phi^\omega,\ y^\omega\neq x^\omega}{\min}|y^\omega-x^\omega|\geq \delta\Big\}\,.
\end{equation}
Analogously, we define
\begin{align}\label{notation_for_thinnes_r.v.}
\begin{split}
& \Phi^{\delta}(E):=\Phi^{\delta}\cap E\,, \quad \quad \ \Phi_{\eps}^{\delta}(E):=\Phi^{\delta}(\eps^{-1}E\big)\,,\\[2pt]
& N^{\delta}(E):=\#\Phi^{\delta}(E)\,, \quad \ \  N_{\eps}^{\delta}(E):=\#\Phi_{\eps}^{\delta}(E)\,.
\end{split}
\end{align}
Lemma \ref{good_and_bad_balls_lemma} below is a statement on the asymptotic random geometry of the perforations and is a straightforward adaptation of \cite[Lemma 4.2]{Giunti-Hofer-Velasquez} to our setting. Since we deal with functionals with $q$-growth (rather than quadratic), the critical scale of the perforations for us is $\eps^{n/(n-q)}$ (rather than $\eps^{n/(n-2)}$). This difference in the scale causes some minor changes in the statement of the result, but is of no consequence in its proof for which we refer to \cite{Giunti-Hofer-Velasquez} and omit here.

\medskip

In what follows $\Omega'\in \mathcal{T}$ denotes a set with $\mathbb{P}(\Omega')=1$ that may vary from line to line and depends only on the m.p.p. $(\Phi,\mathcal{R})$. If a property holds true for every $\o\in \O'$ we may equivalently write that it holds $\mathbb{P}$-a.e.\ in $\Omega$ or, in short, \emph{almost surely}.

\begin{lemma}\label{good_and_bad_balls_lemma} Let $(\Phi,\mathcal{R})$ be a  m.p.p. satisfying the assumptions \ref{H1}--\ref{H4}, and let $H_\eps^\omega$, for $\eps>0$ and $\omega\in \Omega$, be the family of random holes associated to the m.p.p. defined as in \eqref{random_holes}.

There exist $\eps_0:=\eps_0(n,m,q)>0$, random variables $(r_\eps)$ with $r_\eps: \Omega\to \R_+$, and a set $\Omega' \in\mathcal T$ with $\mathbb P(\Omega')=1$ with the following properties. For every $\omega\in \Omega'$,
\begin{equation}\label{r_eps_are_negligible}
\lim_{\eps\to 0^+}  r_\eps^\omega=0\,, \ 
\end{equation}
and for every $\omega\in \Omega'$ and every $\eps\in (0,\eps_0]$ there exists a set $I^\omega_{\eps,b}\subset \Phi^\omega_\eps(D)$ such that, defining  
$I^\omega_{\eps,g}:= \Phi^\omega_\eps(D) \setminus I^\omega_{\eps,b}$ and
\begin{equation}\label{def:good-bad}
H_{\eps,b}^\omega:=\bigcup_{x_i^\omega\in I^\omega_{\eps,b}} \overline{B}_{\alpha_\eps\rho^\omega_i}(\eps x_i^\omega)\, 
,\ \ D_{\eps,b}^\omega:=\bigcup_{x_i^\omega\in I^\omega_{\eps,b}} \overline{B}_{2\alpha_\eps\rho^\omega_i}(\eps x_i^\omega)\,, 
\quad 
H_{\eps,g}^\omega:=\bigcup_{x_i^\omega\in I^\omega_{\eps,g}} \overline{B}_{\alpha_\eps\rho^\omega_i}(\eps x_i^\omega)\,, 
\end{equation}
we have 
\begin{equation}\label{c:cap}
\mathrm{dist}(H_{\eps,g}^\omega,  D_{\eps,b}^\omega)\geq \frac{\eps r_\eps^\omega}{2}\,,
\quad 
\lim_{\eps\to 0^+}\eps^n \#  I^\omega_{\eps,b} =0\,,
\quad 
\lim_{\eps\to 0^+} \eps^n\sum_{x_i^\omega\in I_{\eps,b}^\omega}(\rho_i^\omega)^{n-q}=0\,,
\end{equation}
\begin{equation}\label{c:13_marzo}
\min_{\substack{x_l^\omega,x_i^\omega\in I_{\eps,g}^\omega\\\ x_l^\omega\neq x_i^\omega}}|x_l^\omega-x_i^\omega|\geq 2 r_\eps^\omega\,, \quad 
\max_{x_i^\omega\in I_{\eps,g}^\omega}\alpha_\eps \rho^\omega_i\leq \frac{\eps r_\eps^\omega}{2}\,, \quad \lim_{\eps\to 0^+} \eps^n\# I_{\eps,g}^\omega=\langle N(Q)\rangle\L^n(D)\,.
\end{equation}
Finally, if $\delta>0$, for the thinning process defined in \eqref{notation_for_thinning}, for every $\omega\in \Omega'$ there holds 
\begin{equation}\label{few_points_close_to_safety_region}
\lim_{\eps\to 0^+} \eps^n\#\{x_i^\omega\in \Phi_{\eps}^{2\delta,\omega}(D)\colon \mathrm{dist}(\eps x_i^\omega, D_{\eps,b}^\omega)\leq \delta\eps\}=0\,.
\end{equation}
\end{lemma}

\begin{figure}
\begin{center}
\begin{tikzpicture}
\node[inner sep=0pt] (dec) at (0,0)
    {\includegraphics[width=.6\textwidth]{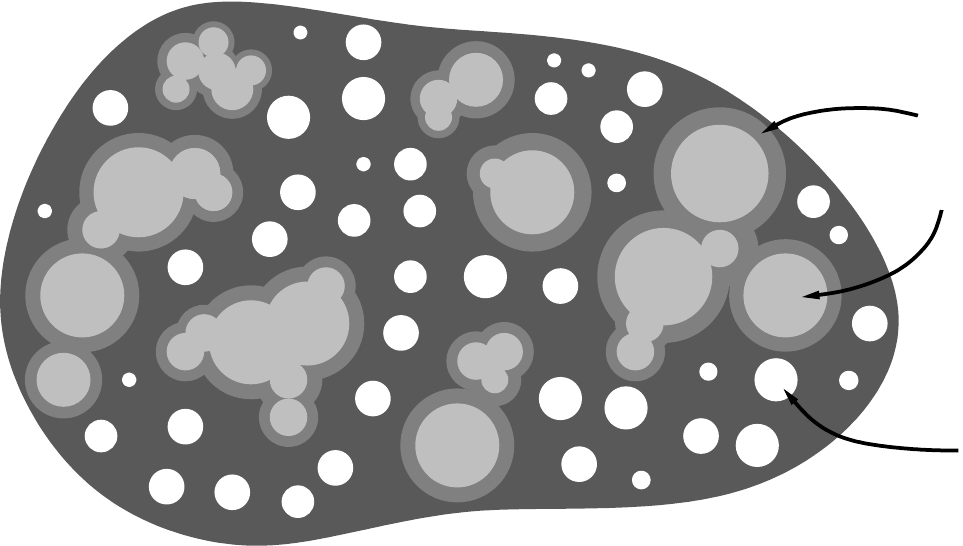}};
\draw(4.4,-1.6)node[right]{good hole};
\draw(4.4,.7)node[right]{bad hole};
\draw(4.2,1.5)node[right]{safety layer};
\end{tikzpicture}
\end{center}
\caption{Domain decomposition into good holes $H^\omega_{\e,g}$, bad holes $H^\omega_{\e,b}$, and safety layer $D^\omega_{\e,b}$.}
\label{fig:decomposition}
\end{figure}

In what follows we refer to the sets $H^\omega_{\eps,g}$ and $H^\omega_{\eps,b}$ in \eqref{def:good-bad} as \textit{good} and \textit{bad perforations}, respectively, while 
$D^\omega_{\eps,b}$ is referred to as the \textit{safety layer}. Note that $H_\eps^\omega=H^\omega_{\eps,g}\cup H^\omega_{\eps,b}$. In short, the good perforations are $(\eps r^\o_\eps/2)$-separated from the safety layer, they are $\eps r^\o_\eps$-separated from one another, their radii are bounded by $(\eps r^\o_\eps)/2$, and  asymptotically they are the only relevant set of perforations. An illustration of the geometry is shown in Figure \ref{fig:decomposition}.

Note that, by \eqref{Capq-increasing}-\eqref{Capq-subadditive}, proceeding as in the proof of \eqref{explicit_K_model_capacity}, we have that for every $\omega \in \Omega'$ and $z\in \R^m$
\begin{align}\label{GHV-lemma}
\mathrm{Cap}_q(H_{\eps,b}^{\omega}, D_{\eps,b}^\omega)&=
\mathrm{Cap}_q\bigg(\bigcup_{x_i^\omega\in I^\omega_{\eps,b}} B_{\alpha_\eps\rho^\omega_i}(\eps x_i^\omega), D_{\eps,b}^\omega\bigg)
 \leq \sum_{x_i^\omega\in I^\omega_{\eps,b}} \mathrm{Cap}_q\big(B_{\alpha_\eps\rho^\omega_i}(\eps x_i^\omega), D_{\eps,b}^\omega\bigg)\nonumber\\
& \leq \sum_{x_i^\omega\in I^\omega_{\eps,b}} \mathrm{Cap}_q\big(B_{\alpha_\eps\rho^\omega_i}(\eps x_i^\omega), B_{2\alpha_\eps\rho^\omega_i}(\eps x_i^\omega)\big)\nonumber\\
&\sim \sum_{x_i^\omega\in I^\omega_{\eps,b}}\big((\alpha_\eps\rho^\omega_i)^{(q-n)/(q-1)}-(2\alpha_\eps\rho^\omega_i)^{(q-n)/(q-1)}\big)^{1-q}\nonumber\\
&\sim \eps^n\sum_{x_i^\omega\in I^\omega_{\eps,b}}(\rho_i^\omega)^{n-q}\to 0\, \quad \text{as } \eps \to 0^+\,,
\end{align}
where we have used \eqref{c:cap}. Condition \eqref{GHV-lemma} (with $\mathrm{Cap}_q$ replaced by the classical harmonic capacity) is explicitly stated in \cite[Lemma 4.2]{Giunti-Hofer-Velasquez} instead of the last equality in \eqref{c:cap}; however the analogue of \eqref{c:cap} can be found in the proof of their result (cf.\ equation (4.58) and the one above it therein). Moreover by following the steps in the proof of \cite[Lemma 4.2]{Giunti-Hofer-Velasquez} it is easy to verify that in our case the random variables $(r_\eps^\omega)$ can be chosen as
\[r_\eps^\omega:=\big(\eps^{n/n-q}\max_{x_i^\omega\in \Phi_\eps^\omega(D)}\rho_i^\omega\big)^{\frac{1}{n}}\vee \eps^{\alpha/4} \  \text{for some }\alpha\in(0,\tfrac{q}{n-q})\,.\] 

Following \cite{Giunti-Hofer-Velasquez}, we now introduce for $\omega\in \Omega'$ the subset of $I^\omega_{\eps,g}$ given by the (centres of the) balls that are deterministically spaced apart from one another and from the safety layer, and have uniformly bounded rescaled  radii $\rho_i^\omega$. More precisely, for $M\in \N$ fixed we define 
\begin{equation}\label{M_good_centers}
G_{\eps,M}^{\omega}:=\big\{x_i^\omega\in I_{\eps,g}^{\omega}\colon d_{\eps,i}^{ \omega}\geq \eps/M\ \ \text{and } \rho_i^{\omega}\leq M\big\}\,,
\end{equation} 
where, for $x_i^\omega\in \Phi^\omega_\eps(D)$, we set
\begin{equation}\label{distance_of_x_i_from_rest}
d_{\eps, i}^{ \omega }:=\min\Big\{\mathrm{dist}(\eps x_i^\omega, D_{\eps,b}^{ \omega }), \tfrac{1}{2}\min_{x_l^\omega\neq x_i^\omega}\eps|x_l^\omega-x_i^\omega|,\eps\Big\}\,.
\end{equation} 
If $(\eps_j)\searrow 0$ we adopt the shorthand notation $G^\omega_{j,M}:=G^\omega_{\eps_j,M}$\,.

 Compared to the good centres $I_{\eps,g}^{\omega}$, where the same scale $\eps r^\o_\eps$ controlled  both the size of the perforations and their separation, for the balls centred at $G_{\eps,M}^{\omega}$ the separation is of order $\eps$, while the size is much smaller, of the critical order $\alpha_\eps$.

Below we show that the family of points $G_{\e,M}^{\omega}$ satisfies, almost surely, properties \ref{c:dis}-\ref{c:card} before \eqref{former:i}. Hence we can deduce a probabilistic version of Lemma  \ref{joining_lemma_on_randomly_perforated_domains} for $\mathbb{P}$-a.e. $\omega\in \Omega$, where sequences will be modified around balls with centres in $G_{\e,M}^{\omega}$.

\begin{lemma}[Probabilistic Joining Lemma]\label{l:probabilistic_joining_lemma}
There exists a set $\Omega' \in \mathcal T$ with $\mathbb P(\Omega')=1$ satisfying the following property. Let $(\eps_j) \searrow 0$, $M\in \mathbb{N}$, $\theta\in(0,1)$ (with $\frac{1}{\theta}\in \N$) and $k\in \N$ be fixed.  
Let $ (u_j)\subset W_0^{1,q}(D;\R^m)$ be such that 
\begin{equation*}
u_j\rightharpoonup u \ \text{ weakly in } W^{1,q}(D;\R^m)\,,
\end{equation*}
for some $u \in W_0^{1,q}(D;\R^m)$, and let $({G}^\omega_{j,M})_j$ be  collections of points in $\R^n$ defined as in \eqref{M_good_centers}, with $\omega\in \Omega$.  
Then for every $\omega \in \Omega'$ and for $x^\omega_{ j,i}\in {G}^\omega_{j,M}$ there exists $k^\omega_{j,i}\in\{0,\dots,k-1\}$ and corresponding annuli  $C_{\eps_j,\theta, M}^{k^\omega_{j,i}}(\e_j x^\omega_{j,i})$, such that 
we can construct a sequence $(w^\omega_j)\subset W_0^{1,q}(D;\R^m)$ satisfying the following properties:
\begin{enumerate}[label=(\roman*)]
\item \label{properties_of_the_modification_w_j_M-1-omega} $w^\omega_j\equiv u_j \ \text{ in } D\setminus \bigcup_{x^\omega_{j,i} \in {G}^\omega_{j,M}} C_{\eps_j,\theta, M}^{k^\omega_{j,i}}(\e_j x^\omega_{j,i})$\,;
\smallskip
\item \label{properties_of_the_modification_w_j_M-2-omega} $w_j^\omega\equiv \bar u^\omega_{j,i} \ \text{ on } \partial B_{ \bar \sigma^\omega_{j,i}}(\eps_j x^\omega_{j,i})$, where 
\begin{equation*}
\bar u^\omega_{j,i}:=\fint_{C_{\eps_j,\theta, M}^{k^\omega_{j,i}}(\e_j x^\omega_{j,i})} u_j\,\mathrm{d}x\,, \qquad  
\bar \sigma^\omega_{j,i}:=\frac{3}{4}2^{-k^\omega_{j,i}}\frac{\theta\eps_j}{M}\,;
\end{equation*}
\item \label{w_j_M_to_u_weakly-omega} $w^\omega_j\rightharpoonup u \ \text{ weakly in } W^{1,q}(D;\R^m)\,;$
\smallskip
\item \label{c:jl_estimate-omega}$\Big|\int_{D} f(\nabla w^\omega_j)\,\mathrm{d}x-\int_{D}f(\nabla u_j)\,\mathrm{d}x\Big|\leq \frac{C}{k}\,$, for some $C>0$ depending on $c_2$, $n$, $m$, $q$, $D$, and $\sup_{j\in \N}\|\nabla u_j\|_{L^q(D;\R^{m\times n})}$.
\end{enumerate}

If, additionally, the sequence $(|\nabla u_j|^q)$ is equi-integrable, then also $(|\nabla w^\omega_j|^q)$ is equi-integrable, and in the definition of $\bar \sigma^\omega_{j,i}$ one can take $k^\omega_{j,i} =0$ for all $x^\omega_{j,i} \in {G}^\omega_{j,M}$, up to replacing \ref{c:jl_estimate-omega} with the following estimate
\begin{equation*}
\Big|\int_{D} f(\nabla w^\omega_j)\,\mathrm{d}x-\int_{D}f(\nabla u_j)\,\mathrm{d}x\Big|\leq C_k\theta^n+\frac{C}{k}\,,
\end{equation*}
where $C_k>0$ can blow up as $k\to +\infty$.
\end{lemma}

\begin{proof}
Let $M\in \mathbb N$ be fixed; by Lemma \ref{good_and_bad_balls_lemma}, there exists $\Omega_M \in\mathcal T$ with $\mathbb P(\Omega_M)=1$ such that for every $\omega\in \Omega_M$ the collection of points $G_{j,M}^{\omega}$ satisfies properties  \ref{c:dis}-\ref{c:card} (note that \ref{c:dis} and \ref{c:card} follow immediately from \eqref{M_good_centers} and \eqref{distance_of_x_i_from_rest}, and by Remark \ref{on_boundary_conditions}). Hence, as in \eqref{former:i}, we get that for every $\omega \in \Omega_M$
\begin{equation}\label{former:iG}
\sup_{j\in \N}(\beta_n M^{-n}\varepsilon_j^n) \# {G}^\omega_{j,M} \leq \mathcal{L}^n(D)\,.
\end{equation}
Finally, set $\Omega':=\bigcap_{M\in \mathbb N} \Omega_M$; clearly $\mathbb P(\Omega')=1$ and for $\omega\in \Omega'$ fixed, Lemma \ref{joining_lemma_on_randomly_perforated_domains} applied to the family of the \textit{truncated good perforations} $(B_{\theta\eps_j/M}(\eps_j x^\omega_{j,i}))_{x^\omega_{j,i} \in G_{j,M}^{\omega}}$ provides us with a sequence $(w_j^{\omega})\subset W_0^{1,q}(D;\R^m)$ enjoying the desired properties. 
\end{proof}

\begin{remark}\label{exclude_points-close_to_safety}
\normalfont We observe that, by \eqref{M_good_centers} and \eqref{distance_of_x_i_from_rest}, the points in $I^\omega_{\eps,g}$ whose distance from the safety layer $D_{\eps, b}^\omega$ is smaller than $\eps/M$ do not belong to the set $G_{\eps,M}^{\omega}$. This guarantees that 
$$
\bigg(\bigcup_{x^\omega_{j,i} \in G_{j,M}^{\omega}} C_{\eps_j,\theta, M}^{k^\omega_{j,i}}(\e_j x^\omega_{j,i})\bigg)\cap \bigg(\bigcup_{x^\omega_{j,i} \in I_{\eps_j,b}^{\omega}} B_{\alpha_{\eps_j}\rho^\omega_{j,i}}(\e_j x^\omega_{j,i})\bigg)=\emptyset\,,
$$
namely that the annuli around the truncated good perforations where the sequence is modified do not touch the bad balls. 

This request is also present in \cite[Equation (4.65)]{Giunti-Hofer-Velasquez} to ensure that the correctors provided by \cite[Lemma 3.1]{Giunti-Hofer-Velasquez} are well-defined.

Note, however, that the annuli around the truncated good perforations, might in principle intersect balls centred at $I_{\eps_j,g}^\omega\setminus G_{j,M}^{\omega}$.
\end{remark}

\subsection{Strong laws of large numbers for marked point processes.}
In this section we state three generalizations of the strong law of large numbers for marked point processes which are relevant for our problem. The first two results, Lemma \ref{probabilistic_lemma_1} and Lemma \ref{probabilistic_lemma_2}, were originally stated and proven in \cite[Section 5]{Giunti-Hofer-Velasquez}. We recall their statements here for the readers' convenience.

\medskip

First, we need to introduce some notation. Let $(\Phi,\mathcal{Y})$ be a m.p.p.\ in $\R^n\times \R_+$, with $\Phi$ satisfying the assumptions \ref{H1}--\ref{H3} of Subsection \ref{Assumptions_on_m.p.p.} and $\mathcal{Y}:=(Y_i)_{x^\omega_i\in \Phi}$, with $Y_i:\Omega\to \R_+$ measurable, satisfying \ref{H4} therein, with \eqref{finite_average_capacity} replaced by 
\begin{equation}\label{first_moment_condition}
\langle Y\rangle:=\int_0^{+\infty}y h(y)\,\mathrm{d}y<+\infty\,,
\end{equation}
and \eqref{short_range_integrable_correlations} replaced by
\begin{equation}\label{short_range_integrable_correlations-y}
|K(r,y_1,y_2)|\leq \frac{C}{(1+r^\gamma)(1+y_1^{s/n-q})(1+y_2^{s/n-q})}\,.
\end{equation}
(Think of $\mathcal{Y}:=(\rho_i^{n-q})$ for our application.)

\begin{lemma}\label{probabilistic_lemma_1}
Let $Q \subset \R^n$ be a unit cube, $(\Phi,\mathcal{Y})$ a m.p.p.\ as above and $B\subset \R^n$ a bounded set star-shaped with respect to the origin. Then, for $\mathbb P$-a.e. $\omega\in \Omega$, 
\begin{equation}\label{1st_law_of_large_numbers_type}
\lim_{\eps\to 0^+}\ \eps^nN_\eps^{\omega}(B)=\langle N(Q)\rangle\L^n(B)
\end{equation}
and
\begin{equation}\label{1st_law_of_large_numbers_type-ii}
\lim_{\eps\to 0^+}\ \eps^n\hspace{-0.5em}\sum_{x^\omega_i\in \Phi_\eps^{\omega}(B)}Y_i^{\omega}=\langle N(Q)\rangle\langle Y\rangle\L^n(B)\,.
\end{equation}
Moreover, for any bounded set $A\subset \R^n$, the thinning process $\Phi^{\delta}$ defined in \eqref{notation_for_thinning} satisfies 
\begin{equation}\label{limit_of_thinnings}
\lim_{\delta\to 0^+}\ \langle N^{\delta}(A)\rangle =\langle N(A)\rangle\,.
\end{equation}
\end{lemma}
A by-product of Lemma \ref{probabilistic_lemma_1} (cf. \cite[Section 5]{Giunti-Hofer-Velasquez}) is the following.
\begin{lemma}\label{probabilistic_lemma_2}
Let $(\Phi,\mathcal{Y})$ be a m.p.p.\ as above and $B\subset \R^n$ a bounded set star-shaped with respect to the origin. Let $I^\omega_\eps\subset \Phi^\omega_\eps(B)$ be such that, for $\mathbb P$-a.e. $\omega\in \Omega$, 
\begin{equation*}
\lim_{\eps\to 0^+}\ \eps^n\#I_\eps^{\omega}= 0\,.
\end{equation*}
Then, for $\mathbb P$-a.e. $\omega\in \Omega$, 
\begin{equation*}
\lim_{\eps\to 0^+}\ \eps^n\sum_{x^\omega_i\in I_\eps^{\omega}}Y_i^{\omega}=0\,.
\end{equation*}	
\end{lemma}

We conclude this section with a technical result which can be seen as the nonlinear counterpart of \cite[Lemma 5.3]{Giunti-Hofer-Velasquez}. This result will then be used to prove a stochastic Riemann-sum approximation for the capacitary term appearing in the $\Gamma$-limit $\F_0$ in \eqref{deterministic_functionals}.

\medskip

Let $M\in \N$ be fixed  and let $\kappa: \R^m\times \R_+\to \R$ be a Borel function bounded from below, such that  
\begin{equation}\label{c:k-gc}
\kappa(0,\cdot)=0, \quad \kappa(z,y)\leq C_{\kappa} (|z|^{q}y^r+1) \ \text{ for } 0<q, r<n\,, 
\end{equation}
and
\begin{equation}\label{c:k-lip}
|\kappa(z_1,y)-\kappa(z_2,y)| \leq C_{\kappa, M}(1+|z_1|^{q-1}+|z_2|^{q-1})|z_1-z_2|  \ \text{ for every } y\in [0,M]\,, 
\end{equation}
for some $C_\kappa, C_{\kappa, M}>0$. 
We observe that $\kappa(\cdot,y)\in W^{1,\infty}_{\rm{loc}}(\R^m)$, with norm uniformly bounded in $y \in [0,M]$\ie for every open and bounded set  $E\subset\subset \R^m$ and for any $y\in [0,M]$, 
\begin{equation}\label{lip-kappa}
|\kappa(z_1,y)-\kappa(z_2,y)| \leq C_{\kappa,M,E}|z_1-z_2|\ \ \forall \, z_1,z_2\in \overline {E}\,,
\end{equation}
for some $C_{\kappa, M,E}>0$.

Let $(\Phi,\mathcal{Y})$ be a m.p.p. as defined above; we furthermore assume that 
\begin{equation}\label{additional-trunc}
Y_i:\Omega\mapsto [0,M]\,. 
\end{equation}
Finally, let $u\in C^\infty_c(D;\R^m)$ be fixed and let
$X^u_i: D\times \Omega \to \R$ be defined as 
\begin{equation}\label{appendix_marks_for_the_application}
X^u_i(x,\omega):= \kappa(u(x), Y_i^\omega)\,.
\end{equation}
Then, $X^u_i(x, \cdot)$ is $(\mathcal{T},\mathcal{B}(\R))$-measurable and $X^u_i(\cdot, \omega)$ is continuous. 
Let $\mathcal{X}^u:=(X^u_i)_{x^\omega_i\in \Phi}$ denote a new family of \emph{space-dependent} marks, and let $(\Phi, \mathcal{X}^u)$ be the corresponding marked point process in $\R^n\times \R_+$. 
Moreover, we assume for the average function $x\mapsto \langle X^u(x,\cdot)\rangle:=\int_{0}^{+\infty}\kappa(u(x),y)h(y)\,\mathrm{d}y$ that 
\begin{equation}\label{new_first_moment_condition}
\langle X^u\rangle=\int_{0}^{+\infty}\kappa(u(\cdot),y)h(y)\,\mathrm{d}y \in L^\infty(D)\,.
\end{equation} 

\begin{proposition}\label{probabilistic_lemma_3}
Let $Q \subset \R^n$ be a unit cube and $(\Phi, \mathcal{X}^u)$ be the m.p.p. in $\R^n\times \R_+$ defined above, where 
$\Phi$ satisfies {\normalfont\ref{H1}--\ref{H3}}, and $X_i^u$ are defined in \eqref{appendix_marks_for_the_application}. Let $r_\eps>0$ be such that 
\begin{equation}\label{r_eps_less_than_eps}
r_\eps\leq C\eps,\ \text{ for some }C>0\,.
\end{equation}
Then there exists $\Omega' \in \mathcal T$ with $\mathbb P(\Omega')=1$ and a subsequence in $\eps>0$ (not relabelled) such that
\begin{equation}\label{strong_integral_law_of_numbers}
\lim_{\eps\to 0^+}\ \eps^n\sum_{x^\omega_i\in \Phi_{\eps}^{2/M,\omega}(D) 
}\fint_{B_{r_\eps}(\eps x^\omega_i)}X^u_i(x,\omega)\,\mathrm{d}x=\langle N^{2/M}(Q)\rangle \int_{D}\langle X^u(x,\cdot)\rangle\,\mathrm{d}x\,,
\end{equation}	
for every $\omega \in \Omega'$, every $M\in \N$, and every $u\in C_c^\infty(D;\R^m)$. 
\end{proposition}

\begin{proof} 
In what follows, by $\lim_{\eps\to 0^+}$ we mean limits taken up to an $\omega$-independent subsequence.  
	
Let $u\in C_c^\infty(D;\R^m)$ and $M\in \N$ be fixed. We split the proof into a number of steps.

\smallskip

\textit{Step 1: Properties of the space-dependent marks.} We observe that $(X^u_i)_{x^\omega_i\in \Phi_\e^\omega(D)}$ satisfy the following properties.
\begin{enumerate}[label=(P$\arabic*$)]
\item\label{P1} For\ every $\omega\in \Omega$  and $x^\omega_i \in \Phi^\omega_{\eps}(D)$ we have  
\begin{equation*}
X^u_i(x, \omega)=0\; \text{ for every }\, x \in \partial D\,.
\end{equation*}
\item\label{P2} There exists $\Lambda:=\Lambda(q,r,\kappa,\|u\|_\infty, M) \in (0,+\infty)$ such that
\begin{equation}\label{c:bound}
\sup_{\omega\in \Omega}\ \sup_{x^\omega_i\in \Phi^\omega_\eps(D)}\|X^u_i(\cdot,\omega)\|_{L^\infty(D)}\leq\Lambda\,.
\end{equation}
\item\label{P3} For every $\omega\in \Omega$ the functions $X^u_i(\cdot, \omega)$ belong to $W^{1,\infty}(D)$, with Lipschitz norm uniformly bounded in $i$ (and similarly for their expected value $\langle X^u(x,\cdot)\rangle$). 
\end{enumerate}
	
\smallskip

Property \ref{P1} follows immediately from $u|_{\partial D}\equiv0$, since for every $\omega\in \Omega$, $x_i^\omega\in \Phi^\omega_{\eps}(D)$, and $x\in \partial D$, by \eqref{c:k-gc} we have
\begin{align*}
\begin{split}
X^u_i(x,\omega)=\kappa(u(x),Y_i^\omega)=\kappa(0,Y_i^\omega)=0\,.
\end{split}
\end{align*}
	
Property \ref{P2} is a consequence of \eqref{c:k-gc}, \eqref{additional-trunc}, and the fact that $u\in C_c^\infty(D;\R^m)$. Indeed, for every $\omega\in \Omega$, $x^\omega_i\in \Phi^\omega_\eps(D)$, and $x\in D$ we have 
\begin{equation*}
|X^u_i(x,\omega)|=|\kappa(u(x),Y_i^\omega)|\leq C_\kappa(|u(x)|^{q}(Y_i^\omega)^{r}+1) \leq C_\kappa (\|u\|_{\infty}^{q}M^{r}+1)\,,
\end{equation*}
thus \eqref{c:bound} is satisfied with $\Lambda:=C_\kappa (\|u\|_{\infty}^{q}M^{r}+1)$, which is independent of both $\omega$ and $i$. 
	
Finally, by \eqref{lip-kappa} and \eqref{additional-trunc}, for every  $\omega\in \Omega$,  $x^\omega_i\in \Phi^\omega_\eps(D)$, and $x,x'\in D$ we deduce that
\begin{align}\nonumber
|X^u_i(x,\omega)-X^u_i(x',\omega)|\ &= |\kappa(u(x), Y_i^\omega)-\kappa(u(x'), Y_i^\omega)|\\[3pt]
&\leq C _{\kappa,M,\|u\|_\infty}|u(x)-u(x')|
\leq  C _{\kappa,M,\|u\|_\infty}\|\nabla u\|_{\infty}|x-x'| \label{equiLipschitzianity1}\,,
\end{align}
and therefore \ref{P3}. We observe that by passing to the expected value, the analogue of \eqref{equiLipschitzianity1} holds true for $\langle X^u(x,\cdot)\rangle$ as well.  	

\smallskip

\textit{Step 2: Replacing $r_\eps$ with $\eps$ in \eqref{strong_integral_law_of_numbers}.} 
If $r_\eps$ satisfies \eqref{r_eps_less_than_eps}, by a change of variables and by \eqref{equiLipschitzianity1} and \eqref{notation_for_thinnes_r.v.} we get 
\begin{align*}\nonumber
&\bigg|\eps^n\sum_{x^\omega_i\in \Phi_{\eps}^{2/M,\omega}(D)}\fint_{B_{r_\eps}(\eps x^\omega_i)}X^u_i(x,\omega)\,\mathrm{d}x-\eps^n\sum_{x^\omega_i\in\Phi_{\eps}^{2/M,\omega}(D)}\fint_{B_{\eps}(\eps x^\omega_i)}X^u_i(x,\omega)\,\mathrm{d}x\bigg| \nonumber \\[3pt]
&\quad \quad =\frac{\eps^n}{\beta_n}\bigg|\sum_{x^\omega_i\in \Phi_{\eps}^{2/M,\omega}(D)}\int_{B_1(0)}\big(X^u_i(\eps x^\omega_i+r_\eps z,\omega)-X^u_i(\eps x^\omega_i+\eps z,\omega)\big)\,\mathrm{d}z\bigg| \nonumber \\[3pt]
&\quad \quad  \leq C|r_\eps-\eps|\eps^n N^{2/M,\omega}_{\eps}(D)\leq C\eps\big(\eps^n N^\omega_{\eps}(D)\big)\underset{\eps\to 0^+}{\longrightarrow}0\,
\end{align*}
for $\PP$-a.e.\ $\o\in \O$, thanks to \eqref{1st_law_of_large_numbers_type}. 
Therefore, to prove \eqref{strong_integral_law_of_numbers} it suffices to show that there exists $ \Omega_M \in \mathcal T$ with $\mathbb P( \Omega_M)=1$ such that, up to subsequences 
\begin{equation}\label{modified_version_at_eps_radius}
\lim_{\eps\to 0^+}\ \sum_{x^\omega_i\in \Phi_{\eps}^{2/M,\omega}(D)}\int_{B_{\eps}(\eps x^\omega_i)}X^u_i(x,\omega)\,\mathrm{d}x=\beta_n\langle N^{2/M}(Q)\rangle\int_{D}\langle X^u(x,\cdot)\rangle\,\mathrm{d}x\,,
\end{equation}
for every $\o\in\O_M$ and every  $u\in C_c^\infty(D;\R^m)$.  

\smallskip

\textit{Step 3: Reducing \eqref{modified_version_at_eps_radius} to a dense subset $\mathcal D \subset C_c^\infty(D;\R^m)$.}  Let $\mathcal D \subset C_c^\infty(D;\R^m)$ be countable and dense with respect to the strong $L^q(D;\R^m)$-topology. Assume that there exists $ \Omega_M \in \mathcal T$ with $\mathbb P( \Omega_M)=1$ such that, up to subsequences, \eqref{modified_version_at_eps_radius} holds true for every $\o\in \O_M$ and every  $u \in \mathcal D$.

Let now $\omega\in \Omega_M$ and $u\in C_c^\infty(D;\R^m)$, and let $(u_k)\subset \mathcal{D}$ be such that $u_k\rightarrow u$ strongly in $L^q(D;\R^m)$ as $k\to +\infty$. We have
\begin{align}\nonumber
&\bigg|\sum_{x^\omega_i\in \Phi_{\eps}^{2/M,\omega}(D)}\int_{B_{\eps}(\eps x^\omega_i)}X^u_i(x,\omega)\,\mathrm{d}x- \beta_n\langle  N^{2/M}(Q)\rangle\int_{D}\langle X^u(x,\cdot)\rangle\,\mathrm{d}x\bigg|
\\\nonumber
& \leq  
\bigg|\sum_{x^\omega_i\in \Phi_{\eps}^{2/M,\omega}(D)}\int_{B_{\eps}(\eps x^\omega_i)}X^{u_k}_i(x,\omega)\,\mathrm{d}x-\beta_n\langle N^{2/M}(Q)\rangle\int_{D}\langle X^{u_k}(x,\cdot)\rangle\,\mathrm{d}x\bigg|
\\\nonumber
&\qquad + 
\sum_{x^\omega_i\in \Phi_{\eps}^{2/M,\omega}(D)}\int_{B_{\eps}(\eps x^\omega_i)}\Big|X^u_i(x,\omega)-X^{u_k}_i(x,\omega)\Big|\,\mathrm{d}x\\
\label{c:claim-dense-sub}
 &\qquad+ \beta_n\langle N^{2/M}(Q)\rangle \int_{D}\Big|\langle X^u(x,\cdot)\rangle -\langle X^{u_k}(x,\cdot)\rangle\Big| \,\mathrm{d}x\,.
\end{align}
The first term in the right-hand side of \eqref{c:claim-dense-sub} converges to zero as $\eps\to 0^+$ by assumption, since $u_k\in \mathcal{D}$. For the second term, by the  definition of $X^u_i$ and by \eqref{c:k-lip} we find that for every $x_i^\omega \in \Phi_{\eps}^{2/M,\omega}(D)$
\begin{align*}
& \quad \int_{B_{\eps}(\eps x^\omega_i)}\Big|X^u_i(x,\omega) - X^{u_k}_i(x,\omega)\Big|\,\mathrm{d}x 
 \leq C_{\kappa,M}\int_{B_{\eps}(\eps x^\omega_i)}\Big(1+ |u|^{q-1}+|u_k|^{q-1}\Big)|u-u_k|\,\mathrm{d}x 
\\
&\qquad \lesssim \Big( \eps^{n/q}+ \|u\|_{L^q(B_{\eps}(\eps x^\omega_i))}+\|u_k\|_{L^q(B_{\eps}(\eps x^\omega_i))}\Big)^{q-1}\|u-u_k\|_{L^q(B_{\eps}(\eps x^\omega_i))}\\
&\qquad \lesssim (\eps^{n/q}+ \|u\|_{L^q(B_{\eps}(\eps x^\omega_i))})^{q-1}\|u-u_k\|_{L^q(B_{\eps}(\eps x^\omega_i))}\,,
\end{align*}
where in the last inequality we used the fact that $u_k\rightarrow u$ strongly in $L^q(D;\R^m)$. 
Adding up the previous inequality over all $x_i^\omega \in \Phi_{\eps}^{2/M,\omega}(D)$, using the fact that each ball $B_\eps(\eps x^\omega_j)$ overlaps with at most a finite number (depending only on $M$) of other balls of the family $(B_\eps(\eps x^\omega_i))_{x^\omega_i\in \Phi_{\eps}^{2/M,\omega}(D)}$, and by a discrete H\"older inequality, we deduce
\begin{align}\label{u_k_goes_to_u_1}
&\sum_{x^\omega_i\in\Phi_{\eps}^{2/M,\omega}(D)} \int_{B_{\eps}(\eps x^\omega_i)}\Big|X^u_i(x,\omega) - X^{u_k}_i(x,\omega)\Big| \nonumber\\
&\lesssim \sum_{x^\omega_i\in\Phi_{\eps}^{2/M,\omega}(D)} (\eps^{n/q}+ \|u\|_{L^q(B_{\eps}(\eps x^\omega_i))})^{q-1}\|u-u_k\|_{L^q(B_{\eps}(\eps x^\omega_i))}\nonumber\\
& \leq \Big(\sum_{x^\omega_i\in \Phi_{\eps}^{2/M,\omega}(D)} (\eps^{n/q}+ \|u\|_{L^q(B_{\eps}(\eps x^\omega_i))})^{q}\Big)^{\frac{q-1}{q}}\Big(\sum_{x^\omega_i\in\Phi_{\eps}^{2/M,\omega}(D)}\int_{B_{\eps}(\eps x^\omega_i)}|u_k-u|^q\Big)^{\frac{1}{q}}\nonumber\\
& \lesssim \Big(\eps^n N_{\eps}^{2/M,\omega}(D)+\sum_{x^\omega_i\in \Phi_{\eps}^{2/M,\omega}(D)}\int_{B_\eps(\eps x^\omega_i)}|u|^q\Big)^{\frac{q-1}{q}}\|u_k-u\|_{L^q(D)}\nonumber\\
& \leq \big((\eps^nN_{\eps}^{\omega}(D))^{q-1/q}+\|u\|_{L^q(D)}^{q-1}\big)\|u_k-u\|_{L^q(D)}\,.
\end{align}
Hence, by \eqref{u_k_goes_to_u_1}, the second term in the right-hand side of \eqref{c:claim-dense-sub} converges to zero as $k\to +\infty$. To conclude, we show that the third term in the right-hand side of \eqref{c:claim-dense-sub} converges to zero as $k\to +\infty$ as well. By the definition \eqref{new_first_moment_condition} and \eqref{c:k-lip} we estimate
\begin{align}\label{u_k_to_u_2}
&\int_{D}\Big|\langle X^u(x,\cdot)\rangle - \langle X^{u_k}(x,\cdot)\rangle\Big| \,\mathrm{d}x\leq \int_{D}\int_{0}^{+\infty} \big|\kappa(u(x),y)-\kappa(u_k(x),y)\big|h(y)\,\mathrm{d}y\,\mathrm{d}x\nonumber \\
&\qquad \leq C_{\kappa,M}\int_D(1+|u_k|^{q-1}+|u|^{q-1})|u-u_k|
\leq C_{\kappa,M}(1+\|u\|_{L^q(D)}^{q-1})\|u_k-u\|_{L^q(D)}\,,
\end{align} 
where we used that $\int_{\R_+}h(y)\,\mathrm{d}y=1$ (see \ref{H4}), and that $u_k \to u$ strongly in $L^q$. This term is indeed infinitesimal as $k\to +\infty$, so the proof of this step is complete.

\smallskip

\textit{Step 4: Proving \eqref{modified_version_at_eps_radius} in $\mathcal D$. } 
Let $v\in \mathcal{D}$. We estimate
\begin{equation*}
\Big|\sum_{x^\omega_i\in  \Phi_{\eps}^{2/M,\omega}(D)}\int_{B_{\eps}(\eps x^\omega_i)}X_i^{v}(x,\o)\dx-\beta_n\langle N^{2/M}(Q)\rangle\int_{D}\langle X^{v}(x,\cdot)\rangle\dx\Big| \leq 
a_\e^{v}(\o)+b_\e^{v}(\o)\,,
\end{equation*}
where
\begin{equation}\label{def:a_eps_o}
a_\e^{v}(\o):= \Big|\sum_{x^\omega_i\in \Phi_{\eps}^{2/M,\omega}(D)}\int_{B_\eps(\eps x^\omega_i)}\langle X^{v}(x,\cdot)\rangle\dx-\beta_n
\langle N^{2/M}(Q)\rangle\int_{D}\langle X^{ v}(x,\cdot)\rangle\dx\Big|\,,
\end{equation}
and 
\begin{equation}\label{def:b_eps_o}
b_\e^{ v}(\o):= \Big|\sum_{x^\omega_i\in \Phi_{\eps}^{2/M,\omega}(D)}\int_{B_{\eps}(\eps x^\omega_i)}\big(X_i^{ v}(x,\omega)-\langle X^{ v}(x,\cdot)\rangle\big)\dx\Big|\,.
\end{equation}
We claim that there exists $\Omega_M^v\in \mathcal{T}$ with $\mathbb{P}(\Omega_M^v)=1$, so that up to a deterministic subsequence in $\eps>0$, 
\begin{equation}\label{c:claim0}
\lim_{\e \to 0^+}a_\e^{ v}(\o)=\lim_{\e \to 0^+} b_\e^{ v}(\o)= 0 \quad \forall \, \omega\in \Omega_M^v\,.
\end{equation}
Note that since $\mathcal{D}$ is countable, by setting $\Omega_M:=\underset{v\in \mathcal{D}}{\bigcap}\Omega_M^v$, we obtain directly \eqref{modified_version_at_eps_radius} in $\mathcal D$.

In what follows, to not overburden the notation, we omit the possible extraction of a subsequence in $\eps>0$, as long as this can be chosen independently of the realisation $\omega\in \Omega$.  

\smallskip

\textit{Substep 4.1: Rewriting \eqref{c:claim0}.} We start by rewriting $a_\e^{v}$. Since, by \ref{P1}, for every $\o\in \O$ and every $x^\omega_i\in \Phi^\omega_{\eps}(D)$ we have
$X_i^{v}(x,\omega)= \langle X^{v}(x,\cdot)\rangle=0$ for every $ x\in \partial D$, by setting $X_i^{v}(x,\omega)=\langle X^{ v}(x,\cdot)\rangle=0$ for every $x\in \R^n\setminus D$, we obtain functions defined in the whole of $\R^n$. We now tessellate $\R^n$ into unitary cubes $\{Q_j\}_{j\in \N}$ with $Q_j:=Q(z_j)$ and $\{z_j\}_{j\in \N}\equiv \Z^n$, and observe that 
\begin{equation}\label{points_eps_close}
|x-\eps z_j|\leq (1+\sqrt{n}/2)\eps, \ \ \text{for every} \ x\in B_\eps(\eps x^\omega_i),\, x^\omega_i\in  \Phi^{2/M,\omega}(Q_j)\,. 
\end{equation}
Set  $\N_\eps(D):=\{j\in \N\colon \eps Q_j\cap D\neq\emptyset\}$; noticing that $\eps Q_j=Q_\e(\e z_j)$, by \eqref{def:a_eps_o}  
we split
\begin{equation}\label{av:split}
a_\e^{v}(\omega)=\Big|\sum_{j\in \N_\eps(D)}\Big(\sum_{x^\omega_i\in \Phi_{\eps}^{2/M,\omega}(\eps Q_j)}\int_{B_\eps(\eps x^\omega_i)}\langle X^{v}(x,\cdot)\rangle\dx-\beta_n
\langle N^{2/M}(Q)\rangle \int_{\eps Q_j}\langle X^{v}(x,\cdot)\rangle\dx\Big)\Big|\,,
\end{equation}
where we have used that $\langle X^{ v}(x,\cdot)\rangle=0$ for every $x\in \R^n\setminus D$.
Now, for every $j\in \N_\eps(D)$ we write 
\begin{align}
&\sum_{x^\omega_i\in \Phi_{\eps}^{2/M,\omega}(\eps Q_j)}\int_{B_\eps(\eps x^\omega_i)}\langle X^{v}(x,\cdot)\rangle\dx-\beta_n
\langle N^{2/M}(Q)\rangle \int_{\eps Q_j}\langle X^{v}(x,\cdot)\rangle\dx \nonumber
\\
&=\sum_{x^\omega_i\in \Phi_{\eps}^{2/M,\omega}(\eps Q_j)}\Bigg(\int_{B_\eps(\eps x^\omega_i)} \Big(
\langle X^{v}(x,\cdot)\rangle - \langle X^{v}(\eps z_j,\cdot)\rangle 
- \beta_n \int_{\eps Q_j} \Big(
\langle X^{v}(x,\cdot)\rangle - \langle X^{v}(\eps z_j,\cdot)\rangle 
\Big)\dx\Bigg)\nonumber\\\label{sum:third:av}
&\qquad + \beta_n \big(N^{2/M,\omega}(Q_j)
-\langle N^{2/M}(Q)\rangle\big) \int_{\eps Q_j}\langle X^{v}(x,\cdot)\rangle \dx\,.
\end{align}
By \eqref{equiLipschitzianity1} (for $\langle X^{v}\rangle$) and \eqref{points_eps_close} we estimate for every $x\in B_\eps(\eps x^\omega_i)$ with $x^\omega_i\in \Phi^{2/M,\omega}(Q_j)=\Phi_{\eps}^{2/M,\omega}(\eps Q_j)$, and similarly for $x\in \eps Q_j$, 
\begin{equation}\label{estimate:BQ}
|\langle X^{v}(x,\cdot)\rangle - \langle X^{v}(\eps z_j,\cdot)\rangle| \leq C _{\kappa,M,\|v\|_\infty}\|\nabla v\|_{\infty}|x-\eps z_j| \leq C _{\kappa,M,\|v\|_\infty}\|\nabla v\|_{\infty} (1+\sqrt{n}/2)\eps\,.
\end{equation}
Hence, from \eqref{av:split}, \eqref{sum:third:av} and \eqref{estimate:BQ}, we can estimate 
\begin{align*}
a_\e^{v}(\omega) &\lesssim  C _{\kappa,M,\|v\|_\infty}\|\nabla v\|_{\infty} \eps \sum_{j\in \N_\eps(D)} \eps^nN_{\eps}^{2/M,\omega}(\eps Q_j)  \\
 &\quad + \Big|\sum_{j\in \N_\eps(D)}\big( N^{2/M,\omega}(Q_j)-\langle N^{2/M}(Q)\rangle\big)\int_{\eps Q_j}\langle X^{v}(x,\cdot)\rangle\dx\Big|\,,\\
 &\lesssim  C(M, v) \eps (\eps^nN_{\eps}^{2/M,\omega}(D))  + \Big|\sum_{j\in \N_\eps(D)}\big( N^{2/M,\omega}(Q_j)-\langle N^{2/M}(Q)\rangle\big)\int_{\eps Q_j}\langle X^{v}(x,\cdot)\rangle\dx\Big|\,.
\end{align*}
In view of \eqref{1st_law_of_large_numbers_type}, to prove \eqref{c:claim0} for $\alpha_\eps^v$, it then suffices to show that there exists $\Omega_M^v\in \mathcal{T}$ with $\mathbb{P}(\Omega_M^v)=1$, such that  
\begin{equation}\label{alternative_2nd_term_in_integral_LLN_goes_to_0}
\lim_{\eps\to 0^+}S_{1,\eps}(\omega)=0 \quad  \forall \, \omega\in \Omega_M^v\,,
\end{equation}
where
\begin{equation}\label{S_1_eps_def}
S_{1,\eps}(\o):=\sum_{j\in \N_\eps(D)} \alpha^j_{\e}\big( N^{2/M,\omega}(Q_j)-\langle N^{2/M}(Q)\rangle\big), \quad \alpha^j_{\eps}:=\int_{\eps Q_j}\langle X^{v}(x,\cdot)\rangle\dx\,.
\end{equation}
We now rewrite \eqref{c:claim0} for $b_\e^{v}$ in \eqref{def:b_eps_o}. First we split
\begin{align*}
b_\e^{v}(\o)&=\Big|\sum_{j\in \N_\eps(D)}\sum_{x^\omega_i\in \Phi_{\eps}^{2/M,\omega}(\eps Q_j)}\int_{B_\eps(\eps x^\omega_i)}\tilde X_i^{v}(x,\omega)\dx \Big|\,,
\end{align*}
where we set $\tilde X_i^{v}(x,\o):=X_{i}^{v}(x,\o)-\langle X^{v}(x,\cdot)\rangle$. Then, proceeding as for $a_\eps^v(\omega)$, by \eqref{equiLipschitzianity1} for $\tilde X_i^{v}$ and \eqref{points_eps_close}, we get, by a similar inequality as \eqref{estimate:BQ},
\begin{align*}
b_\e^{v}(\o) &\lesssim C(M, v)(\eps^nN_\eps^{2/M,\omega}(D))\eps+\eps^n\Big|\sum_{j\in \N_\eps(D)}\sum_{x^\omega_i\in \Phi_{\eps}^{2/M,\omega}(\eps Q_j)}
\tilde X_i^{v}(\e z_j, \o) \Big|\,.
\end{align*} 
Hence, again by \eqref{1st_law_of_large_numbers_type}, to prove \eqref{c:claim0} for $b_\eps^v$, it suffices to show that there exists $\Omega_M^v\in \mathcal{T}$ with $\mathbb{P}(\Omega_M^v)=1$, such that  
\begin{equation}\label{alternative_1st_term_in_integral_LLN_goes_to_0}
\lim_{\eps\to 0^+}S_{2,\eps}(\omega)=0 \quad \forall \, \omega\in \Omega_M^v\,,
\end{equation} 
where
\begin{equation}\label{S_2_eps_def}
S_{2,\eps}(\o):=\eps^n \sum_{j\in \N_\eps(D)}\sum_{x^\omega_i\in \Phi_{\eps}^{2/M,\omega}(\eps Q_j)}
\tilde X_i^{ v}(\e z_j,\o)\,.
\end{equation}

\smallskip

\textit{Substep 4.2: Rewriting \eqref{alternative_2nd_term_in_integral_LLN_goes_to_0} and \eqref{alternative_1st_term_in_integral_LLN_goes_to_0}.} Assume that there exists a subsequence $(\eps_k)_{k}$ independent of the realisation $\omega\in \Omega$, with $\eps_k\searrow 0$ as $k\to +\infty$, such that
\begin{equation}\label{summability_of_second_moments}
\sum_{k=1}^\infty \langle S^2_{l,\eps_k}\rangle<+\infty\,, \  \text{for}\; l=1,2\,.
\end{equation} 
We now show that then \eqref{alternative_2nd_term_in_integral_LLN_goes_to_0} and \eqref{alternative_1st_term_in_integral_LLN_goes_to_0} follow. Indeed, for $l=1,2$ and for every $\delta>0$, by the Chebyshev inequality, \eqref{summability_of_second_moments} yields 
\begin{equation*}
\sum_{k=1}^\infty \mathbb{P}(\{\omega\in \Omega\colon |S_{l,\eps_k}(\omega)|>\delta\})\leq \delta^{-2}\sum_{k=1}^\infty\langle S^2_{l,\eps_k}\rangle<+\infty\,.
\end{equation*}
Therefore, the Borel-Cantelli Lemma ensures that for $l=1,2$ and every $\delta>0$
\begin{equation*}
\mathbb{P}\Big(\bigcap_{m=1}^\infty\bigcup_{k=m}^{\infty}\{\omega\in \Omega\colon |S_{l,\eps_k}(\omega)|>\delta\}\Big)=0\iff \mathbb{P}\Big(\bigcup_{m=1}^\infty\bigcap_{k=m}^{\infty}\{\omega\in \Omega\colon |S_{l,\eps_k}(\omega)|\leq\delta\}\Big)=1\,.
\end{equation*}
Hence, equivalently
\[\PP(\{\o \in \Omega \colon |S_{l,\eps_k}(\omega)|\leq\delta\; \text{for all but a finite number of $k\in \mathbb{N}$}\})=1\,, 
\]
for every $\delta>0$, which ensures that 
\[\mathbb{P}\big(\{\omega\in \Omega\colon \lim_{k \to +\infty}S_{l,\eps_k}(\omega)=0\}\big)=1 \ \ \text{for $l=1,2$}\,,
\]
as desired. As already remarked, the set of events where the last equality holds depends on $M\in \N$ and $v\in \mathcal{D}$, and hence it is an admissible $\Omega_{M}^v$ as in the claims.

\smallskip

\textit{Substep 4.3: Proving \eqref{summability_of_second_moments} for $l=1$.} 
We note that, by \ref{P2}, the deterministic ($v$-dependent) coefficients $(\alpha_{\eps}^j)$ defined in \eqref{S_1_eps_def} satisfy  
\begin{equation}\label{bound_on_beta_e_j}
|\alpha_\eps^j|\leq \|\langle X^{v} \rangle\|_{L^\infty(D)}\eps^n\leq \Lambda\, \eps^n\,,
\end{equation}
with $\Lambda>0$ as in \eqref{c:bound}. Recalling that $\langle X^{v}\rangle\in W^{1,\infty}(D)$ (recall \ref{P3}) and $\langle X^{v}\rangle=0$ in $\R^n\setminus D$, we get
\begin{align}\label{calculation_for_l=1-0}
S^2_{1,\eps}(\o)&=\sum_{j\in \N_\eps(D)}(\alpha_\eps^j)^2\big(N^{2/M,\omega}(Q_j)- \langle N^{2/M}(Q)\rangle \big)^2\nonumber\\
&\quad +\sum_{j\neq j'\in \N_\eps(D)}\alpha_{\eps}^j\alpha_{\eps}^{j'}
\big(N^{2/M,\omega}(Q_j)- \langle N^{2/M}(Q)\rangle\big)\big(N^{2/M,\omega}(Q_{j'})- \langle N^{2/M}(Q)\rangle \big)\,.
\end{align}
Note that the thinning process $\Phi^{2/M}$ inherits the properties \ref{H1}--\ref{H3} in Subsection \ref{Assumptions_on_m.p.p.} from $\Phi$. In particular, the stationarity condition \ref{H1} for $\Phi^{2/M}$ and \eqref{second_moment_condition} yield 
\[\langle N^{2/M}(Q_j)\rangle=\langle N^{2/M}(Q)\rangle, \quad 
\langle (N^{2/M}(Q_j))^2\rangle=\langle (N^{2/M}(Q))^2\rangle\quad \text{for every } j \in \N_\e(D)\,,\] 
and 
\begin{equation}\label{QjQjp}
\langle N^{2/M}(Q_j)N^{2/M}(Q_{j'})\rangle \leq \langle (N^{2/M}(Q))^2\rangle\leq \lambda^2 \quad \text{for every }\,  j,j' \in \N_\e(D)\,.
\end{equation}
Then, taking the expected value in \eqref{calculation_for_l=1-0}, by \eqref{bound_on_beta_e_j}, \eqref{QjQjp} and \eqref{second_moment_condition}, we have 
\begin{align*}
\langle S^2_{1,\eps}\rangle &\lesssim \eps^{2n}\Big(\lambda^2\#\N_\eps(D)+\sum_{j\neq j'\in \N_\eps(D)}|\langle N^{2/M}(Q_j)N^{2/M}(Q_{j'})\rangle- \langle N^{2/M}(Q)\rangle^2|\Big)\\
&\lesssim_\lambda\  \eps^{n}+\eps^{2n}\sum_{\underset{|z_j-z_{j'}|\leq 2\sqrt {n}}{j\neq j'\in \N_\eps(D)}}|\langle N^{2/M}(Q_j)N^{2/M}(Q_{j'})\rangle- \langle N^{2/M}(Q)\rangle^2|\\
&\qquad + \eps^{2n}\sum_{\underset{|z_j-z_{j'}|> 2\sqrt {n}}{j\neq j'\in \N_\eps(D)}}|\langle N^{2/M}(Q_j)N^{2/M}(Q_{j'})\rangle- \langle N^{2/M}(Q)\rangle^2|\,,
\end{align*}
where we have also used the fact that $\#\N_\eps(D)\lesssim \L^n(D)\eps^{-n}$, and split the sum over $j\neq j'\in \N_\eps(D)$ into  contributions ``close to the diagonal'' and ``far from the diagonal''. For the ``close'' contribution, by \eqref{QjQjp} and \eqref{second_moment_condition}  we estimate 
$$
\sum_{\underset{|z_j-z_{j'}|\leq 2\sqrt {n}}{j\neq j'\in \N_\eps(D)}}|\langle N^{2/M}(Q_j)N^{2/M}(Q_{j'})\rangle- \langle N^{2/M}(Q)\rangle^2| \lesssim \lambda^2 \sum_{j\in \N_{\eps}(D)}\#\{j'\colon |z_j-z_{j'}|\leq 2\sqrt{n}\}\lesssim_\lambda \#\N_\eps(D)\,.
$$
For the ``far'' contribution, by  \eqref{strong_mixing} applied, by stationarity, to the random variables $N^{2/M}(Q({z_j}-z_{j'})$ and $N^{2/M}(Q(0))$, we get
\begin{align*}
\sum_{\underset{|z_j-z_{j'}|> 2\sqrt {n}}{j\neq j'\in \N_\eps(D)}}|\langle N^{2/M}(Q_j)N^{2/M}(Q_{j'})\rangle- \langle N^{2/M}(Q)\rangle^2|
&\lesssim \sum_{\underset{|z_j-z_{j'}|> 2\sqrt {n}}{j\neq j'\in \N_\eps(D)}}\frac{C\langle(N^{2/M}(Q))^2\rangle}{1+(|z_{j}-z_{j'}|-\sqrt{n})^\gamma} \\
&\lesssim_\lambda \sum_{\underset{|z_j-z_{j'}|> 2\sqrt {n}}{j\neq j'\in \N_\eps(D)}}\frac{1}{|z_{j}-z_{j'}|^\gamma}\,,
\end{align*}
where we have used the elementary fact that $\mathrm{diam}(Q)=\sqrt{n}$. Hence, we estimate 
\begin{align}\label{mean_value_calculation_for_l=1}
\langle S^2_{1,\eps}\rangle \lesssim_\lambda \eps^n+\eps^{2n}\ \sum_{\underset{|z_j-z_{j'}|> 2\sqrt {n}}{j\neq j'\in \N_\eps(D)}}\frac{1}{|z_{j}-z_{j'}|^\gamma}\,.
\end{align}
To conclude the proof of \eqref{summability_of_second_moments} for $l=1$, we are now left to estimate the second summand in \eqref{mean_value_calculation_for_l=1}.  To this end, note that for every $j\in \N$,
\begin{equation*}
\{j'\in \N\colon |z_j-z_{j'}|>2\sqrt{n}\}\subset  \bigcup_{\ell=1}^{ \infty}\{k\in \N\colon z_k\in \partial Q_{2\ell}(z_j)\}\,.
\end{equation*}
Since $|z_j-z_k|\geq \ell$ for every $k\in \N$ such that $z_k\in \partial Q_{2\ell}(z_j)$, and  
$$\#\{k\in \N\colon z_k\in \partial Q_{2\ell}(z_j)\}\leq (2\ell)^n-(2(\ell-1))^n\lesssim \ell^{n-1}\,,$$
we can estimate 
\begin{align}\label{eq:layers}
\sum_{\underset{|z_j-z_{j'}|> 2\sqrt {n}}{j\neq j'\in \N_\eps(D)}}\frac{1}{|z_{j}-z_{j'}|^\gamma} &\leq 
\sum_{j\in  \N_\eps(D)} \Big(\sum_{\ell\geq 1} \sum_{\underset{z_k\in \partial Q_{2\ell}(z_j)}{k\in \N}}\frac{1}{|z_j-z_k|^\gamma}\Big)\nonumber\\
&\leq \sum_{j\in  \N_\eps(D)}\sum_{\ell\geq 1}\frac{1}{\ell^\gamma}\#\{k\in \N\colon z_k\in \partial Q_{2\ell}(z_j)\} \nonumber\\ 
&\lesssim \left(\sum_{\ell\geq 1} \frac{1}{\ell^{1+(\gamma-n)}}\right)\#\N_\eps(D)\lesssim_{\gamma}\L^n(D) \eps^{-n}\,,
\end{align}
where we used the assumption that $\gamma>n$. Gathering \eqref{mean_value_calculation_for_l=1} and \eqref{eq:layers}, if we choose a deterministic $(\eps_k)\searrow 0$ such that $\sum_{k\in \N} \eps^n_{k}<+\infty$, we obtain
\begin{align*}
\sum_{k=1}^\infty \langle S^2_{1,\eps_k}\rangle\lesssim_{\Lambda,\lambda, D,\gamma} \ \sum_{k=1}^{\infty}\eps_k^{n}<+\infty\,,
\end{align*} 
hence the desired summability condition \eqref{summability_of_second_moments} for $S_{1,\eps_k}$.  
	
\smallskip

\textit{Substep 4.4: Proving \eqref{summability_of_second_moments} for $l=2$.} To simplify the presentation, it is convenient to introduce some shorthand notation for $S_{2,\e}$. Namely, for $j\in \N$ and $\o \in \Omega$, set
\[
Z_j(\o):=\sum_{x^\o_i\in \Phi^{2/M,\o}(Q_j)}\widetilde X_i^{v}(\eps z_j,\o)\,
\]
so that, by \eqref{S_2_eps_def}, we have 
\[
S^2_{2,\e}(\o)=\e^{2n}\bigg(\sum_{j\in \N_\e(D)} Z_j(\o)\bigg)^2=\eps^{2n}\sum_{j\in \N_\eps(D)}Z^2_{j}(\o)+\eps^{2n}\sum_{j\neq j'\in \N_\eps(D)}Z_j(\o)Z_{j'}(\o)\,.
\]
By the definition of $\widetilde X_i^{v}$ after \eqref{S_1_eps_def} and by \eqref{c:bound} for $\widetilde X_i^{v}$, we can estimate the diagonal term in $S^2_{2,\eps}(\o)$ as 
\begin{align*}
\sum_{j\in \N_\eps(D)}Z^2_{j}(\o)\leq \sum_{j\in \N_\eps(D)}\Big(\sum_{x^\o_i\in \Phi^{2/M,\o}(Q_j)}\|\tilde X_i^{v}(\cdot,\o)\|_{L^\infty}\Big)^2\lesssim_\Lambda  \sum_{j\in \N_{\eps}(D)}(N^{2/M,\o}(Q_j))^2\,.
\end{align*}
Therefore, for the expected value of $S^2_{2,\e}$, by \eqref{second_moment_condition}, we get 
\begin{align}\label{mean_value_calculation_for_l=2}
\langle S^2_{2,\eps}\rangle\lesssim_\Lambda \lambda^2\eps^{2n}\#\N_\eps(D)+ \eps^{2n}\sum_{j\neq j'\in \N_\eps(D)}\langle Z_j Z_{j'}\rangle
\lesssim_{\Lambda,\lambda, D}\ \eps^{n}+\eps^{2n}\sum_{j\neq j'\in \N_\eps(D)}\langle Z_j Z_{j'}\rangle\,,
\end{align}
where we have used again that $\#\N_\eps(D)\lesssim \L^n(D)\eps^{-n}$. 
We now estimate the last sum in the right-hand side of \eqref{mean_value_calculation_for_l=2}. We start by splitting again the sum into contributions ``close to the diagonal'' and ``far from the diagonal'', as
\begin{equation}\label{sum_split}
\sum_{j\neq j'\in \N_{\eps}(D)}\langle Z_j Z_{j'}\rangle=\sum_{\underset{|z_j-z_{j'}|\leq 2\sqrt{n}}{j\neq j'\in \N_{\eps}(D)}}\langle Z_jZ_{j'}\rangle+\sum_{\underset{|z_j-z_{j'}|> 2\sqrt{n}}{j\neq j'\in \N_{\eps}(D)}}\langle Z_j Z_{j'} \rangle\,.
\end{equation}
By the definition of $Z_j$ and by \ref{P2}, we can estimate 
$$
Z_j Z_{j'} = \sum_{x^\omega_i\in \Phi^{2/M,\o}(Q_j)}\ \sum_{x^\omega_{i'}\in \Phi^{2/M,\o}(Q_{j'})}\tilde X_i^{v}(\eps z_j,\cdot)\tilde X_{i'}^{v}(\eps z_{j'},\cdot)\leq \Lambda^2 N^{2/M}(Q_j)N^{2/M}(Q_{j'})\,,
$$
hence, by \eqref{QjQjp} and \eqref{second_moment_condition}, 
\begin{align}\label{Z_jZ_j'_trivial bound}
\langle Z_j Z_{j'}\rangle  \leq \Lambda^2 \langle N^{2/M}(Q_j)N^{2/M}(Q_{j'})\rangle 
\leq (\Lambda\lambda)^2\,. 
\end{align}
For the contribution close to the diagonal in \eqref{sum_split}, using \eqref{Z_jZ_j'_trivial bound} we obtain
\begin{align}\label{sum_split-close}
\sum_{\underset{|z_j-z_{j'}|\leq 2\sqrt{n}}{j\neq j'\in \N_{\eps}(D)}}\langle Z_jZ_{j'}\rangle
&=\sum_{j\in \N_{\eps}(D)}\sum_{\underset{ j'\neq j, \, |z_j-z_{j'}|\leq 2\sqrt{n}}{j' \in \N_{\eps}(D)}}\langle Z_j Z_{j'}\rangle
\nonumber\\
&\leq (\Lambda\lambda)^2\sum_{j\in \N_{\eps}(D)}\#\{z_{j'}\in \Z^n\colon |z_{j'}-z_{j}|\leq 2\sqrt{n}\}
\lesssim_{\Lambda,\lambda,n}\#\N_\eps(D) \,.
\end{align} 
Hence, from \eqref{mean_value_calculation_for_l=2}, by \eqref{sum_split} and \eqref{sum_split-close}, we have
\begin{equation}\label{break}
\langle S^2_{2,\eps}\rangle\lesssim_{\Lambda,\lambda,n,D} \ \eps^n + \eps^{2n}\sum_{\underset{|z_j-z_{j'}|> 2\sqrt{n}}{j\neq j'\in \N_{\eps}(D)}}\langle Z_j Z_{j'} \rangle\,.
\end{equation}
We are now left to estimate the sum in the right-hand side of  \eqref{break}, namely the contribution far from the diagonal. Let then $j,j'\in \N_\e(D)$ be such that $|z_j-z_{j'}|>2\sqrt{n}$ and recall that 
\begin{equation}\label{the beast}
\langle Z_j Z_{j'} \rangle = \bigg\langle\sum_{x^\omega_i\in \Phi^{2/M,\o}(Q_j)}\ \sum_{x^\omega_{i'}\in \Phi^{2/M,\o}(Q_{j'})}\tilde X_i^{v}(\eps z_j,\cdot)\tilde X_{i'}^{v}(\eps z_{j'},\cdot)\bigg\rangle\,. 
\end{equation}
Now let $x\in Q_j$ and $x'\in Q_{j}'$ be fixed; for later use it is convenient to estimate the conditional expectation of the random variable $\tilde X_i^{v}(\eps z_j,\omega)\tilde X_{i'}^{v}(\eps z_{j'},\omega)$, conditioned to the locations of the corresponding points $(x_i^\omega, x^\omega_{i'})$ of the point process being fixed, and given by the vector $(x, x')$. Namely we estimate 
\[
\Big\langle\tilde X_i^{v}(\eps z_j,\cdot)\tilde X_{i'}^{v}(\eps z_{j'},\cdot) | (x_i, x_{i'})=(x, x')\Big\rangle\,. 
\] 
Recalling the definition of $X_i^v$ in \eqref{appendix_marks_for_the_application}, we have 
$$
 \tilde X_i^{v}(\eps z_j,\omega) = \kappa(v(\eps z_j), Y_i^\omega) - \langle \kappa(v(\eps z_j), Y_i) \rangle =: \tilde\kappa(v(\eps z_j), Y_i^\o)\,.
$$
Hence, by assumption \ref{H4} in Subsection \ref{Assumptions_on_m.p.p.} on the marks $(Y_i)_{x^\omega_i\in \Phi}$, we have 
\begin{align}\nonumber
\langle \tilde X_i^{v}(\eps z_j,\cdot)\tilde X_{i'}^{v}(\eps z_{j'},\cdot) | (x_i, x_{i'})=(x, x') \rangle 
&= \int_{[0,M]^2}\tilde\kappa(v(\eps z_j), y)\tilde\kappa(v(\eps z_{j'}), y') f_2\big((x,y),(x',y')\big)
\dy\dy' \nonumber\\
&= \int_{[0,M]^2}\tilde\kappa(v(\eps z_j), y)\tilde\kappa(v(\eps z_{j'}), y') h(y) h(y')\dy\dy' \nonumber \\
&+ \int_{[0,M]^2}\tilde\kappa(v(\eps z_j), y)\tilde\kappa(v(\eps z_{j'}), y')K(|x-x'|,y,y')
\dy\dy'\,.\label{XiXip}
\end{align}
\EEE
Note that 
$$
\int_{[0,M]^2}\tilde\kappa(v(\eps z_j), y)\tilde\kappa(v(\eps z_{j'}), y') h(y) h(y')\dy\dy' = 
\langle \tilde\kappa(v(\eps z_j), Y_i)\rangle 
\langle \tilde\kappa(v(\eps z_{j'}), Y_{i'})\rangle =0\,,
$$
since by definition $\tilde \kappa$ has zero average. To estimate the last integral in \eqref{XiXip}, first note that for every $x\in Q_{j}$ and $x'\in Q_{j'}$ there holds
\begin{align}\label{good_cuboids_spacing}
|z_j-z_{j'}|\leq |z_j-x|+|x- x'|+|x'-z_{j'}|\leq \sqrt{n}+|x-x'|\leq \tfrac{1}{2}|z_{j}-z_{j'}|+|x-x'|\,,  
\end{align}
since $|z_j-z_{j'}|>2\sqrt{n}$. Hence, for every $x\in Q_{j}$ and $x'\in Q_{j'}$, we have that, by \eqref{short_range_integrable_correlations-y} and \eqref{good_cuboids_spacing},
\begin{align*}
|K(|x-x'|,y,y')|&\leq \frac{C}{(1+|x-x'|^\gamma)(1+y^{s/n-q})(1+(y')^{s/n-q})}\\
&\lesssim_{\gamma} 
\frac{C}{(1+|z_j-z_{j'}|^\gamma)(1+y^{s/n-q})(1+(y')^{s/n-q})}\,.
\end{align*}
From \eqref{XiXip}, using \ref{P2} and recalling that $s>n-q$, we then deduce that
\begin{align}\label{l=2_4_estimate}
\langle \tilde X_i^{v}(\eps z_j,\cdot)\tilde X_{i'}^{v}(\eps z_{j'},\cdot) | (x_i, x_{i'})=(x, x') \rangle 
&\lesssim_\gamma\
 \frac{C \Lambda^2}{(1+|z_j-z_{j'}|^\gamma)}\bigg(\int_{[0,M]}\frac{1}{1+y^{s/n-q}}\dy\bigg)^2 \nonumber\\
& \lesssim_\gamma \
\frac{C \Lambda^2}{(1+|z_j-z_{j'}|^\gamma)}
\bigg(1+\int_{1}^{+\infty}\frac{\mathrm{d}y}{y^{s/n-q}}\,\bigg)^2  \nonumber
\\
&\lesssim_{\Lambda,\gamma}\ \frac{1}{|z_j-z_{j'}|^\gamma}\,,
\end{align}  
for every $x \in Q_j$ and $x' \in Q_{j'}$. 
\EEE
We are now ready to estimate \eqref{the beast}.

To this end,  it is convenient to introduce the auxiliary random variables 
\[
W^{k,\ell}_{j,j'}(\o):=
\sum_{ \substack {x^\o_i\in \Phi^{2/M,\o}(Q_j), \\ x_{i'}^\o\in \Phi^{2/M,\o}(Q_{j'})}}
\tilde X_{i}^{v}(\eps z_j,\o)\tilde X_{i'}^{v}(\eps z_{j'},\o)\,\chi_{\{N^{2/M,\o}(Q_j)=k,N^{2/M,\o}(Q_{j'})=\ell\}} 
\]
so that, from \eqref{the beast},
\begin{equation}\label{beastkell}
\langle Z_j Z_{j'} \rangle = \sum_{k,\ell=0}^\infty\langle W^{k,\ell}_{j,j'} \rangle\,.
\end{equation}
In this case, we write more explicitly 
$$
\Phi^{2/M,\o}(Q_j)=\{a_{1}^\o,\dots, a_{k}^\o\}, \quad   
\Phi^{2/M,\o}(Q_{j'})=\{b_{1}^\o,\dots, b_{\ell}^\o\}\,,
$$
and define the vector-valued random variables (where we chose an ordering of the elements of the sets)
$$
\underline a_{k}:= (a_{1},\dots, a_{k}):\Omega \longrightarrow (Q_j)^k, \quad \underline b_{\ell}:=(b_{1},\dots, b_{\ell}):\Omega  \longrightarrow (Q_{j'})^\ell.
$$
By the properties of conditional expectation
\begin{equation*}
\langle W^{k,\ell}_{j,j'} \rangle = 
\Big\langle 
\langle W^{k,\ell}_{j,j'} |  (\underline a_{k}, \underline b_{\ell}) 
\rangle
\Big\rangle\,.
\end{equation*}
Hence, denoting with $\Psi^{k,\ell}:(\R^n)^k\times (\R^n)^\ell\to \R_{+}$ the conditional density of $W^{k,\ell}_{j,j'}$ relative to 
$(\underline a_{k}, \underline b_{\ell})$, we have 
\begin{align}\label{28June1050}
\langle W^{k,\ell}_{j,j'} \rangle 
&= \int_{\Omega}\langle {W}^{k,\ell}_{j,j'}| (\underline a^\o_{k}, \underline b^\o_{\ell})\rangle \, \d\mathbb{P}(\o)\nonumber\\
&= \int_{(Q_j)^k\times (Q_{j'})^\ell}\langle {W}^{k,\ell}_{j,j'}|  (\underline a^\o_{k}, \underline b^\o_{\ell})=(\underline\xi_{k},  \underline\zeta_{\ell})\rangle\Psi^{k,\ell}(\underline\xi_{k},  \underline\zeta_{\ell})\, \d (\underline\xi_{k} ,\underline\zeta_{\ell})\,,
\end{align}
where 
$$
\underline\xi_{k} = (\xi_1,\dots, \xi_k) \in (Q_j)^k, \quad \underline\zeta_{\ell}:=(\zeta_{1},\dots,\zeta_{\ell})\in (Q_{j'})^\ell\,.
$$
We now write the inner expectation in the integral above more explicitly. Up to ordering the components of $\underline \xi_k$ and $\underline \zeta_\ell$, we have
\begin{multline}\label{l=2_2_estimate}
\langle {W}^{k,\ell}_{j,j'}|  (\underline a^\o_{k}, \underline b^\o_{\ell}) =(\underline\xi_{k}, \underline \zeta_{\ell})\rangle \\
=\Big\langle \sum_{i,i'=1}^{k,\ell}\Big(\tilde X_{i}^{v}(\eps z_j,\cdot)\tilde X_{i'}^{v}(\eps z_{j'},\cdot) |  (x_i, x_{i'})=(\xi_{i}, \zeta_{i'})\Big) \chi_{\{N^{2/M}(Q_j)=k,N^{2/M}(Q_{j'})=\ell\}} \Big\rangle 
\\
=\sum_{i,i'=1}^{k,\ell}\Big\langle \Big(\tilde X_{i}^{v}(\eps z_j,\cdot)\tilde X_{i'}^{v}(\eps z_{j'},\cdot) |  (x_i, x_{i'})=(\xi_{i}, \zeta_{i'})\Big) \chi_{\{N^{2/M}(Q_j)=k,N^{2/M}(Q_{j'})=\ell\}} \Big\rangle\,,
\end{multline}
where now we can use the linearity of the expectation inside since the sum is deterministic. Moreover, there holds
\begin{align}\label{28June1052}
&\Big\langle \Big(\tilde X_{i}^{v}(\eps z_j,\cdot)\tilde X_{i'}^{v}(\eps z_{j'},\cdot) |  (x_i, x_{i'})=(\xi_{i}, \zeta_{i'})\Big)\chi_{\{N^{2/M}(Q_j)=k,N^{2/M}(Q_{j'})=\ell\}} \Big\rangle\nonumber\\
&\quad = \mathbb{P}(N^{2/M}(Q_j)=k,N^{2/M}(Q_{j'})=\ell)\Big\langle\Big(\tilde X_i^{v}(\eps z_j,\cdot)\tilde X_{i'}^{ v}(\eps z_{j'},\cdot)|  (x_i, x_{i'})=(\xi_{i}, \zeta_{i'})\Big)\Big\rangle\,.
\end{align}
\EEE
Then, gathering \eqref{28June1050}-\eqref{28June1052} and the estimate \eqref{l=2_4_estimate} we have 
\begin{align*}
\langle W^{k,\ell}_{j,j'} \rangle 
& \lesssim_{\Lambda,\gamma} \  \mathbb{P}(N^{2/M}(Q_j)=k,N^{2/M}(Q_{j'})=\ell)\sum_{i,i'=1}^{k,\ell}\frac{1}{|z_j-z_{j'}|^\gamma}\,\nonumber\\
&\lesssim_{\Lambda,\gamma} \ k\ell  \,\mathbb{P}(N^{2/M}(Q_j)=k,N^{2/M}(Q_{j'})=\ell)\frac{1}{|z_j-z_{j'}|^\gamma}\,,
\end{align*}
where we also used that 
\[\int_{(\R^n)^{k}\times (\R^n)^\ell}\Psi^{k,\ell}(\underline\xi_{k},  \underline\zeta_{\ell})\, \d (\underline\xi_{k} ,\underline\zeta_{\ell})=1\,.\]
Hence, by \eqref{beastkell}, \eqref{QjQjp} and \eqref{second_moment_condition} we conclude that 
\begin{align*}
\langle Z_j Z_{j'} \rangle &\lesssim \frac{1}{|z_j-z_{j'}|^\gamma} \sum_{k,\ell=0}^\infty  k\ell \, \mathbb{P}(N^{2/M}(Q_j)=k,N^{2/M}(Q_{j'})=\ell)\nonumber\\
&= \frac{1}{|z_j-z_{j'}|^\gamma} \langle N^{2/M}(Q_j)N^{2/M}(Q_{j'})\rangle 
\leq \frac{1}{|z_j-z_{j'}|^\gamma}\langle (N^{2/M}(Q))^2\rangle 
\lesssim_{\lambda}\frac{1}{|z_j-z_{j'}|^\gamma}\,.
\end{align*}
Finally, from \eqref{break} we get
\begin{equation*}
\langle S^2_{2,\eps}\rangle\lesssim_{\Lambda,\lambda,n,D} \eps^n + \eps^{2n}\sum_{\underset{|z_j-z_{j'}|> 2\sqrt{n}}{j\neq j'\in \N_{\eps}(D)}}\frac{1}{|z_j-z_{j'}|^\gamma}\,,
\end{equation*}
hence the claim follows again by \eqref{eq:layers}, arguing as in the end of Substep 4.3.
\end{proof}

\begin{remark} \normalfont
Proposition \ref{probabilistic_lemma_3} is the key result in the identification of the limit capacitary term in the stochastic $\Gamma$-convergence result Theorem \ref{main_thm}. This identification will be done via successive approximations, and Section \ref{discrete_capacity_approximation} (in particular Proposition \ref{main_prop_for_approximating_the_capacity}) will be devoted to this. 

We want to flag up that while the statement of Proposition \ref{probabilistic_lemma_3} contains the thinning process $\Phi_\eps^{2/M}$, the results in Section \ref{discrete_capacity_approximation} are formulated for centres in the set $G_{\eps,M}$ in \eqref{M_good_centers}, and $G_{\eps,M}\subset \Phi_\eps^{2/M}$. The technical reason for using $\Phi_\eps^{2/M}$ in Proposition \ref{probabilistic_lemma_3} is that the  thinning process inherits the properties \ref{H1}--\ref{H3} from $\Phi$, while $G_{\eps,M}$, which  depends further on the random safety layer $D_{\eps,b}^\omega$, does not.  However, by  \eqref{few_points_close_to_safety_region}, this choice is of no consequence in the proof of the results of the next sections. 
\end{remark}


\section{Discrete approximation of the limit capacitary term}\label{discrete_capacity_approximation}

The main result of this section is Proposition \ref{main_prop_for_approximating_the_capacity} below, where we state that the capacitary term in \eqref{deterministic_functionals} can be obtained as the limit of a ``random'' Riemann sum of the auxiliary capacities \eqref{K_j_rho_truncation_of_capacity}, where the sum is restricted to the perforations centred in the set defined in \eqref{M_good_centers}.

\begin{proposition}\label{main_prop_for_approximating_the_capacity}
Let $ (\eps_j)\searrow 0$, and let $(u_j), u\in W_0^{1,q}(D;\R^m)\cap L^\infty(D;\R^m)$ satisfy 
\begin{equation}\label{u_j_bounded_weakly_convergent}
L:=\sup_{j\in \N}\|u_j\|_{L^{\infty}(D;\R^m)}<+\infty\,,
\end{equation}
and $u_j\rightharpoonup u$ weakly in $W^{1,q}(D;\R^m)$. 
Then, there exist $\Omega' \in \mathcal T$ with $\mathbb P(\Omega')=1$ such that for every $\omega\in \Omega'$ we have (possibly along $\omega$-independent subsequences) 
\begin{equation}\label{main_discrete_to_continuum_convergence}
\lim_{M\to +\infty}\lim_{\theta\to 0^+}\lim_{j\to +\infty}\eps_j^n\sum_{x^\omega_{j,i}\in G_{j,M}^{\omega}}\varphi^j_{\theta,\rho_{j,i}^{ \omega }}({\bar u}_{j,i}^\omega)=\langle N(Q)\rangle\int_{D}\varphi(u)\,\mathrm{d}x\,,
\end{equation}
where $G_{j,M}^{\omega}$ is as in \eqref{M_good_centers},  $\varphi_{\theta,\rho_{j,i}^{\omega}}^j$ as in \eqref{K_j_rho_truncation_of_capacity}, with $\rho$ replaced by $\rho_{j,i}^{ \omega }$ (namely the mark associated to $x_{j,i}^\omega$),  ${\bar u}_{j,i}^\omega$ as in Lemma \ref{l:probabilistic_joining_lemma}, and $\varphi$ is defined in \eqref{expected_g_capacity}.
\end{proposition}

The rest of this section is devoted to the proof of Proposition \ref{main_prop_for_approximating_the_capacity}. This will be carried out in a number of intermediate steps, via successive approximations.

\medskip

In what follows we assume that the sequence $(u_j)$ converges to $u$ pointwise $\L^n$-a.e. in $D$, which can be achieved up to passing to a subsequence. Moreover, throughout this section we always assume that $1/\theta\in \N$, and that $j\in \N$ is so large that \eqref{larger_j} holds true. Finally, we do not relabel the ($\omega$-independent) subsequences along which the limits as $j\to +\infty$ and as $\theta\to 0^+$ are taken.

\medskip

Our first step, Lemma \ref{approximation_with_or_without_j} below, shows that in the discrete approximation \eqref{main_discrete_to_continuum_convergence}, the auxiliary capacity $\varphi^j_{\theta,\rho}$ can be replaced by its limit in $j$, namely by the capacity $\varphi_{\theta,\rho}$ defined in \eqref{local_uniform_convergence_phi_j}.

\begin{lemma}\label{approximation_with_or_without_j} 
Under the same assumptions and notational conventions as in Proposition \ref{main_prop_for_approximating_the_capacity}, there exists $\Omega' \in \mathcal T$ with $\mathbb P(\Omega')=1$ such that for every $\omega\in \Omega'$, we have 
\begin{equation}\label{differece_with_without_j_only_negligible}
\lim_{\theta\to 0^+}\lim_{j\to +\infty}\eps_j^n\sum_{x^\omega_{j,i}\in G_{j,M}^{\omega}}\big|\varphi^j_{\theta,\rho_{j,i}^{ \omega}}({\bar u}_{j,i}^\omega)-\varphi_{\theta,\rho_{j,i}^{ \omega}}({\bar u}_{j,i}^\omega)\big|=0\,,
\end{equation}
for every $M\in \N$, where $\varphi_{\theta,\rho_{j,i}^{ \omega }}$ is defined via \eqref{local_uniform_convergence_phi_j}, with $\rho$ replaced by $\rho_{j,i}^{ \omega }$.\EEE
\end{lemma}
\vspace{-0.5em}
\begin{proof}
Let $M\in \N$ be fixed. For $j\in \N$ and $\theta\in (0,1)$, we define the function $\beta_{\theta}^j:(0,M]\mapsto \R_+$ as 
\begin{equation}\label{moduli_of_continuity}
\beta_{\theta}^j(\rho):=\|\varphi^j_{\theta,\rho}-\varphi_{\theta,\rho}\|_{L^\infty(B^m_L(0))}\,.
\end{equation}
In view of Corollary \ref{local_uniform_convergence_for_capacity_densities} and the fact that in \eqref{moduli_of_continuity} an essential sup is involved, the functions $\beta_{\theta}^j$ are  measurable. Moreover, by \eqref{local_uniform_convergence_phi_j}, $\beta_{\theta}^j \to 0$ pointwise in $(0,M]$ as $j\to +\infty$.

Let now $\tau\in (0,1)$ be fixed. By the Egoroff Theorem there exists a measurable set $J_{\theta,M,\tau}\subset (0,M]$ such that
\begin{equation}\label{application_of_Egorrof}
\mathcal{L}^1(J_{\theta,M,\tau})>M-\tau \quad \text{and } \  \| \beta_{\theta}^j\|_{L^\infty(J_{\theta,M,\tau})}\rightarrow 0 \; \text{ as }\; j\to +\infty\,. 
\end{equation}
We now start estimating the sum in \eqref{differece_with_without_j_only_negligible}. 
For fixed $j\in \N$ and $\omega\in \Omega$,  we get
\begin{align}\label{1st_estimate-start}
\sum_{x^\omega_{j,i}\in  G_{j,M}^{\omega}}|\varphi^j_{\theta,\rho_{j,i}^{\omega}}({\bar u}_{j,i}^\omega)-
\varphi_{\theta,\rho_{j,i}^{\omega}}({\bar u}_{j,i}^\omega
)|
\leq&\   \sum_{x^\omega_{j,i}\in G_{j,M}^\omega}\big(|\varphi^j_{\theta,\rho_{j,i}^{\omega}}({\bar u}_{j,i}^\omega)|+|\varphi_{\theta,\rho_{j,i}^{\omega}}({\bar u}_{j,i}^\omega
)|\big)\chi_{(0,M]\setminus J_{\theta,M,\tau}}(\rho_{j,i}^\omega)\nonumber \\
&+\sum_{x^\omega_{j,i}\in G_{j,M}^\omega}\big|\varphi^j_{\theta,\rho_{j,i}^\omega}({\bar u}_{j,i}^\omega)-\varphi_{\theta,\rho_{j,i}^\omega}({\bar u}_{j,i}^\omega
)\big|\chi_{J_{\theta,M,\tau}}(\rho_{j,i}^\omega)\,,
\end{align}
since by \eqref{M_good_centers} $\rho^\omega_{j,i}\leq M$ for $x^\omega_{j,i}\in G_{j,M}^\omega$. Moreover, for every $x^\omega_{j,i}\in G_{j,M}^\omega$, and for $j\in \N$ large enough, by \eqref{lower_bound_for_the_truncated_capacities}, \eqref{upper_bound_for_the_truncated_capacities}, \eqref{c:c-star} and \eqref{phi_theta_rho__bounds} we obtain
\begin{align}\label{bound_abs_val_of_phi_j}
|\varphi^j_{\theta,\rho_{j,i}^\omega}(\bar u_{j,i}^\omega)|+|\varphi_{\theta,\rho_{j,i}^\omega}(\bar u_{j,i}^\omega)|
&\lesssim |\bar u_{j,i}^\omega|^q(\rho_{j,i}^\omega)^{n-q}+\theta^n\,, 
\end{align}
for some $C>0$. Hence, from \eqref{1st_estimate-start}, using \eqref{bound_abs_val_of_phi_j}, \eqref{u_j_bounded_weakly_convergent} and \eqref{moduli_of_continuity}, and since $G_{j,M}^\omega\subset \Phi_{\eps_j}^\omega(D)$, we get
\begin{align}\label{1st_estimate_on_phi_j_phi_wo_j}
&\sum_{x^\omega_{j,i}\in  G_{j,M}^{\omega}}|\varphi^j_{\theta,\rho_{j,i}^{\omega}}({\bar u}_{j,i}^\omega)-
\varphi_{\theta,\rho_{j,i}^{\omega}}({\bar u}_{j,i}^\omega
)|\nonumber\\
\lesssim&\ \sum_{x^\omega_{j,i}\in  G_{j,M}^\omega}\big(L^q(\rho_{j,i}^\omega)^{n-q}+\theta^n\big)\chi_{(0,M]\setminus J_{\theta,M,\tau}}(\rho_{j,i}^\omega)+\sum_{x^\omega_{j,i}\in  G_{j,M}^\omega} \beta^j_{\theta}(\rho_{j,i}^\omega)\chi_{J_{\theta,M,\tau}}(\rho_{j,i}^\omega)\nonumber \\
\lesssim&\ L^q\sum_{x^\omega_{j,i}\in \Phi_{\eps_j}^\omega(D)}(\rho_{j,i}^\omega)^{n-q}\chi_{(0,M]\setminus J_{\theta,M,\tau}}(\rho_{j,i}^\omega)
+\big(\theta^n+\| \beta^j_{\theta}\|_{L^{\infty}(J_{\theta,M,\tau})}\big)(\#G^\omega_{j,M})\,.
 \end{align}
Define now the random variables $Y_{j,i}:\Omega\to \R_+$ as
\[
Y_{j,i}^\omega:= (\rho_{j,i}^\omega)^{n-q}\chi_{(0,M]\setminus J_{\theta,M,\tau}}(\rho_{j,i}^\omega);
\]
then Lemma \ref{probabilistic_lemma_1} (see \eqref{1st_law_of_large_numbers_type-ii}), applied to  $(Y^\omega_{j,i})$, ensures that there exists $\Omega_{M,\theta,\tau}\in \mathcal{T}$ with $\mathbb{P}(\Omega_{M,\theta,\tau})=1$ such that for every $\omega\in \Omega_{M,\theta,\tau}$, 
\begin{equation}\label{c:lln}
\lim_{j\to +\infty}\ \eps_j^n\hspace{-0.5em}\sum_{x^\omega_{j,i}\in \Phi_{\eps_j}^{\omega}(D)}Y_{j,i}^{\omega}=\langle N(Q)\rangle \L^n(D)\langle Y\rangle
= \langle N(Q)\rangle \L^n(D) \int_0^{+\infty}\rho^{n-q}\chi_{(0,M]\setminus J_{\theta,M,\tau}}(\rho) h(\rho) \,\mathrm{d}\rho\,.
\end{equation}
Note that, by \eqref{finite_average_capacity}, the function $y: \R_+ \to \R_+$ defined as $y(\rho):=\rho^{n-q}h(\rho)$ satisfies $y\in L^1(\R_+)$.

Let $\Omega'\in \mathcal{T}$ be defined as
\begin{equation}\label{1_full_prob_set}
\Omega':=\bigcap_{M,\theta^{-1},\tau^{-1}\in \N}\Omega_{M,\theta,\tau};
\end{equation}
clearly $\mathbb{P}(\Omega')=1$, since $\Omega'$ is a countable intersection of sets of probability 1.
By \eqref{application_of_Egorrof}, \eqref{1st_estimate_on_phi_j_phi_wo_j}, \eqref{c:lln} and \eqref{1st_law_of_large_numbers_type}, for every $\omega\in \Omega'$ and for every $T>0$, we find
\begin{align*}
\underset{j\to+\infty}{\overline{\lim}}\eps_j^n\sum_{x^\omega_{j,i}\in G_{j,M}^\omega}|\varphi^j_{\theta,\rho_{j,i}^\omega}({\bar u}_{j,i}^\omega)-\varphi_{\theta,\rho_{j,i}^\omega}({\bar u}_{j,i}^\omega)|&\lesssim_{L,\Phi,D}\int_{(0,M]\setminus J_{\theta,M,\tau}}\rho^{n-q}h(\rho)\,\mathrm{d}\rho+\theta^n\nonumber \\
&\lesssim T\L^1((0,M]\setminus J_{\theta,M,\tau})+\int_{\{y(\rho)>T\}}y(\rho)\,\mathrm{d}\rho+\theta^n\nonumber \\
&\lesssim_{L} T\tau+\int_{\{y(\rho)>T\}}y(\rho)\,\mathrm{d}\rho+\theta^n\,,
\end{align*}
which holds true for every $\tau\in (0,M)$, and $\theta\in(0,1)$. By taking first the limit as $\tau\to 0^+$, then as $T\to +\infty$, by using the fact that $y\in L^1(\R_+;\R_+)$ and finally by taking the limit as $\theta\to 0^+$, we obtain \eqref{differece_with_without_j_only_negligible} in the set $\Omega'$ defined in \eqref{1_full_prob_set}, and hence the claim.
\end{proof}

In the following lemma we derive a qualitative result for the capacity $\varphi_{\theta,\rho}$.

\begin{lemma}\label{averages_on_balls}
Under the same assumptions and notational conventions as in Proposition \ref{main_prop_for_approximating_the_capacity}, there exists $\Omega' \in \mathcal T$ with $\mathbb P(\Omega')=1$ such that for every $\omega\in \Omega'$, we have 
\begin{equation}\label{differece_between_constant_and_average}
\lim_{j\to +\infty}\eps_j^n\sum_{x_{j,i}^\omega\in G_{j,M}^\omega}\Big|\varphi_{\theta,\rho_{j,i}^{ \omega }}({\bar u}_{j,i}^\omega)-\fint_{B_{\theta\eps_j/M}(\eps_j x^\omega_{j,i})}\varphi_{\theta,\rho_{j,i}^{ \omega}}(u_j)\,\mathrm{d}x\Big|=0\,,
\end{equation}
and
\begin{equation}\label{differece_between_average_u_j_u}
\lim_{j\to +\infty}\eps_j^n\sum_{x^\omega_{j,i}\in  G_{j,M}^\omega}\Big|\fint_{B_{\theta\eps_j/M}(\eps_j x^\omega_{j,i})}\varphi_{\theta,\rho_{j,i}^{ \omega}}(u_j)\,\mathrm{d}x-\fint_{B_{\theta\eps_j/M}(\eps_jx^\omega_{j,i})}\varphi_{\theta,\rho_{j,i}^{\omega}}(u)\,\mathrm{d}x\Big|=0\,,
\end{equation}
for every fixed $M\in \N$ and $\theta\in(0,1)$ with $1/\theta\in \N$. 
\end{lemma}
\begin{proof}
Let $j\in \N$ be large enough so that \eqref{larger_j} holds true, and  $k\in \N$ be fixed. Then, by \eqref{phi_theta_rho_Lipschitz_estimate} and by the bound  \eqref{u_j_bounded_weakly_convergent} that extends to $({\bar u}^\omega_{j,i})$, for every $\omega \in \Omega$  we can estimate
\begin{align}\label{const_vs_mean_value_1st_est}
&\eps_j^n\sum_{x^\omega_{j,i}\in G_{j,M}^\omega}\Big|\varphi_{\theta,\rho_{j,i}^\omega}({\bar u}_{j,i}^\omega)-\fint_{B_{\theta\eps_j/M}(\eps_j x^\omega_{j,i})}\varphi_{\theta,\rho_{j,i}^\omega}(u_j)\,\mathrm{d}x\big|\nonumber\\
&\leq \eps_j^n\sum_{x^\omega_{j,i}\in G_{j,M}^\omega}\fint_{B_{{\theta\eps_j/M}}(\eps_j x^\omega_{j,i})}\big|\varphi_{\theta,\rho_{j,i}^\omega}(u_j)-\varphi_{\theta,\rho_{j,i}^\omega}({\bar u}_{j,i}^\omega)\big| \,\mathrm{d}x\nonumber\\ 
&\leq \eps_j^n\sum_{x^\omega_{j,i}\in G_{j,M}^\omega}C_{M}\fint_{B_{\theta\eps_j/M}(\eps_j x^\omega_{j,i})}\big(\theta^{n(q-1)/q}+|u_j|^{q-1}+|{\bar u}_{j,i}^\omega|^{q-1}\big)|u_j-{\bar u}_{j,i}^\omega|\,\mathrm{d}x\nonumber \\
&\lesssim C_M\big(L^{q-1}+\theta^{n(q-1)/q}\big)\eps_j^n\sum_{x^\omega_{j,i}\in G_{j,M}^\omega}\fint_{B_{\theta\eps_{j}/M}(\eps_j x^\omega_{j,i})}|u_j-{\bar u}_{j,i}^\omega|\,\mathrm{d}x\,,
\end{align}
where we used that, by \eqref{M_good_centers}, $\rho_{j,i}^\omega\leq M$ for every $x^\omega_{j,i}\in G_{j,M}^\omega$. Let now  
${C}_{j,i}^{\omega}:= C_{\eps_j,\theta,M}^{k_{j,i}^\omega}(\eps x_{j,i}^\omega)$ be the annulus defined as in Lemma \ref{l:probabilistic_joining_lemma} and,  
for every $ x^\omega_{j,i}\in G_{j,M}^\omega$, set 
$$
\tilde u_{j,i}^{\omega}:=\fint_{B_{\theta\eps_{j}/M}(\eps_jx^\omega_{j,i})}u_j\,\mathrm{d}x\,.
$$ 
By H\"older's inequality, and since $ {C}_{j,i}^{\omega}\subset B_{\theta\eps_{j}/M}(\eps_jx^\omega_{ j,i})$, we can  estimate the integral in the right-hand side of \eqref{const_vs_mean_value_1st_est} as
\begin{align*}
&\fint_{B_{\theta\eps_{j}/M}(\eps_jx^\omega_{j,i})}|u_j-{\bar u}_{j,i}^\omega|\, \mathrm{d}x
\leq \Big(\fint_{B_{\theta\eps_{j}/M}(\eps_jx^\omega_{ j,i})}|u_j-{\bar u}_{j,i}^\omega|^q\, \mathrm{d}x\Big)^{1/q}
\\[2pt]
&\lesssim \Big(\fint_{B_{\theta\eps_{j}/M}(\eps_j{x^\omega_{j,i}})}|u_j-\tilde u_{j,i}^{\omega}|^q\, \mathrm{d}x+\fint_{{C}_{j,i}^{\omega}}|\tilde u_{j,i}^{\omega}-{\bar u}_{j,i}^{\omega}|^q\, \mathrm{d}x\Big)^{1/q}\nonumber\\[2pt]
&\lesssim\bigg(\fint_{B_{\theta\eps_{j}/M}(\eps_j x^\omega_{j,i})}|u_j-\tilde u_{j,i}^{\omega}|^q+\frac{1}{\L^n({C}_{j,i}^{\omega})}\int_{B_{\theta\eps_{j}/M}(\eps_jx^\omega_{j,i})}|u_j-\tilde u_{j,i}^{\omega}|^q+\fint_{{C}_{j,i}^{\omega}}|u_j- {\bar u}_{j,i}^{\omega}|^q\bigg)^{1/q}.
\end{align*}
Hence, by Poincar{\' e}'s inequality, and since $\mathcal{L}^n({C}_{j,i}^{\omega})$ and $\diam({C}_{j,i}^{\omega})$ differ from the corresponding quantities of $B_{\theta\eps_{j}/M}(\eps_jx^\omega_{j,i})$ only by a multiplicative constant (which can be bounded from above and below uniformly in $i,j$ and $\omega$), we have 
\begin{align*}
\fint_{B_{\theta\eps_{j}/M}(\eps_jx^\omega_{j,i})}|u_j-{\bar u}_{j,i}^\omega|\, \mathrm{d}x \lesssim_{k} (\theta\eps_j/M)^{1-\tfrac{n}{q}}\Big(\int_{B_{\theta\eps_{j}/M}(\eps_jx^\omega_{ j,i})}|\nabla u_j|^q\Big)^{1/q}\,.
\end{align*}
By adding up the estimate above over all $i$ such that $x^\omega_{ j,i}\in G_{j,M}^\omega$ and by the discrete H\"older inequality, we have
\begin{align}\label{intermediate-LL}
\sum_{x^\omega_{j,i}\in G_{j,M}^\omega}\fint_{B_{\theta\eps_{j}/M}(\eps_jx^\omega_{j,i})}|u_j-{\bar u}_{j,i}^\omega|\,
&\lesssim_{k} (\theta\eps_j/M)^{1-\tfrac{n}{q}}\sum_{x^\omega_{j,i}\in G_{j,M}^\omega}\Big(\int_{B_{\theta\eps_{j}/M}(\eps_jx^\omega_{j,i})}|\nabla u_j|^q\Big)^{1/q}\,\nonumber\\
&\lesssim_{k} (\theta\eps_j/M)^{1-\tfrac{n}{q}} (\# G_{j,M}^{\omega})^{1-\tfrac{1}{q}}\Big(\sum_{x^\omega_{j,i}\in G_{j,M}^\omega}\int_{B_{\theta\eps_{j}/M}(\eps_jx^\omega_{ j,i})}\big|\nabla u_j|^q\Big)^{\tfrac{1}{q}}\nonumber\\
&\lesssim_{k} (\theta\eps_j/M)^{1-\tfrac{n}{q}} (\# G_{j,M}^{\omega})^{1-\tfrac{1}{q}}\|\nabla u_j\|_{L^q(D)}\,,
\end{align}
where in the last step we used the fact that the good perforations $(B_{\theta\eps_j/M}(\eps_j x^\omega_{j,i}))_{x^\omega_{j,i}\in G_{j,M}^\omega}$ are pairwise disjoint subsets of $D$. By  \eqref{const_vs_mean_value_1st_est} and 
\eqref{intermediate-LL} we then obtain the bound
\begin{align}\label{const_vs_mean_value_2nd_est}
\eps_j^n\sum_{x^\omega_{j,i}\in G_{j,M}^\omega}\Big|\varphi_{\theta,\rho_{j,i}^{\omega}}({\bar u}_{j,i}^\omega)-\fint_{B_{\theta\eps_{j}/M}(\eps_j x^\omega_{j,i})}\varphi_{\theta,\rho_{j,i}^{\omega}}(u_j)\big|
\lesssim_{\theta, M, L, k} \ \eps_j\big(\eps_j^n\# G_{j,M}^\omega\big)^{1-\tfrac{1}{q}}\|\nabla u_j\|_{L^q(D)}\,.
\end{align}
In view of \eqref{former:iG} and since $(u_j)$ is weakly convergent in $W^{1,q}(D;\R^m)$, by taking in \eqref{const_vs_mean_value_2nd_est} the limit (superior) as $j\to+\infty$, we deduce \eqref{differece_between_constant_and_average}. 

The proof of \eqref{differece_between_average_u_j_u} follows similarly. In particular, using again \eqref{phi_theta_rho_Lipschitz_estimate} and \eqref{u_j_bounded_weakly_convergent}, we estimate 
\begin{align*}
&\eps_j^n\sum_{x^\omega_{j,i}\in G_{j,M}^\omega}\Big|\fint_{B_{\theta\eps_{j}/M}(\eps_jx^\omega_{ j,i})}\Big(\varphi_{\theta,\rho_{j,i}^{\omega}}(u_j)-\varphi_{\theta,\rho_{j,i}^{ \omega}}(u)\Big)\Big|\,\mathrm{d}x\\
&\leq C_M\eps_j^n\sum_{x^\omega_{j,i}\in G_{j,M}^\omega}\fint_{B_{\theta\eps_{j}/M}(\eps_j x^\omega_{j,i})}\big(\theta^{n(q-1)/q}+|u_j|^{q-1}+|u|^{q-1}\big)|u_j-u|\,\mathrm{d}x\nonumber\\
&\lesssim C_M\big(L^{q-1}+\theta^{n(q-1)/q}\big)\eps_j^n(\theta\eps_{j}/M)^{-n}\sum_{x^\omega_{j,i}\in G_{j,M}^\omega}\int_{B_{\theta\eps_{j}/M}(\eps_jx^\omega_{ j,i})}|u_j-u|\,\mathrm{d}x\nonumber\\
&\lesssim_{L,M,\theta,D}\  \int_{D}|u_j-u|\,\mathrm{d}x \lesssim_{L, M, \theta, D} \|u_j-u\|_{L^q(D)}\,.
\end{align*}
Then, since $u_j$ is weakly convergent to $u$ in $W^{1,q}(D;\R^m)$,  taking again the limit (superior)  as $j\to +\infty$  in the above estimate, we obtain \eqref{differece_between_average_u_j_u}\,. 
\end{proof}

We are now ready for the proof of Proposition \ref{main_prop_for_approximating_the_capacity}. 

\begin{proof}[Proof of Proposition \ref{main_prop_for_approximating_the_capacity}]
In view of Lemmata \ref{approximation_with_or_without_j} and \ref{averages_on_balls} it suffices to show that  there exists $\Omega'\in \mathcal{T}$ with $\mathbb{P}(\Omega')=1$ and subsequences in $j\in \N$ and $\theta\in (0,1)$ (not relabelled), such that for all $\omega\in \Omega'$, and every $u\in W_0^{1,q}(D;\R^m)\cap L^\infty(D;\R^m)$,
\begin{equation}\label{reduced_limit_to_show}  
\lim_{M\to +\infty}\lim_{\theta\to 0^+}\lim_{j\to +\infty} \eps_j^n\sum_{x^\omega_{j,i}\in G_{j,M}^\omega}\fint_{B_{\theta\eps_{j}/M}(\eps_j x^\omega_{j,i})}\varphi_{\theta,\rho_{j,i}^{ \omega }}(u)\,\mathrm{d}x=\langle N(Q)\rangle\int_{D}\varphi(u)\,\mathrm{d}x\,.
\end{equation}
We will prove \eqref{reduced_limit_to_show} in a number of steps.

\smallskip

\textit{Step 1: Extending the sum in \eqref{reduced_limit_to_show} to a thinning process with well-separated perforations upon truncation.} For $M\in \N$, $j\in \N$ and ~$\omega\in \Omega$, let  $\Phi^{2/M,\omega}_{\eps_j}(D)$  denote the thinning process defined in \eqref{notation_for_thinnes_r.v.}. For every $x^\omega_{j,i}\in \Phi^{2/M,\omega}_{\eps_j}(D) $ we define the truncated radius of the corresponding perforation as $\rho_{j,i,M}^{\omega}:=\min\{\rho_{j,i}^{\omega},M\}$. (Note that for $x^\omega_{j,i}\in G_{j,M}^\omega$ we have that $\rho_{j,i}^{\omega}\leq M$, so the truncation is not necessary.) We claim that  there exists $\Omega'\in \mathcal{T}$ with $\mathbb{P}(\Omega')=1$ such that for every $\omega\in \Omega'$, 
\begin{equation}\label{limits_in_the whole_process}
\lim_{
\substack{
j\to +\infty \\
\theta\to 0^+\\
M\to +\infty
}}
\eps_j^n\Big|\sum_{x^\omega_{j,i}\in G_{j,M}^\omega}\fint_{B_{\theta\eps_{j}/M}(\eps_jx^\omega_{j,i})}\varphi_{\theta,\rho_{j,i}^{\omega }}(u) - \sum_{x^\omega_{j,i}\in \Phi^{2/M,\omega}_{\eps_j}(D)} \fint_{B_{\theta\eps_{j}/M}(\eps_jx^\omega_{ j,i})}\varphi_{\theta,\rho_{j,i,M}^{ \omega }}(u)\Big|=0\,.
\end{equation}

To see this, we start by estimating the sum over the centres in the thinning process 
 $\Phi^{2/M,\omega}_{\eps_j}(D)$  which are not in the good set $G_{j,M}^\omega$. For fixed $j\in \N$, by \eqref{phi_theta_rho__bounds}, \eqref{u_j_bounded_weakly_convergent}, \eqref{notation_for_r.v.},  Lemma \ref{good_and_bad_balls_lemma}, and recalling the definitions of $I_{\eps_j,b}^\omega$ and $I_{\eps_j,g}^\omega$ therein, we have
\begin{align}\label{for_whole_process_1st_estimate}
&\eps_j^n \bigg|\sum_{x^\omega_{ j,i}\in \Phi^{2/M,\omega}_{\eps_j}(D) \setminus G_{j,M}^{\omega} } \fint_{B_{\theta\eps_{j}/M}(\eps_jx^\omega_{ j,i })}\varphi_{\theta,\rho_{j,i,M}^{ \omega}}(u)\,\mathrm{d}x \bigg|
\nonumber\\
&\lesssim \eps_j^n\!\!\sum_{x^\omega_{ j,i}\in \Phi^{2/M,\omega}_{\eps_j}(D)\setminus G_{j,M}^{\omega} } \Big((\rho_{j,i,M}^{\omega })^{n-q}\fint_{B_{\theta \eps_{j}/M}(\eps_jx^\omega_{ j,i})}|u|^q\,\mathrm{d}x+\theta^n\Big)\nonumber\\
&\lesssim L^q\eps_j^n\!\!\sum_{x^\omega_{ j,i}\in  \Phi^{2/M,\omega}_{\eps_j}(D)\setminus G_{j,M}^{\omega} } (\rho_{j,i,M}^{\omega})^{n-q}+\theta^n(\eps_j^nN^\omega_{\eps_j}(D))\nonumber\\
&\lesssim_L \eps_j^n\!\!\sum_{x^\omega_{j,i}\in I_{\eps_j,b}^{\omega}}(\rho_{j,i}^{\omega})^{n-q}+\eps_j^n\!\!\sum_{x^\omega_{j,i}\in ( \Phi^{2/M,\omega}_{\eps_j}(D)\cap I_{\eps_j,g}^\omega)\setminus G_{j,M}^\omega}(\rho_{j,i,M}^{\omega})^{n-q}+\theta^n(\eps_j^nN^\omega_{\eps_j}(D))\,.
\end{align}
By \eqref{c:cap}, \eqref{1st_law_of_large_numbers_type} (and the fact that we will send $\theta\to 0^+$), we just need to focus on the second term in the right-hand side of \eqref{for_whole_process_1st_estimate} and show it is infinitesimal as $j\to +\infty$, $\theta\to 0^+$ and $M\to +\infty$. We split
$$
(\Phi^{2/M,\omega}_{\eps_j}(D) \cap I_{\eps_j,g}^\omega)\setminus G_{j,M}^\omega = I_{j,M}^{\omega}\cup J_{j,M}^{\omega}\,,
$$
where (recalling \eqref{M_good_centers}, \eqref{distance_of_x_i_from_rest}) 
\begin{align}\label{new_index_sets}
\begin{split}
I_{j,M}^{\omega}&:=\{ x^\omega_{j,i}\in  \Phi^{2/M,\omega}_{\eps_j}(D)\cap I_{\eps_j,g}^\omega \colon \rho_{j,i}^{\omega}>M\}\,,\\
J_{j,M}^{\omega}&:=\big\{ x^\omega_{j,i}\in \Phi^{2/M,\omega}_{\eps_j}(D)\cap I_{\eps_j,g}^\omega \colon \rho_{j,i}^{\omega}\leq M,\quad \mathrm{dist}(\eps_j x^\omega_{j,i}, D_{\eps_j,b}^{\omega})<{\eps_j}/{M}\big\}\,,
\end{split}
\end{align}
so that 
\begin{equation}\label{splitt}
\eps_j^n\!\!\sum_{x^\omega_{j,i}\in ( \Phi^{2/M,\omega}_{\eps_j}(D)\cap I_{\eps_j,g}^\omega)\setminus G_{j,M}^\omega}(\rho_{j,i,M}^{\omega})^{n-q}=\eps_j^n\!\!\sum_{x^\omega_{j,i}\in I_{j,M}^{\omega}}(\rho_{j,i,M}^{\omega})^{n-q}+\eps_j^n\!\!\sum_{x^\omega_{j,i}\in J_{j,M}^{\omega}}(\rho_{j,i,M}^{\omega})^{n-q}.
\end{equation}
For the sum over $I_{j,M}^{\omega}$, by \eqref{1st_law_of_large_numbers_type-ii} (applied to the random variables $Y_{j,i}^\omega:= (\rho_{j,i}^\omega)^{n-q}\chi_{(M,+\infty)}(\rho_{j,i}^\omega)$)  there exists $\Omega_M\in \mathcal{T}$ with $\mathbb{P}(\Omega_M)=1$, such that for every $\omega\in\Omega_M$ we have 
\begin{align}\label{for_whole_process_2nd_estimate}
\underset{j\to+\infty}{\overline{\lim}} \eps_j^n
\sum_{x_{j,i}^\omega\in I_{j,M}^{\omega} }
(\rho_{j,i}^{\omega})^{n-q}&\leq 
\underset{j\to+\infty}{\overline{\lim}} \eps_j^n\sum_{x^\omega_{j,i}\in \Phi_{\eps_j}^{ \omega }(D)}(\rho_{j,i}^{\omega})^{n-q}\chi_{(M,+\infty)}(\rho_{j,i}^{\omega })\nonumber\\
&=\langle N(Q)\rangle\L^n(D)
\int_M^{+\infty} \rho^{n-q}h(\rho)\,\mathrm{d}\rho\,.
\end{align}
Let $\widetilde\Omega:=\underset{M\in \N}{\bigcap}\Omega_M$; note that $\mathbb{P}(\widetilde \Omega)=1$. Let $\omega\in \widetilde \Omega$; by \eqref{for_whole_process_2nd_estimate} and \eqref{finite_average_capacity} we deduce that 
\begin{equation}\label{I_j_M_estimate}
\underset{M\to+\infty}{\overline{\lim}}\ \underset{j\to+\infty}{\overline{\lim}} \eps_j^n \sum_{x^\omega_{j,i}\in I_{j,M}^{\omega} }(\rho_{j,i}^{\omega})^{n-q}\leq \langle N(Q)\rangle\L^n(D)\underset{M\to +\infty}{\overline{\lim}} \int_M^{+\infty} \rho^{n-q}h(\rho)\,\mathrm{d}\rho=0\,.
\end{equation}
For the sum over $J_{j,M}^{\omega}$ in \eqref{splitt}, we first note that by \eqref{few_points_close_to_safety_region} (with $\delta:=M^{-1}$), there exists $\Omega^M\in \mathcal{T}$ with $\mathbb{P}(\Omega^M)=1$, such that for every $\omega\in \Omega^M$, 
\begin{equation}\label{balls-in-J}
\lim_{j\to+\infty}\eps_j^n\#J_{j,M}^{\omega}=0\,.
\end{equation}
Hence, we can apply Lemma \ref{probabilistic_lemma_2} to the m.p.p. $(\Phi,\mathcal{Y})$, with $Y_{j,i}^\omega:=(\rho_{j,i}^{\omega})^{n-q}$, over the set of indices $J_{j,M}^\omega$ in \eqref{new_index_sets}, to obtain that for every $\omega \in \Omega^M$
\begin{equation}\label{J_j_M_estimate}
 \underset{j\to+\infty}{{\lim}}\eps_j^n\sum_{x^\omega_{ j,i}\in J_{j,M}^{\omega}} (\rho_{j,i}^{\omega})^{n-q}=0\,.
\end{equation}
Now we set $\widetilde\Omega':=\tilde \Omega\cap\big({\bigcap}_{M\in \N}\Omega^M\big)$; then $\widetilde\Omega'\in \mathcal T$, $\mathbb{P}(\widetilde\Omega')=1$, and $\widetilde\Omega'$ depends only on the m.p.p. $(\Phi,\mathcal{R})$, but is independent of all the parameters $\eps_j,\theta,M$ and the function $u$. For every $\omega\in \widetilde\Omega'$, by taking the limit (superior) in \eqref{for_whole_process_1st_estimate} as $j\to +\infty$, $\theta\to 0^+$ and then $M\to +\infty$ in this order, and by \eqref{splitt}, \eqref{I_j_M_estimate}, \eqref{J_j_M_estimate}, \eqref{c:cap}, and \eqref{1st_law_of_large_numbers_type}, we arrive at 
$$
\underset{M\to+\infty}{\overline{\lim}}\ \underset{\theta\to 0^+}{\overline{\lim}}\ \underset{j\to+\infty}{\overline{\lim}}
\eps_j^n\!\!\sum_{x^\omega_{j,i}\in  \Phi^{2/M,\omega}_{\eps_j}(D) \setminus  G_{j,M}^{\omega}} \fint_{B_{\theta\eps_{j}/M}(\eps_jx^\omega_{ j,i})}\varphi_{\theta,\rho_{j,i,M}^{ \omega}}(u)\,\mathrm{d}x=0\,,
$$
which proves \eqref{limits_in_the whole_process}.

\smallskip

\textit{Step 2: Proof of a simplified version of \eqref{reduced_limit_to_show} by \eqref{limits_in_the whole_process}.}  In view of \eqref{limits_in_the whole_process}, to prove \eqref{reduced_limit_to_show} it suffices to show that there exists $\Omega'\in \mathcal{T}$ with $\mathbb{P}(\Omega')=1$ and subsequences in $j\in \N$ and $\theta\in (0,1)$ (not relabelled), such that for all $\omega\in \Omega'$, and every $u\in W_0^{1,q}(D;\R^m)\cap L^\infty(D;\R^m)$
\begin{equation}\label{whole_process_reduced_limit_to_show}  
\lim_{M\to +\infty}\lim_{\theta\to 0^+}\lim_{j\to +\infty} \eps_j^n\sum_{x^\omega_{ j,i}\in  \Phi^{2/M,\omega}_{\eps_j}(D)}\fint_{B_{\theta\eps_{j}/M}(\eps_jx^\omega_{ j,i})}\varphi_{\theta,\rho_{j,i,M}^{\omega }}(u)\,\mathrm{d}x=\langle N(Q)\rangle\int_{D}\varphi(u)\,\mathrm{d}x\,.
\end{equation}	 
First of all, by a standard density argument (along the very same lines as the arguments between \eqref{modified_version_at_eps_radius}--\eqref{u_k_to_u_2}, using the bounds \eqref{phi_theta_rho__bounds}-- \eqref{tilde_phi_rho_Lipschitz_estimate}), and by the Dominated Convergence Theorem we may suppose without loss of generality that $u\in C_{c}^{\infty}(D;\R^m)$. 

To prove \eqref{whole_process_reduced_limit_to_show} we apply Proposition \ref{probabilistic_lemma_3} to the m.p.p. $(\Phi, \mathcal{X}^{u})$, for the space-dependent marks 
\begin{equation}\label{marks_for_the_application}
X_{j,i}^{u}(x,\omega):=\varphi_{\theta,\rho^{\omega}_{j,i,M}}(u(x))-\varphi_{\theta,\rho^{\omega}_{j,i,M}}(0)\,,
\end{equation}
and for the radii $r_{\eps_j}:=\theta\eps_j/M$, which satisfy \eqref{r_eps_less_than_eps} for fixed $M\in \N$, $\theta\in(0,1)$. 
We first check that the marks satisfy the assumptions in the proposition. Comparing \eqref{marks_for_the_application} with \eqref{appendix_marks_for_the_application}, we have that in our case 
\begin{equation}\label{matching-k}
\kappa(z,y) := \varphi_{\theta,y}(z)-\varphi_{\theta,y}(0), \quad \text{for } (z,y)\in \R^m\times \R_+,
\end{equation}
and that $Y_{j,i}^\omega=\rho^{\omega}_{j,i,M}$. Due to the truncation to $M$ the marks $Y_{j,i}^\omega$ satisfy \eqref{first_moment_condition}, since $h\in L^1(\R_+; \R_+)$. Moreover, due to the truncation to $M$ they also satisfy \eqref{additional-trunc}. We now check that the function $\kappa$ in \eqref{matching-k} satisfies the requirements for Proposition \ref{probabilistic_lemma_3}. First of all, $\kappa$ is Carath\'eodory and hence measurable, since by Corollary \ref{local_uniform_convergence_for_capacity_densities} the function $z\mapsto \varphi_{\theta,y}(z)$ is continuous, while $y\mapsto \varphi_{\theta,y}(z)$ is increasing, and hence measurable. Moreover, $\kappa(0,y)=0$ and $\kappa$ is bounded from below by \eqref{phi_theta_rho__bounds}, which implies that $\kappa(z,y)\geq -(C_2+C_4)$ for every $y\in [0,M]$ and $\theta\in (0,1)$. By \eqref{phi_theta_rho__bounds} it is immediate to see that \eqref{c:k-gc} holds true with $C_{\kappa}=\max\{C_2+C_4,C_3\}$ and $r:=n-q$, for every $\theta\in (0,1)$;  also, \eqref{phi_theta_rho_Lipschitz_estimate} guarantees the validity of \eqref{c:k-lip}, again for every $\theta\in (0,1)$. Finally, by \eqref{c:k-gc} we have that, for every $x\in \R^n$, 
\begin{align*}
|\langle X^{u}(x,\cdot)\rangle|&=\int_{0}^{+\infty}|\kappa(u(x),y)|h(y)\,\mathrm{d}y \leq C_\kappa |u(x)|^q\int _{0}^{+\infty}y ^{n-q} h(y)\,\mathrm{d}y + C_\kappa\\
&\leq C_\kappa |u(x)|^q M^{n-q} + C_\kappa,
\end{align*} 
where we have used that $\int_{\R_+}h=1$. Hence 
\eqref{new_first_moment_condition} is also satisfied. 

By applying \eqref{strong_integral_law_of_numbers} to the m.p.p. $(\Phi, \mathcal{X}^{u})$ and \eqref{1st_law_of_large_numbers_type-ii} to the m.p.p. $(\Phi^{2/M},\varphi_{\theta,\rho\wedge M}(0))$, we have that for fixed $M\in \N$ and $\theta\in (0,1)$ (with $\frac{1}{\theta}\in \N$), there exists $\Omega_{M,\theta} \in \mathcal T$ with $\mathbb P(\Omega_{M,\theta})=1$ and a subsequence in $j\in \N$ (not relabelled) such that for every $\omega\in \Omega_{M,\theta}$, 
\begin{eqnarray*}
\lim_{j\to +\infty}\eps_j^n\sum_{x^\omega_{j,i}\in \Phi^{2/M,\omega}_{\eps_j}(D)}\Big(\fint_{B_{\theta\eps_{j}/M}(\eps_jx^\omega_{ j,i})}(\varphi_{\theta,\rho_{j,i,M}^{ \omega}}(u)-\varphi_{\theta,\rho^\omega_{j,i,M}}(0))\,\mathrm{d}x\Big)\\
=\langle N^{2/M}(Q) \rangle \int_{D}\big(\langle \varphi_{\theta,\rho\wedge M}(u)\rangle-\langle\varphi_{\theta,\rho\wedge M}(0)\rangle \big)\,\mathrm{d}x\,,
\end{eqnarray*}
and
\begin{equation*}
\lim_{j\to +\infty}\eps_j^n\sum_{x^\omega_{j,i}\in \Phi^{2/M,\omega}_{\eps_j}(D)}\varphi_{\theta,\rho^\omega_{j,i,M}}(0)=
\langle N^{2/M}(Q) \rangle \L^n(D)\langle\varphi_{\theta,\rho\wedge M}(0)\rangle\,,
\end{equation*}
which gives
\begin{equation}\label{final_application_of_integral_lln}
\lim_{j\to +\infty}\eps_j^n\sum_{x^\omega_{j,i}\in \Phi^{2/M,\omega}_{\eps_j}(D)}\fint_{B_{\theta\eps_{j}/M}(\eps_jx^\omega_{ j,i})}\varphi_{\theta,\rho_{j,i,M}^{ \omega}}(u(x))\,\mathrm{d}x= \langle N^{2/M}(Q) \rangle \int_{D}\langle \varphi_{\theta,\rho\wedge M}(u)\rangle(x)\,\mathrm{d}x\,.
\end{equation}
Finally, we set 
\[\Omega':=\bigcap_{M,\frac{1}{\theta}\in\N}\Omega_{M,\theta}\,,\] 
which by its definition only depends on the m.p.p. $(\Phi,\mathcal{R})$, and satisfies $\mathbb{P}(\Omega')=1$. Then for every $\omega\in \Omega'$ 
using \eqref{theta_local_uniform_convergence_phi_j}, Proposition  \ref{identification_of_capacity_density}, the fact that for $\L^n$-a.e. $x\in D$ and $\L^1$-a.e. $\rho>0$, 
$$\varphi_{\rho\wedge M}(u(x))\rightarrow \varphi_{\rho}(u(x)) \ \text{ as } M\to +\infty\,,$$ 
and passing to the limit as $\theta\to 0^+$ and then $M\to +\infty$ in \eqref{final_application_of_integral_lln}, by the Dominated Convergence Theorem and \eqref{limit_of_thinnings} we obtain \eqref{whole_process_reduced_limit_to_show} and conclude the proof.
\end{proof}

\begin{remark}
\normalfont In the periodic setting, an important assumption for the validity of the analogue of Proposition \ref{main_prop_for_approximating_the_capacity}, i.e., \cite[Proposition 4.3]{Ansini-Braides}, was that the corresponding radii for the application of the joining lemma therein were all equal to a constant multiple of the lattice spacing (with constant less that $1/2$) (cf. the Erratum \cite{Erratum-Ansini-Braides}). Such a fact is also reflected here in the statement of Proposition \ref{main_prop_for_approximating_the_capacity}, where the radii of all the perforations are equal to $\theta\eps_j/M$, $M\in \N$, $\theta \in(0,1)$, for all $ x^\omega_{j,i}\in G_{j,M}^\omega$.

\end{remark}

\section{Proof of Theorem \ref{main_thm}}\label{s:main-thm-proof}

\subsection{The $\Gamma$-liminf inequality}
In this subsection we prove the following result.

\begin{proposition}\label{newprop_inf}
Let $(\eps_j) \searrow 0$. There exists $\Omega'\in\mathcal{T}$ with $\mathbb{P}(\Omega')=1$, and a deterministic subsequence of $(\eps_j)$ (not relabelled), such that for every $\omega\in \Omega'$ and every $(u_j), u \in W_0^{1,q}(D;\R^m)$ satisfying $u_j\rightarrow u$ in $L^{1}(D;\R^m)$, we have 
\begin{equation}\label{liminf_inequality_to_show}
\liminf_{j\to +\infty}\mathcal{F}_{\eps_j}^{\omega}(u_j)\geq \mathcal{F}_0(u)\,,
\end{equation}
where 
$\mathcal{F}_{\eps_j}^{\omega}$ and $ \mathcal{F}_0$ are defined in \eqref{random_functionals} and \eqref{deterministic_functionals} respectively. 
\end{proposition}

\begin{proof}
To prove \eqref{liminf_inequality_to_show} we follow closely the arguments in the periodic case (cf. \cite[Section 5]{Ansini-Braides}, and \cite{Erratum-Ansini-Braides}), and adapt them to our stochastic setting.

Without loss of generality, we assume that $\liminf_{j\to +\infty}\mathcal{F}_{\eps_j}^{\omega}(u_j)<+\infty$. 
Then 
$$u_j\equiv0\quad \text{on } (H_{\eps_j}^{\omega}\cap D)\cup \partial D\,,
$$
where $H^\omega_{\eps_j}$ is defined as in \eqref{random_holes}, and by \eqref{growth_of_f}
\begin{equation}\label{weak_convergence_for_liminf}
u_j\rightharpoonup u \quad \text{weakly in } W^{1,q}(D;\R^m) \text{ as } j\to +\infty\,.
\end{equation}

\textit{Step 1: Truncating $(u_j)$ and $u$.} \EEE
Using \cite[Lemma 3.5]{vitalibraides1996homogenization}, for every $L\in \N$ and $\eta>0$, there exists $R_L\geq L$ and a Lipschitz function $\Psi_L:\R^m\mapsto \R^m$ with Lipschitz constant at most $1$, satisfying 
\begin{equation*}
\Psi_L(z):=\begin{cases} z \quad \text{if }\ |z|<R_L\,,\\
0 \quad \text{if }\ |z|>2R_L\,,
\end{cases}
\end{equation*}	 
and such that for every $\omega \in \Omega$, 
\begin{equation}\label{eta_error_in_energy}
\liminf_{j\to+\infty}\mathcal{F}_{\eps_j}^{\omega}(u_j)\geq \liminf_{j\to+\infty}\mathcal{F}_{\eps_j}^{ \omega}(\Psi_L(u_j))-\eta\,.
\end{equation}	
Moreover, setting 
\begin{equation}\label{truncated_functions}
u_j^L:=\Psi_L(u_j)\,, \quad u^L:=\Psi_L(u)\,,
\end{equation}
we also have that
\begin{align}\label{properties_of_truncations}
\begin{split}
(i)&\quad u_j^L\equiv 0 \ \text{on } (H_{\eps_j}^{\omega}\cap D)\cup \partial D\,,\quad \quad  \sup_{j\in \N}\|u_j^L\|_{L^\infty(D)}\leq 2R_L\,,\\
(ii)&\quad u_j^L\underset{j\to +\infty}{\rightharpoonup} u^L\,\ \text{ and }\ u^L\underset{L\to +\infty}{\rightharpoonup } u\ \ \text{weakly in } W^{1,q}(D;\R^m)\,.
\end{split}
\end{align}

\textit{Step 2: Applying the Probabilistic Joining Lemma to $(u_j^L)$.} 
Let now $M, k\in \N$, $\theta\in (0,1)$ (with $\frac{1}{\theta}\in \N$) be fixed, and take $j\in \N$ large enough so that 
\begin{equation}\label{k_M_large_j}
\theta K_{j}>2^kM^2\,,
\end{equation}
where $K_j$ is as in \eqref{K_j_radii}. For $\omega\in \Omega$, let $G_{j,M}^\omega$ be the set defined in \eqref{M_good_centers}, with $\eps$ replaced by $\e_j$. By Lemma \ref{l:probabilistic_joining_lemma} there exists $\Omega'\in \mathcal{T}$ with $\mathbb{P}(\Omega')=1$ such that for every $\omega\in\Omega'$ and every $x^\omega_{j,i}\in G_{j,M}^{\omega}$ there exists $k_{j,i}^{\omega}\in \{0,\dots, k-1\}$ and corresponding annuli $C^\omega_{j,i}:=C_{\eps_j,\theta, M}^{k^\omega_{j,i}}(\e_j x^\omega_{j,i})$, such that 
we can construct a sequence $(w^\omega_j)\subset W_0^{1,q}(D;\R^m)$ satisfying properties \ref{properties_of_the_modification_w_j_M-1-omega} -\ref{c:jl_estimate-omega}, with ${\bar u}_{j,i}^{L,\omega}$ and ${\bar \sigma}_{j,i}^\omega$ defined as in \ref{properties_of_the_modification_w_j_M-2-omega} therein.

Note that the modified sequence $(w^\omega_j)$ vanishes on the perforations centred in $G_{j,M}^\omega$, and on $\partial D$. Indeed, by Lemma \ref{l:probabilistic_joining_lemma} \ref{properties_of_the_modification_w_j_M-1-omega}, we have that 
\begin{equation}\label{i-now-with-L}
w^\omega_j\equiv u^L_j \ \text{ in } D\setminus \bigcup_{x^\omega_{j,i} \in {G}^\omega_{j,M}} C_{j,i}^{\omega}\,.
\end{equation}
Moreover, by \eqref{eps_scaling}, \eqref{K_j_radii}, \eqref{k_M_large_j}, for $j\in \N$ large enough (with respect to $k,M,\theta$) and for every $x^\omega_{j,i}\in G_{j,M}^\omega$ we have
$$\alpha_{\eps_j}\rho_{j,i}^\omega\leq (\eps_{j}/K_j) M\leq 2^{-k}\theta\eps_j/M\,,$$
and hence
\begin{equation*}
B_{\alpha_{\eps_j}\rho^{\omega}_{j,i}}(\eps_jx^\omega_{ j,i})\subset B_{2^{-k}\theta\eps_j/M}(\eps_jx^\omega_{ j,i})\, \ \ \forall \, x^\omega_{j,i}\in G_{j,M}^{\omega}\,.
\end{equation*}
It then follows that 
$$
\bigcup_{x^\omega_{j,i} \in {G}^\omega_{j,M}} \overline{B}_{\alpha_{\eps_j}\rho^{\omega}_{j,i}}(\eps_jx^\omega_{ j,i}) \subset D\setminus \bigcup_{x^\omega_{j,i} \in {G}^\omega_{j,M}} C_{j,i}^{\omega}
$$
and so by \eqref{i-now-with-L} and by \eqref{properties_of_truncations}$(i)$, 
\begin{equation}\label{modifications_vanish_on_the_holes}
w_{j}^{\omega}\equiv u_j^L\equiv 0\,\ \text{on } \bigg(\bigcup_{x^\omega_{j,i} \in {G}^\omega_{j,M}} B_{\alpha_{\eps_j}\rho^{\omega}_{j,i}}(\eps_jx^\omega_{ j,i})\cap D\bigg)\cup \partial D\,.
\end{equation}
Moreover, since $w_j^\omega\underset{j\to \infty}{\rightharpoonup} u$ weakly in $W^{1,q}(D;\R^m)$,
\begin{equation}\label{uniform_W_1_q_bound_ind_of_omega-w}
\sup_{j\in \N}\|w_{j}^{\omega}\|_{W^{1,q}(D;\R^m)} \leq C<+\infty\,,
\end{equation} 
where the constant $C>0$ might depend on $L, u$ but not on $\omega\in \Omega'$.

\smallskip

\textit{Step 3: Splitting of the energy.}  By Lemma \ref{l:probabilistic_joining_lemma} \ref{c:jl_estimate-omega} we have
\begin{align}\label{en-split:wj}
\mathcal{F}_{\eps_j}^{\omega}(u^L_j)= \int_{D}f(\nabla u_{j}^{L})\,\mathrm{d}x\geq \int_{D\setminus E_j^{\theta,M, \omega}}f(\nabla w_{j}^{\omega})\,\mathrm{d}x+\int_{E_{j}^{\theta,M, \omega}}f(\nabla w_{j}^{\omega})\,\mathrm{d}x-\frac{C}{k},
\end{align}
where we set 
\begin{equation*}
E_{j}^{\theta,M,\omega}:=\bigcup_{x^\omega_{ j,i}\in G_{j,M}^{\omega}}B_{{\bar \sigma}_{j,i}^{\omega}}(\eps_jx^\omega_{ j,i})\,.
\end{equation*}
In what follows we deal with the two integrals in the right-hand side of \eqref{en-split:wj} separately. We will show that the integral outside $E_{j}^{\theta,M,\omega}$ will result in the first integral in $\mathcal{F}_0$, while the integral on $E_{j}^{\theta,M,\omega}$ will give the capacitary term in the limit.

\smallskip

\textit{Step 3.1: The energy contribution outside $E_{j}^{\theta,M,\omega}$.} For every $\omega\in \Omega'$  we define the auxiliary sequence
\begin{equation*}
v_{j}^{\omega}:=\begin{cases}
w_{j}^{\omega}\quad &\text{in }\ D\setminus E_j^{\theta,M,\omega}\,,\\[4pt] 
{\bar u}_{j,i}^{L,\omega} \quad &\text{in } \overline{B}_{{\bar \sigma_{j,i}^\omega}}\ \ \text{for each } x^\omega_{ j,i}\in G_{j,M}^{\omega}\,.
\end{cases}
\end{equation*}	 
Note that $v_{j}^{\omega}$ is a modification of $w_j^\omega$ in $D$ at no additional cost in terms of $f$, since $v_{j}^{\omega}$ is piecewise constant on $E_j^{\theta,M, \omega}$ and $f(0)=0$ by assumption, so that
\begin{equation}\label{aaaab}
\int_{D\setminus E_j^{\theta,M, \omega}}f(\nabla w_{j}^{\omega})\,\mathrm{d}x
= \int_{D}f(\nabla v_{j}^{\omega})\,\mathrm{d}x\,.
\end{equation}
Moreover, by construction the sequence $(v_{j}^{\omega})$ is bounded in $W^{1,q}(D;\R^m)$ independently of $\omega \in \Omega'$, since by \eqref{properties_of_truncations} and \eqref{uniform_W_1_q_bound_ind_of_omega-w},
\begin{equation}\label{uniform_W_1_q_bound_ind_of_omega}
\sup_{j\in \N}\|v_{j}^{\omega}\|_{W^{1,q}(D;\R^m)}\lesssim R_L+\sup_{j\in \N}\|w_{j}^{\omega}\|_{W^{1,q}(D;\R^m)} \leq C<+\infty\,.
\end{equation} 
We now show that there exists a set $\Omega'\in \mathcal{T}$ with $\mathbb{P}(\Omega')=1$ (which may be different from the one in Step 2, but depends only on the m.p.p. $(\Phi,\mathcal{R})$, and hence will not be relabelled) such that for every $\omega\in \Omega'$ 
\begin{equation}\label{v_j_M_L_omega_convergence_is_ok}
v_{j}^{\omega}\rightharpoonup u^L\, \quad \text{weakly in } W^{1,q}(D;\R^m) \text{ as } j\to + \infty\,.
\end{equation}
Note that any subsequence $(v_j^{\omega})$ is bounded in $W^{1,q}(D;\R^m)$ by \eqref{uniform_W_1_q_bound_ind_of_omega}, and hence admits a (not-relabelled)  convergent subsequence, converging weakly in $W^{1,q}(D;\R^m)$ to a function $v^\omega$, that a priori might be probabilistic. It therefore suffices to show that  there exists $\Omega'\in \mathcal{T}$ with $\mathbb{P}(\Omega')=1$ such that for every $\omega\in \Omega'$, 
\begin{equation*}
v^\omega\equiv u^L \quad \L^n-\mathrm{a.e.\ in }\ D\,.
\end{equation*}
We have that
$$
v_j^\omega=w_j^\omega=u_j^L \quad \text{in }\ 
D\setminus \underset{x^\omega_{j,i}\in \Phi_{j}^{2/M,\omega}(D)}{\bigcup}B_{\theta\eps_j/M}(\eps_j x^\omega_{j,i})\subset D\setminus \underset{x^\omega_{j,i}\in G_{j,M}^\omega}{\bigcup}B_{\theta\eps_j/M}(\eps_j x^\omega_{j,i})\,.
$$
In particular, this means that 
\begin{equation}\label{chi_j_weakly_converging_to_constant}
(v_j^\omega-u_j^L)\chi_j^\omega=0, \quad \text{where } \quad \chi_j^\omega:=\chi_{D\setminus \underset{x^\omega_{j,i}\in \Phi_{j}^{2/M,\omega}(D)}{\bigcup}B_{\theta\eps_j/M}(\eps_j x^\omega_{j,i})}\,.
\end{equation}
By the stationarity assumption \ref{H1} and \ref{H3} of Subsection \ref{Assumptions_on_m.p.p.} for $\Phi^{2/M}$ (these properties being transferred to it by the corresponding ones of $\Phi$), we have that  $\Phi^{2/M}$ is ergodic. Hence, Birkhoff's Ergodic Theorem (see, e.g.,  \cite[Property 2.10]{PelScaZep}) guarantees that, since $\theta\in (0,1)$, there exists 
a set $\Omega'\in \mathcal{T}$ with $\PP(\Omega')=1$ such that for every $\omega\in \Omega'$ 
\begin{equation}\label{weak_limit_of_chi_j}
\chi_j^\omega\rightharpoonup K>0\, \quad \text{weakly* in }L^\infty\,, 
\end{equation} 
where $K>0$ is a deterministic constant. From \eqref{chi_j_weakly_converging_to_constant}, \eqref{weak_limit_of_chi_j}, \eqref{properties_of_truncations}($ii$) and by the convergence (up to subsequences) of $v_j^\omega$ to $v^\omega$, we then conclude that 
$$
0=\int_{D} \chi_j^\omega |v_j^\omega-u_j^L| \,\mathrm{d}x 
\underset{j\to +\infty}{\to} K \int_{D} |v^\omega-u^L| \,\mathrm{d}x\,,
$$
and hence $v^\omega\equiv u^L$ $\mathcal{L}^n$-a.e. in $D$, as desired. This proves \eqref{v_j_M_L_omega_convergence_is_ok}.

Then, by taking the $\liminf$ as $j\to +\infty$ in \eqref{aaaab}, we obtain that for every $\omega$ in a set $\Omega'\in \mathcal{T}$ with $\mathbb{P}(\Omega')=1$ (which depends only on the m.p.p. $(\Phi,\mathcal{R})$, but is independent of all the parameters and the functions involved in the arguments), 
\begin{align}\label{first_energy_estimate_from_below}	
\liminf_{j\to+\infty}\int_{D\setminus E_j^{\theta,M, \omega} } f(\nabla w_{j}^{\omega})\,\mathrm{d}x=\liminf_{j\to+\infty}\int_D f(\nabla v_{j}^{\omega})\,\mathrm{d}x
\geq\int_{D}Qf(\nabla u^L)\,\mathrm{d}x\,,
\end{align}
where we have used \eqref{v_j_M_L_omega_convergence_is_ok}, and the fact that  the functional $\int_D Qf(\nabla \cdot)\,\mathrm{d}x$ is the lower semicontinuous envelope of $\int_D f(\nabla \cdot)\,\mathrm{d}x$ with respect to the weak $W^{1,q}(D;\R^m)$-topology (cf. \cite{Ansini-Braides}).

\smallskip

\textit{Step 3.2: The energy contribution in $E_{j}^{\theta,M, \omega}$.}  Let $j\in \N$ be large enough so that \eqref{k_M_large_j} is satisfied, where $K_j$ is as in \eqref{K_j_radii}, and let $x^\omega_{j,i}\in G_{j,M}^{\omega}$. We define $\zeta_{j,i}^{\omega}: B_{\theta K_{j}}(0)\to \R^m$ as 
\begin{equation}\label{suitable_test_function}
\zeta_{j,i}^{\omega}(y):=
\begin{cases}
\medskip
\displaystyle w_{j}^{\omega}(\eps_jx^\omega_{ j,i}+\alpha_{\eps_j}y)\ \ &\text{in } B_{\bar\sigma^\omega_{j,i}/\alpha_{\eps_j}}(0)\,,\\
\bar u_{j,i}^{L,\omega} &\text{in } B_{\theta K_{j}}(0)\setminus B_{\bar\sigma^\omega_{j,i}/\alpha_{\eps_j}}(0)\,,
\end{cases}
\end{equation}
where $\bar\sigma^\omega_{j,i}$ and $\bar u_{j,i}^{L,\omega}$ are defined as in Lemma \ref{l:probabilistic_joining_lemma}$(ii)$. In particular, \EEE
\begin{equation}\label{zeta_on_the_boundary}
\zeta_{j,i}^{\omega}|_{\partial B_{\bar\sigma^\omega_{j,i}/\alpha_{\eps_j}}(0)}\equiv w_{j}^{\omega}|_{\partial B_{\bar \sigma_{j,i}^{\omega}}(\eps_jx^\omega_{j,i})}\equiv \bar u_{j,i}^{L,\omega}\,,
\end{equation}
and
$B_{\bar\sigma^\omega_{j,i}/\alpha_{\eps_j}}(0)=B_{\tfrac{3}{4}2^{-k^\omega_{j,i}}\theta K_{j}/M}(0)\subset B_{\theta K_{j}}(0)$. Moreover, since $x^\omega_{ j,i}\in G_{j,M}^{\omega}$, by \eqref{k_M_large_j} we have that 
$$
\alpha_{\eps_j}\rho_{j,i}^\omega \leq \alpha_{\eps_j}M \leq 
\frac{\eps_j}{K_j} M <  \eps_j M \frac{\theta 2^{-k}}{M^2} 
\leq \frac{3}{4}2^{-(k-1)}\frac{\theta \eps_j}{M}\leq \bar\sigma_{j,i}^\omega\,.
$$
By \eqref{suitable_test_function}, \eqref{zeta_on_the_boundary}, and \eqref{modifications_vanish_on_the_holes}, we have
\begin{equation}\label{properties_of_the_test_function_z}
\zeta_{j,i}^\omega-\bar u_{j,i}^{L,\omega}\in W_0^{1,q}(B_{\theta K_{j}}(0);\R^m)\,, \quad \zeta_{j,i}^{\omega}|_{\overline{B}_{\rho_{j,i}^{ \omega }}(0)}\equiv w_{j}^{\omega}\big|_{\overline{B}_{\alpha_{\eps_j}\rho_{j,i}^{ \omega }}(\eps_jx^\omega_{ j,i })}\equiv 0\,,
\end{equation}
i.e., $\zeta_{j,i}^{\omega }$ is a competitor in the minimisation problem defining $\varphi^j_{\theta,\rho_{j,i}^{\omega}}(\bar u_{j,i}^{L,\omega})$ (see   \eqref{K_j_rho_truncation_of_capacity}). 

By a change of variables we rewrite the energy contribution in \eqref{en-split:wj} relative to the set $E_j^{\theta,M, \omega}$, as
\begin{align*}
\int_{E_j^{\theta,M, \omega }}f(\nabla w_{j}^{\omega})\,\mathrm{d}x&=\sum_{x^\omega_{j,i}\in G_{j,M}^{\omega}}\int_{B_{\bar{\sigma}_{j,i}^\omega}(\eps_jx^\omega_{ j,i})} f(\nabla w_{j}^{\omega} (x))\,\mathrm{d}x\\
&=\alpha_{\eps_j}^n\sum_{x^\omega_{ j,i}\in G_{j,M}^{\omega}}\int_{B_{\bar{\sigma}_{j,i}^\omega/\alpha_{\eps_j}}(0)} f\big(\nabla w_{j}^{\omega}(\eps_jx^\omega_{ j,i}+\alpha_{\eps_j}y)\big)\,\mathrm{d}y\\
&=\eps_j^{n}\sum_{x^\omega_{j,i}\in G_{j,M}^\omega}\int_{B_{\bar{\sigma}_{j,i}^\omega/\alpha_{\eps_j}}(0)} \alpha_{\eps_j}^{q}f\big(\alpha_{\eps_j}^{-1}\nabla \zeta_{j,i}^{ \omega }(y)\big)\,\mathrm{d}y\,.
\end{align*}
Recalling that $f(0)=0$, since $\nabla\zeta_{j,i}^{\omega}\equiv 0$ in $B_{\theta K_{j}}(0)\setminus B_{\bar{\sigma}_{j,i}^\omega/\alpha_{\eps_j}}(0)$, for every $x^\omega_{j,i}\in G_{j,M}^{\omega}$ we have 
\begin{align*}
\int_{B_{\bar{\sigma}_{j,i}^\omega/\alpha_{\eps_j}}(0)} \alpha_{\eps_j}^{q}f\big(\alpha_{\eps_j}^{-1}\nabla \zeta_{j,i}^{ \omega }(y)\big)\,\mathrm{d}y &=
\int_{B_{\theta K_{j}}(0)} \alpha_{\eps_j}^{q}f\big(\alpha_{\eps_j}^{-1}\nabla \zeta_{j,i}^{ \omega}(y)\big)\,\mathrm{d}y\\
& \geq \int_{B_{\theta K_{j}}(0)} g_j(\nabla \zeta_{j,i}^{\omega}(y))\,\mathrm{d}y
\geq \varphi^j_{\theta,\rho_{j,i}^{ \omega }}(\bar u_{j,i}^{L,\omega})\,,
\end{align*}
where we have used that $f\geq Qf$, \eqref{j_capacity_densities}, \eqref{K_j_rho_truncation_of_capacity}, and \eqref{properties_of_the_test_function_z}. In conclusion,
\begin{equation}\label{capacity_energy_main_estimate_from_below}
\int_{E_j^{\theta,M, \omega }}f(\nabla w_{j}^{\omega})\,\mathrm{d}x 
\geq  \eps_j^{n}\sum_{x^\omega_{j,i}\in G_{j,M}^\omega}\varphi^j_{\theta,\rho_{j,i}^{ \omega }}(\bar u_{j,i}^{L,\omega})\,.
\end{equation}
Taking now in \eqref{capacity_energy_main_estimate_from_below} the liminf as $j\to +\infty$, then (up to a not-relabelled subsequence) $\theta \to 0^+$ and then $M\to +\infty$, and using \eqref{properties_of_truncations} and Proposition \ref{main_prop_for_approximating_the_capacity}, we deduce that there exists $\Omega'\in\mathcal{T}$ with $\mathbb{P}(\Omega')=1$ (depending only on $(\Phi,\mathcal{R})$) such that for every $\omega\in\Omega'$, 
\begin{equation}\label{liminf_capacity_inequality}
\lim_{M\to +\infty}\lim_{\theta\to 0^+}\liminf_{j\to +\infty}\int_{E_j^{\theta,M, \omega }}f(\nabla w_{j}^{\omega})\,\mathrm{d}x\geq \langle N(Q)\rangle \int_D\varphi(u^L)\, \mathrm{d}x\,.
\end{equation} 

\smallskip

\textit{Step 4: Conclusion.}  Taking in \eqref{en-split:wj} the liminf as $j\to +\infty$,
and then the (subsequential) limits as $\theta \to 0^+$ and $M\to +\infty$, and since the left-hand side of \eqref{en-split:wj} is independent of $\theta,M$, in view of \eqref{first_energy_estimate_from_below} and \eqref{liminf_capacity_inequality} we obtain that for $\mathbb{P}$-a.e. $\omega\in \Omega$, 
\begin{align}\label{almost-almost}
\liminf_{j\to +\infty}\mathcal{F}_{\eps_j}^{\omega}(u^L_j)\geq \int_{D}Qf(\nabla u^L)\,\mathrm{d}x+\langle N(Q)\rangle \int_D\varphi(u^L)-\frac{C}{k}.
\end{align}
Combining \eqref{almost-almost}, \eqref{eta_error_in_energy} and \eqref{truncated_functions} we conclude that for $\mathbb{P}$-a.e. $\omega\in \Omega$ 
\begin{equation}\label{liminf_almost_done}
\liminf_{j\to +\infty}\mathcal{F}_{\eps_j}^{\omega }(u_j)\geq\int_{D}Qf(\nabla u^L)\,\mathrm{d}x+\langle N(Q)\rangle \int_D\varphi(u^L)\,\mathrm{d}x-\frac{C}{k}-\eta\,.
\end{equation}
By the arbitrariness of $k,\eta,L>0$ we first let $k\to +\infty$ in \eqref{liminf_almost_done}, then $\eta\to 0^+$, and finally $L\to +\infty$. Then,  \eqref{properties_of_truncations}$(ii)$, the lower semicontinuity of the functional $\int_D Qf(\nabla\cdot)\,\mathrm{d}x$ with respect to the weak $W^{1,q}(D;\R^m)$-topology and the continuity of $\int_D\varphi(\cdot)\,\mathrm{d}x$ with respect to strong $L^q(D;\R^m)$-topology, guarantee the existence of a set $\Omega'\in\mathcal{T}$ with $\mathbb{P}(\Omega')=1$ (depending only on $(\Phi,\mathcal{R})$) such that for every $\omega\in\Omega'$ and every $(u_j), u\in W^{1,q}_0(D;\R^m)$ with $u_j\underset{j\to+\infty}{\rightharpoonup}u$ weakly in $W^{1,q}(D;\R^m)$, \EEE  
$$\liminf_{j\to +\infty}\mathcal{F}_j^{ \omega }(u_j)\geq \int_{D}Qf(\nabla u)\,\mathrm{d}x+\langle N(Q)\rangle \int_D\varphi(u)\,\mathrm{d}x=\mathcal{F}_0(u)\,,$$
which concludes the proof of the proposition.
\end{proof}

\subsection{The $\Gamma$-limsup inequality}
In this subsection we prove the following result.

\begin{proposition}\label{Gamma_limsup_ineq}
Let $(\eps_j) \searrow 0$. Let $\Omega'\in\mathcal{T}$ with $\mathbb{P}(\Omega')=1$ be such that the conclusions of Proposition \ref{newprop_inf} hold true for every $\omega\in \Omega'$, and let $u\in W_0^{1,q}(D;\R^m)$. Then there exists a sequence $(u_j)\subset W_0^{1,q}(D;\R^m)$ satisfying $u_j\rightarrow u$ in $L^{1}(D;\R^m)$, and such that 
\begin{equation}\label{limsup_inequality_to_show}
\limsup_{j\to +\infty}\mathcal{F}_{\eps_j}^{\omega}(u_j)\leq \mathcal{F}_0(u)\,,
\end{equation}
where  
$\mathcal{F}_{\eps_j}^{\omega}$ and $ \mathcal{F}_0$ are defined in \eqref{random_functionals} and \eqref{deterministic_functionals} respectively. 
\end{proposition}

Before embarking on the proof of the proposition, we present some auxiliary lemmata, that can be thought of as partially constructing ``correctors'' in the spirit of \cite[Lemma 3.1]{Giunti-Hofer-Velasquez}. These auxiliary results will be crucial in the construction of an admissible recovery sequence, namely that vanishes on the entire set of perforations $H_\eps^\omega \cap D$.

\begin{figure}
\includegraphics[scale=0.7]{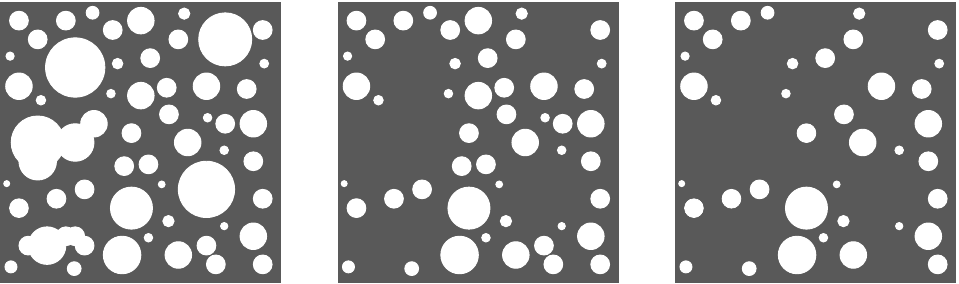}
\caption{From left to right: a portion of the random domain, the holes centred at $G_{\eps,M}^\o$, and the very good holes.}
\label{fig:goodbaadall}
\end{figure}

We first introduce the set of \textit{very good perforations} as follows. Let $\eps>0$, $M\in \N$, $\theta\in (0,1)$ (with $1/\theta\in \N$) be fixed, and assume that 
\begin{equation}\label{even_smaller_theta}
0<\theta<\frac{3}{8M}\,.
\end{equation}
For every $\omega\in\Omega$ we define 
$$
I_{\eps,g,\theta}^\o:=\{ x^\o_i\in I_{\eps,g}^\o \colon \exists \, x^\o_k\in I^\omega_{\eps,g}\setminus \{x_i^\o\} \ \text{with } \partial B_{\theta \eps}(\eps x^\o_i)\cap \overline{B}_{\eps r_\eps^\omega}(\eps x^\o_k)\neq \emptyset\}\,,
$$
where $(r_\eps^\o)$ and $I^\omega_{\eps,g}$ are as in Lemma \ref{good_and_bad_balls_lemma}. 
These are the good centres that are not $(\theta+r_\eps^\o)$-separated from other good centres. We recall that points in  $I^\omega_{\eps,g}$ are at least $2 r_\eps$-separated from one another (cf. \eqref{c:13_marzo}). Recalling \eqref{M_good_centers}, we define the \textit{mildly good centers} as
\begin{equation}\label{eq: mildly_good_centers}
{(MG)}_{\eps,\theta, M}^\omega:=(I^\omega_{\eps,g}\setminus G_{\eps,M}^\omega)\cup (G_{\eps,M}^\omega \cap I_{\eps,g,\theta}^\o)\,,
\end{equation}
and the \textit{very good centers} as
\begin{equation}\label{eq: very_good_centers}
{(VG)}_{\eps,\theta,M}^\omega:=G_{\eps,M}^\omega\setminus {(MG)}_{\eps,\theta, M}^\omega\,.	
\end{equation}
In short, the very good centres are the subset of $G_{\eps,M}^\omega$ for which the corresponding balls are deterministically separated, at scale $\eps$, also from other good balls in $I^\omega_{\eps,g}$, which is a priori not guaranteed (see Remark \ref{exclude_points-close_to_safety}). Figure \ref{fig:goodbaadall} illustrates the identification of the perforations centred at $G^\o_{\eps,M}$, and their subset of very good perforations.

We also use the shorthand notations ${(MG)}_{j,\theta, M}^\omega:={(MG)}_{\eps_j,\theta, M}^\omega$ and ${(VG)}_{j,\theta, M}^\omega:={(VG)}_{\eps_j,\theta, M}^\omega$ for a sequence $(\eps_j)\searrow 0$. 

A fundamental property that we will use is that the mildly good perforations  do not contribute to the limit capacity, as made precise in the following lemma.

\begin{lemma} 
Let $(\eps_j) \searrow 0$. There exists $\Omega' \in \mathcal T$ with $\mathbb P(\Omega')=1$ such that for every $\omega\in \Omega'$ we have (possibly along $\omega$-independent subsequences) 
\begin{equation}\label{eq: mildly_good_zero_capacity}
\lim_{M\to +\infty}\lim_{\theta\to 0^+}\lim_{j\to +\infty}\eps_j^n\sum_{x^\omega_{j,i}\in {(MG)}_{j,\theta,M}^{\omega}}(\rho_{j,i}^\omega)^{n-q}=0\,.
\end{equation}
\end{lemma}

\begin{proof}
We show separately that 
\begin{equation}\label{eq: good_notgreat_zero_capacity}
\lim_{M\to +\infty}\lim_{j\to +\infty}\eps_j^n\sum_{x_{j,i}^\o\in(I_{\eps_j,g}^\omega \setminus G_{j,M}^\omega)} (\rho_{j,i}^\o)^{n-q}=0\,
\end{equation} 
and 
\begin{equation}\label{eq: great_closetogood_zero_capacity}
\lim_{M\to +\infty}\lim_{\theta\to 0^+}\lim_{j\to +\infty} \eps_j^n\sum_{x^\omega_{j,i}\in (G_{j,M}^\omega \cap I_{\eps_j,g,\theta}^\o)}(\rho_{j,i}^\omega)^{n-q}=0\,,
\end{equation}
so that \eqref{eq: mildly_good_zero_capacity} follows by \eqref{eq: mildly_good_centers}.

Let $\Omega'\in \mathcal{T}$ be as in Lemma \ref{good_and_bad_balls_lemma}, and let $M\in \N$, $j\in \N$ and $\theta\in (0,1)$ (with $1/\theta\in \N$) be such that (by \eqref{even_smaller_theta})
\begin{equation}\label{eq: parameters_small}
0<r_{\eps_j}^\omega<\theta<\frac{3}{8M}\,.	
\end{equation} 
To prove \eqref{eq: good_notgreat_zero_capacity}, we first decompose in $\Omega'$,
\begin{align*}
I_{\eps_j,g}^\omega &= (I_{\eps_j,g}^\omega\cap \Phi_{\eps_j}^{2/M,\omega}(D))\cup (I_{\eps_j,g}^\omega\setminus \Phi_{\eps_j}^{2/M,\omega}(D))\\
&= (G_{j,M}^\omega\cup I_{j,M}^\omega\cup J_{j,M}^\omega )\cup (I_{\eps_j,g}^\omega\setminus \Phi_{\eps_j}^{2/M,\omega}(D))\,,
\end{align*}
where $I_{j,M}^\omega$ and $J_{j,M}^\omega$ are defined in \eqref{new_index_sets}, so that we can write 
\begin{equation}\label{rewriting-good-not-great}
I_{\eps_j,g}^\omega \setminus G_{j,M}^\omega = I_{j,M}^\omega\cup J_{j,M}^\omega \cup (I_{\eps_j,g}^\omega\setminus \Phi_{\eps_j}^{2/M,\omega}(D))\,.
\end{equation}
Hence, by \eqref{I_j_M_estimate} and \eqref{J_j_M_estimate}, the claim \eqref{eq: good_notgreat_zero_capacity} follows if we show that 
\begin{equation}\label{corto}
\lim_{M\to +\infty}\lim_{j\to +\infty}\eps_j^n\sum_{x_{j,i}^\o\in I_{\eps_j,g}^\omega\setminus \Phi_{\eps_j}^{2/M,\omega}(D)} (\rho_{j,i}^\o)^{n-q}=0\,.
\end{equation}
Let $R\in \N$ be fixed. Then we can write 
\begin{align}\label{split-R}
&\sum_{x_{j,i}^\o\in I_{\eps_j,g}^\omega \setminus\Phi_{\eps_j}^{2/M,\omega}(D)} (\rho_{j,i}^\o)^{n-q}\nonumber\\
&\qquad \leq 
\sum_{x_{j,i}^\o\in I_{\eps_j,g}^\omega \setminus \Phi_{\eps_j}^{2/M,\omega}(D)} (\rho_{j,i}^\o)^{n-q}\chi_{[0,R]}(\rho_{j,i}^\o)+ \hspace{-0.5em}
\sum_{x_{j,i}^\o\in I_{\eps_j,g}^\omega \setminus \Phi_{\eps_j}^{2/M,\omega}(D)} (\rho_{j,i}^\o)^{n-q}\chi_{(R,+\infty)}(\rho_{j,i}^\o)\nonumber\\
&\qquad \leq R^{n-q} \# 
\big(I_{\eps_j,g}^\omega \setminus \Phi_{\eps_j}^{2/M,\omega}(D)\big) +\sum_{x^\omega_{j,i}\in \Phi_{\eps_j}^{\omega}(D)}(\rho_{j,i}^\omega)^{n-q}\chi_{(R,+\infty)}(\rho_{j,i}^\o). 
\end{align}
Since by \eqref{c:cap} we have that $\eps_j^n \#  I^\omega_{\eps_j,b} \to 0$ as $j\to+\infty$, we deduce that there exists a set $\Omega_M\in \mathcal{T}$, $\Omega_M\subset \Omega'$, with $\mathbb{P}(\Omega_M)=1$, such that 
\begin{equation}\label{a:do}
\lim_{j\to +\infty}\eps_j^n  \# \big(I_{\eps_j,g}^\omega \setminus \Phi_{\eps_j}^{2/M,\omega}(D)\big) = \lim_{j\to +\infty} \eps_j^n \# \big(\Phi_{\eps_j}^\omega(D) \setminus \Phi_{\eps_j}^{2/M,\omega}(D) \big) =(\langle N(Q)\rangle-\langle N^{2/M}(Q)\rangle) \mathcal{L}^n(D)\,,
\end{equation}
for $\o\in \Omega_M$, where we have used \eqref{1st_law_of_large_numbers_type} and \eqref{limit_of_thinnings}. Moreover, by Lemma \ref{probabilistic_lemma_1} applied to $Y_{j,i}^\omega:=(\rho_{j,i}^\omega)^{n-q}\chi_{(R,+\infty)}(\rho_{j,i}^\omega)$ we have that 
\begin{equation}\label{lemma1-tr}
\lim_{j\to+\infty}\ \eps_j^n\hspace{-0.5em}\sum_{x^\omega_{j,i}\in \Phi_{\eps_j}^{\omega}(B)}(\rho_{j,i}^\omega)^{n-q}\chi_{(R,+\infty)}(\rho_{j,i}^\omega)=\langle N(Q)\rangle\langle \rho^{n-q}\chi_{(R,+\infty)}(\rho)\rangle\L^n(D)\,
\end{equation}
in a set of full probability $\Omega_{M,R}$. Hence, 
by \eqref{split-R}, \eqref{a:do} and \eqref{lemma1-tr} it follows that $\mathbb{P}$-a.e. in $\Omega$
\begin{align*}
\varlimsup_{j\to +\infty}\eps_j^n\hspace{-0.5em}\sum_{x_{j,i}^\o\in I_{\eps_j,g}^\omega \setminus \Phi_{\eps_j}^{2/M,\omega}(D)} \hspace{-0.5em} (\rho_{j,i}^\o)^{n-q}&\leq R^{n-q} (\langle N(Q)\rangle-\langle N^{2/M}(Q)\rangle) \mathcal{L}^n(D) \\
&+\langle N(Q)\rangle\langle \rho^{n-q}\chi_{(R,+\infty)}(\rho)\rangle\L^n(D)\,.
\end{align*}
Now, by letting $M\to +\infty$, thanks to \eqref{limit_of_thinnings}, we have 
\begin{align*}
\varlimsup_{M\to +\infty}\varlimsup_{j\to +\infty}\eps_j^n\hspace{-0.5em}\sum_{x_{j,i}^\o\in I_{\eps_j,g}^\omega \setminus \Phi_{\eps_j}^{2/M,\omega}(D)} (\rho_{j,i}^\o)^{n-q}\leq \langle N(Q)\rangle\langle \rho^{n-q}\chi_{(R,+\infty)}(\rho)\rangle\L^n(D)\,.
\end{align*}
Finally, by letting $R\to +\infty$, since $\langle\rho^{n-q}\rangle<+\infty$ by \eqref{finite_average_capacity}, 
we obtain \eqref{eq: good_notgreat_zero_capacity}, for every $\omega$ in a set of probability 1, where the latter can be chosen independently of $M\in \N$ and $\theta\in (0,1)$ (with $1/\theta\in \N$). 
For the proof of \eqref{eq: great_closetogood_zero_capacity}, we claim that by the choice \eqref{eq: parameters_small}
\begin{equation}\label{GG-thin}
G_{j,M}^\omega \cap I_{\eps_j,g,\theta}^\o \subset I_{\eps_j,g}^\o\setminus\Phi_{\eps_j}^{2/M,\omega}(D)\,.
\end{equation}
If \eqref{GG-thin} is proved, by proceeding as above for the proof of \eqref{corto}, we can conclude. So it remains to prove \eqref{GG-thin}.

To this aim, let $x_{j,i}^\omega\in G_{j,M}^\omega \cap I_{\eps_j,g,\theta}^\o$. Then there exists $x_{j,k}^\omega\in I_{\eps_j,g}^\omega\setminus \{x_{j,i}^\omega\}$ (and in fact a posteriori $x_{j,k}^\omega\in I_{\eps_j,g}^\omega\setminus G_{j,M}^\omega$) and $y\in \R^n$ such that
\[y \in\partial B_{\theta \eps_j}(\eps_j x_{j,i}^\omega)\cap \overline{B}_{\eps_j r_{\eps_j}^\omega}(\eps_j x_{j,k}^\omega)\,.\]
This and \eqref{eq: parameters_small} imply that
\begin{align*}
|\eps_j x_{j,i}^\omega-\eps_j x_{j,k}^\omega|&\leq |\eps_j x_{j,i}^\omega-y |+|y -\eps_j x_{j,k}^\omega|\leq (\theta+r_{\eps_j}^\omega)\eps_j<2\theta\eps_j<\frac{3\eps_j}{4M}<\frac{2\eps_j}{M}\,,
\end{align*}
i.e., $x_{j,i}^\omega\in I_{\eps_j,g}^\o\setminus\Phi_{\eps_j}^{2/M,\omega}(D)$\, as we claimed in \eqref{GG-thin}.
\end{proof}

We now present a construction of ``partial correctors'' in the spirit of \cite[Lemma 3.1]{Giunti-Hofer-Velasquez}.

\begin{lemma}\label{construction_of_correctors}
There exists a set $\Omega'\in \mathcal{T}$ with $\mathbb{P}(\Omega')=1$ with the following property. Let $(\eps_j)\searrow 0$, $M\in \N$, and $\theta\in (0,1)$ (with $1/\theta\in \N$) such that \eqref{eq: parameters_small} holds true. For every $\omega\in \Omega'$ there exists a corrector function $\phi^\omega_{j,\theta,M}\in W^{1,q}(D;[0,1])$ satisfying the following properties:
\begin{align}\label{eq:properties_of_correctors}
\begin{split}
(i)\quad \phi^\omega_{j,\theta,M}&\equiv 0 \ \text{ on  }\ H_{\eps_j,b}^\omega\cup \big(\bigcup_{x_{j,i}^\omega\in {(MG)}_{j,\theta,M}^\omega}\overline{B}_{\alpha_{\eps_j}\rho_{j,i}^\omega}(\eps_jx_{j,i}^\omega)\big)\,,\\
(ii)\ \ \phi^\omega_{j,\theta,M}&\equiv 1 \ \text{ on  }\ D\setminus \bigg(D_{\eps_j,b}^\omega \cup \bigcup_{x_{j,i}^\omega\in {(MG)}_{j,\theta,M}^\omega}\overline{B}_{\eps_jr^\omega_{\eps_j}}(\eps_jx_{j,i}^\omega)\bigg)\,,\\
\ (iii)\hspace{2em} \quad&\lim_{M\to +\infty}\lim_{\theta\to 0^+}\lim_{j\to +\infty}\|\phi^\omega_{j,\theta,M}- 1\|_{W^{1,q}(D)}=0\,.
\end{split}
\end{align} 
\end{lemma} 

\begin{proof}
First of all, for every $\omega\in\Omega'\in \mathcal{T}$ as in Lemma \ref{good_and_bad_balls_lemma}, by the definition of the $q$-capacity in \eqref{relative_scalar_q_capacity} there exists $\phi_{0,j}^{\omega}\in W_0^{1,q}(D_{\eps_j,b}^\omega;[0,1])$ such that
\[\phi_{0,j}^{\omega}|_{H_{\eps_j,b}^\omega}\equiv 1\,, \quad \int_{D_{\eps_j,b}^\omega}|\nabla \phi_{0,j}^\omega|^q\,\mathrm{d}x\leq 2\mathrm{Cap}_q(H_{\eps_j,b}^\omega, D_{\eps_j,b}^\omega)\,.\]
By setting $\phi_{1,j}^\omega:=1-\phi_{0,j}^\omega$ (extended to 1 in $D\setminus D_{\eps_j,b}^\omega$) we obtain a function $\phi_{1,j}^{\omega}\in W^{1,q}(D;[0,1])$ such that
\begin{equation}\label{eq: phi_1_j}
\phi_{1,j}^{\omega}|_{H_{\eps_j,b}^\omega}\equiv 0\,, \quad \phi_{1,j}^{\omega}|_{D\setminus D_{\eps_j,b}^\omega}\equiv 1\,, \quad \int_{D}|\nabla \phi_{1,j}^\omega|^q\,\mathrm{d}x\leq 2\mathrm{Cap}_q(H_{\eps_j,b}^\omega, D_{\eps_j,b}^\omega)\,.
\end{equation}
Analogously, for every $x_{j,i}^\omega\in {(MG)}_{j,\theta,M}^\omega$, let $\phi_{0,j,i}^{\omega}\in W_0^{1,q}({B}_{\eps_jr^\omega_{\eps_j}}(\eps_jx_{j,i}^\omega);[0,1])$ be such that 
$$
{\phi_{0,j,i}^{\omega}|}_{\overline{B}_{\alpha_{\eps_j}\rho_{j,i}^\omega}(\eps_jx_{j,i}^\omega)}\equiv 1\,, \quad \int_{{B}_{\eps_jr^\omega_{\eps_j}}(\eps_jx_{j,i}^\omega)}|\nabla \phi_{0,j,i}^\omega|^q\,\mathrm{d}x\leq 2\mathrm{Cap}_q({B}_{\alpha_{\eps_j}\rho_{j,i}^\omega}(\eps_jx_{j,i}^\omega), {B}_{\eps_jr^\omega_{\eps_j}}(\eps_jx_{j,i}^\omega))\,.
$$
Note that, all the functions $\phi_{0,i,j}^{\omega}$ have disjoint supports. Proceeding as in the proof of \eqref{GHV-lemma}, we estimate 
\begin{align*}
\mathrm{Cap}_q({B}_{\alpha_{\eps_j}\rho_{j,i}^\omega}(\eps_jx_{j,i}^\omega), {B}_{\eps_jr^\omega_{\eps_j}}(\eps_jx_{j,i}^\omega)
&\leq c_{n,q}\big((\alpha_{\eps_j}\rho_{j,i}^\omega)^{(q-n)/(q-1)}-(\eps_jr_{\eps_j}^\omega)^{(q-n)/(q-1)}\big)^{1-q} \\
&\lesssim \eps_j^n(\rho_{j,i}^\omega)^{n-q}.
\end{align*}
Set 
$$
\phi_{2,j,\theta,M}^{\omega}:= 1- \sum_{x_{j,i}^\omega\in {(MG)}_{j,\theta,M}^\omega}\phi_{0,j,i}^{\omega}\,;
$$
then $\phi_{2,j,\theta,M}^\omega\in W^{1,q}(D;[0,1])$,
\begin{equation*}
\phi_{2,j,\theta,M}^{\omega}\equiv 0\,\, \text{ on }\bigcup_{x_{j,i}^\omega\in {(MG)}_{j,\theta,M}^\omega}\overline{B}_{\alpha_{\eps_j}\rho_{j,i}^\omega}(\eps_jx_{j,i}^\omega)\,, \,\, \phi_{2,j,\theta,M}^{\omega}\equiv 1 \,\, \text{on }D\setminus \bigcup_{x_{j,i}^\omega\in {(MG)}_{j,\theta,M}^\omega}\overline{B}_{\eps_jr^\omega_{\eps_j}}(\eps_jx_{j,i}^\omega)\,,
\end{equation*}
and, since the functions $\phi_{0,j,i}^{\omega}$ have disjoint supports (cf. \eqref{c:13_marzo}),
\begin{align}\label{eq: phi_2_j-cap}
\int_{D}|\nabla \phi_{2,j,\theta,M}^\omega|^q\,\mathrm{d}x = \sum_{x_{j,i}^\omega\in {(MG)}_{j,\theta,M}^\omega}
\int_{{B}_{\eps_jr^\omega_{\eps_j}}(\eps_jx_{j,i}^\omega)}|\nabla \phi_{0,j,i}^\omega|^q\,\mathrm{d}x
\lesssim \eps_j^n\sum_{x_{j,i}^\omega\in {(MG)}_{j,\theta,M}^\omega}(\rho_{j,i}^\omega)^{n-q}\,.
\end{align}
Setting now
\begin{equation}\label{eq: correctors}
\phi_{j,\theta,M}^\omega:=\phi_{1,j}^\omega\wedge\phi_{2,j,\theta,M}^\omega \in W^{1,q}(D;[0,1])\,,
\end{equation}
it is clear by \eqref{eq: phi_1_j}-\eqref{eq: correctors} that $\phi_{j,\theta,M}^\omega$ satisfies the claims $(i)$ and $(ii)$. To prove $(iii)$, since $\phi_{j,\theta,M}^\omega|_{\partial D}\equiv 1$, by the Poincar{\'e} inequality in $D$ it suffices to verify that 
\begin{equation}\label{eq: correctors_q_energy_to_0}
\lim_{M\to +\infty}\lim_{\theta\to 0^+}\lim_{j\to +\infty}\int_D|\nabla \phi_{j,\theta,M}^\omega|^q\,\mathrm{d}x=0\,.
\end{equation}
By \eqref{eq: phi_1_j} and \eqref{eq: phi_2_j-cap}  it is straightforward to estimate
\begin{align*}
\int_D|\nabla \phi_{j,\theta,M}^\omega|^q\,\mathrm{d}x&\leq\int_{D}|\nabla \phi_{1,j}^\omega|^q\,\mathrm{d}x+\int_{D}|\nabla \phi_{2,j,\theta,M}^\omega|^q\,\mathrm{d}x \nonumber\\
&\lesssim \mathrm{Cap}_q(H_{\eps_j,b}^\omega,D_{\eps_j,b}^\omega)+\eps^n\sum_{x_{j,i}^\omega\in {(MG)}_{j,\theta,M}^\omega}(\rho_{j,i}^\omega)^{n-q}\,.
\end{align*}
 The claim \eqref{eq: correctors_q_energy_to_0} then follows from \eqref{GHV-lemma} and \eqref{eq: mildly_good_zero_capacity}\,. 
\end{proof}

We are now ready for the proof of the limsup inequality.

\begin{proof}[Proof of Proposition \ref{Gamma_limsup_ineq}]
For the proof of \eqref{limsup_inequality_to_show} we follow broadly the approach developed for the periodic setting in \cite[Section 6]{Ansini-Braides} and \cite{Erratum-Ansini-Braides}. We note, however, that extra care has  to be taken here for the construction of an admissible recovery sequence to guarantee that it vanishes in the set $H_\eps^\omega$, due to the lack of separation of the perforations. We give the whole argument in detail for the sake of completeness.

Let $(u_j)\subset W_0^{1,q}(D;\R^m)$ be such that
\begin{equation}\label{sequence_attaining_the_envelope-1}
u_j {\rightharpoonup} \, u\quad \text{weakly in} \  W^{1,q}(D;\R^m) \text{ as } j\to +\infty\,,
\end{equation}
and 
\begin{equation}\label{sequence_attaining_the_envelope-2}
\lim_{j\to +\infty}\int_D f(\nabla u_j)\,\mathrm{d}x=\int_D Qf(\nabla u)\,\mathrm{d}x\,.
\end{equation}
Note that the existence of a sequence $(u_j)$ satisfying \eqref{sequence_attaining_the_envelope-1} and \eqref{sequence_attaining_the_envelope-2} follows by the characterization of $\int_D Qf(\nabla \cdot)\, \mathrm{d}x$ as the lower semicontinuous envelope of $\int_D f(\nabla \cdot)\, \mathrm{d}x$ with respect to the weak $W^{1,q}(D;\R^m)$-topology. By \cite[Lemma C.5, Remark C.6]{Braides-De Franceschi} and \cite[Lemma 1.2]{Fonseca-Mueller-Pedregal} we may also suppose, without loss of generality, that $(u_j)\subset W_0^{1,\infty}(D;\R^m)$ and that $(|\nabla u_j|^q)$ are equi-integrable. We now modify the sequence $(u_j)$ to make it admissible for the energy $\mathcal{F}^\omega_{\eps_j}$. In particular, we need to set it to zero on the perforations and on $\partial D$. 

Let $j\in \N$ be large enough so that \eqref{larger_j} is satisfied, where $K_j$ is defined as in \eqref{K_j_radii}. Let $k,M\in \N$, $\theta\in (0,1)$ (with $1/\theta\in \mathbb{N}$) be fixed and assume that also \eqref{even_smaller_theta} holds true.  
Let now $\theta':=\frac{4\theta M}{3}$; note that $\theta'\in (0, 1/2)$ by  \eqref{even_smaller_theta}. We apply Lemma \ref{l:probabilistic_joining_lemma}, with the $\theta'$ defined above, to the sequence $(u_j)$, but where $G_{j,M}^\omega$ is replaced by the smaller set $(VG)_{j,\theta,M}^\omega$ defined in \eqref{eq: very_good_centers}. Then we find a modified sequence $(w_j^\omega)\subset W_0^{1,q}(D;\R^m)$ satisfying 
\begin{enumerate}[label=(\roman*)]
\item $w^\omega_j\equiv u_j \ \text{ in } D\setminus \bigcup_{x^\omega_{j,i} \in (VG)_{j,\theta,M}^\omega} C_{\eps_j,\theta', M}^{0}(\e_j x^\omega_{j,i})$;
\smallskip
\item $w_j^\omega\equiv \bar u^\omega_{j,i} \ \text{ on } \partial B_{ \bar \sigma^\omega_{j,i}}(\eps_j x^\omega_{j,i})$, for every $x^\omega_{j,i} \in (VG)_{j,\theta,M}^\omega$, where 
\begin{equation}\label{intermediate_values-limsup}
\bar u^\omega_{j,i}:=\fint_{C_{\eps_j,\theta', M}^{0}(\e_j x^\omega_{j,i})} u_j\,\mathrm{d}x\,, \qquad  
\bar \sigma^\omega_{j,i}:=\frac{3}{4}\frac{\theta'\eps_j}{M}=\theta \eps_j\,;
\end{equation}
\item $w^\omega_j\rightharpoonup u \ \text{ weakly in } W^{1,q}(D;\R^m)$;
\item $\Big|\int_{D} f(\nabla w^\omega_j)\,\mathrm{d}x-\int_{D}f(\nabla u_j)\,\mathrm{d}x\Big|\leq C_k(\theta')^n+\frac{C}{k}\,$,  where $C_k>0$ can blow up as $k\to +\infty$;
\item $(w_j^\omega)\subset W^{1,\infty}(D;\R^m)$ and $\|w_j^\omega\|_{L^\infty(D)}\leq \|u_j\|_{L^\infty(D)}$\,.
\end{enumerate}
Note that, since $(|\nabla u_j|^q)$ are equi-integrable, one can choose $k^\omega_{j,i}=0$ in (i)-(ii) above for every $x_{j,i}^\omega\in {(VG)}_{j,\theta,M}^\omega\subset{G}^\omega_{j,M}$. Moreover, following the explicit construction of the modified sequence  in the proof of Lemmata \ref{joining_lemma_on_randomly_perforated_domains} and \ref{l:probabilistic_joining_lemma}, one can see that also $(w_j^\omega)$ has equi-integrable $q$-gradients, and that it satisfies (v).

\smallskip

We now proceed with the proof of the $\limsup$ inequality \eqref{limsup_inequality_to_show} in two steps.

\smallskip

\textit{Step 1: Assuming that $\sup_{j\in \N}\|u_j\|_{L^{\infty}(D)}=:L<+\infty$.} 
Let $\eta>0$ be fixed. By the characterization of $\int_D Qf(\nabla \cdot)\, \mathrm{d}x$ as the lower semicontinuous envelope of $\int_D f(\nabla \cdot)\, \mathrm{d}x$ with respect to the weak $W^{1,q}(D;\R^m)$-topology, and by  \eqref{j_capacity_densities} and \eqref{K_j_rho_truncation_of_capacity}, for every $x^\omega_{ j,i}\in {(VG)}_{j,\theta,M}^\omega\subset  G_{j,M}^{\omega}$ we can pick a test function $\zeta_{j,i}^{\omega}$ satisfying
\begin{equation}\label{zeta_j_i_def}
\zeta_{j,i}^{\omega}-\bar u^\omega_{j,i}\in W_0^{1,q}(B_{\theta K_{j}}(0);\R^m)\,, \quad \zeta_{j,i}^{\omega}\equiv 0\ \ \text{on } \overline{B}_{\rho_{j,i}^{ \omega}}(0)\,,
\end{equation} 
such that 
\begin{equation}\label{candidate_for_capacity_zeta_j_i}
\int_{B_{\theta K_{j}}(0)}\alpha_{\eps_j}^{q}f(\alpha_{\eps_j}^{-1}\nabla \zeta_{j,i}^\omega)\,\mathrm{d}x\leq \varphi^j_{\theta,\rho_{j,i}^{ \omega}}( \bar u^\omega_{j,i})+\eta\,.
\end{equation}
We now set 
\begin{equation}\label{def_of_recovery_sequence}
v_{j,\theta,M}^{\omega}(x):= 
\begin{cases}
\smallskip
\displaystyle
w_j^\omega(x) \ \ &\text{for }\ x\in D\setminus \bigcup_{x^\omega_{ j,i}\in {(VG)}_{j,\theta,M}^\omega}B_{\theta\eps_j}(\eps_jx^\omega_{ j,i})\,,\\
\displaystyle
\zeta_{j,i}^{\omega}\big(\alpha_{\eps_j}^{-1}(x-\eps_jx^\omega_{ j,i})\big) \ \ &\text{for }\ x\in B_{\theta\eps_j}(\eps_jx^\omega_{ j,i})\,, x^\omega_{ j,i}\in {(VG)}_{j,\theta,M}^\omega\,.
\end{cases}
\end{equation}
Recall that by property \eqref{intermediate_values-limsup} and  \eqref{zeta_j_i_def} we have that $w_j^\omega \equiv \bar u^\omega_{j,i}$ on $\partial B_{\theta \eps_j}(\eps_jx^\omega_{ j,i})$ for every $x^\omega_{j,i}\in {(VG)}_{j,\theta,M}^\omega$. Hence $v_{j,\theta,M}^{\omega}\in W_0^{1,q}(D;\R^m)$, since  $(w_j^\omega)\subset W_0^{1,q}(D;\R^m)$ (and actually in $W^{1,\infty}_0(D;\R^m)$), and with no loss of generality we can assume that $B_{\theta\eps_j}(\eps_jx^\omega_{ j,i})\subset\subset D$ for every $x^\omega_{ j,i}\in (VG)_{j,\theta,M}^{\omega}$ by repeating a similar argument as the one presented in Remark \ref{on_boundary_conditions}.  Moreover, by \eqref{zeta_j_i_def} and \eqref{def_of_recovery_sequence} we have that 
\begin{equation}\label{vzero-VG}
v_{j,\theta,M}^\omega=0 \quad  \text{in } \quad \bigcup_{x^\omega_{ j,i}\in (VG)_{j,\theta,M}^{\omega}}\overline{B}_{\alpha_{\eps_j}\rho_{j,i}^\omega}(\eps x_{j,i}^\omega)\,.
\end{equation}
We now further modify $v_{j,\theta,M}^\omega$ to obtain a sequence that vanishes on all perforations. To this aim, we set
\begin{equation}\label{eq: almost_final_recovery_sequences}
U_{j,\theta,M}^\omega:=\phi_{j,\theta,M}^\omega v_{j,\theta,M}^\omega\,,
\end{equation}
where the corrector $\phi_{j,\theta,M}^\omega$ is defined as in Lemma \ref{construction_of_correctors}. 

By \eqref{eq:properties_of_correctors}$(i)$, \eqref{eq: mildly_good_centers}, \eqref{eq: very_good_centers} and \eqref{vzero-VG} we have that $U_{j,\theta,M}^\omega=0$ in $H_{\eps_j}^\omega$; moreover, with no loss of generality (up to enlarging the domain $D$) we also have that $U_{j,\theta,M}^\omega=0$ on $\partial D$. We next show that $U_{j,\theta,M}^\omega \rightharpoonup u$ weakly in $W^{1,q}(D;\R^m)$, as $j\to +\infty$, $\theta \to 0+$ and $M\to +\infty$. We do so by proving a uniform bound for the energy \eqref{random_functionals} of $(U_{j,\theta,M}^\omega)$, which guarantees a bound on its $q$-gradients by  \eqref{growth_of_f}, and hence a bound on the $W^{1,q}$-norm of $(U_{j,\theta,M}^\omega)$ by the Poincar{\'e} inequality on $D$. 

First we set for notational convenience 

\begin{equation}\label{new_bad_sets}
\tilde H_{j,\theta,M}^\omega:=D_{\eps_j,b}^\omega \cup \bigcup_{x_{j,i}^\omega\in {(MG)}_{j,\theta,M}^\omega}\overline{B}_{\eps_jr^\omega_{\eps_j}}(\eps_jx_{j,i}^\omega)\,, \qquad \tilde D_{j,\theta,M}^\omega:=\bigcup_{x_{j,i}^\omega\in {(VG)}_{j,\theta,M}^\omega}B_{\theta\eps_j}(\eps_jx_{j,i}^\omega)\,.
\end{equation}
\EEE
Note that, by \eqref{def:good-bad} and\eqref{c:13_marzo}, 
\begin{align}\label{eq:volume_of_bad_region_to_0}
\L^n(\tilde H_{j,\theta,M}^\omega)&\lesssim \sum_{x_{j,i}^\omega\in I_{\eps_j,b}^\omega}(\alpha_{\eps_j}\rho_{j,i}^\omega)^{n}+\sum_{x_{j,i}^\omega\in {(MG)}_{j,\theta,M}^\omega}(\eps_jr^\omega_{\eps_j})^{n}\nonumber\\
&\leq \sum_{x_{j,i}^\omega\in I_{\eps_j,b}^\omega}(\eps_j^n(\rho_{j,i}^\omega)^{n-q})^{\frac{n}{n-q}}+\sum_{x_{j,i}^\omega\in \Phi_{\eps_j}^\omega(D)}(\eps_jr^\omega_{\eps_j})^{n}\nonumber\\
&\leq \bigg(\eps_j^n\sum_{x_{j,i}^\omega\in I_{\eps_j,b}^\omega}(\rho_{j,i}^\omega)^{n-q}\bigg)^{\frac{n}{n-q}}+(\eps_j^nN^\omega_{\eps_j}(D))(r_{\eps_j}^\omega)^n\,,
\end{align} 
where we used that, for $\alpha\geq 1$ and $t_i \geq 0$ for every $i\in \N$,
\[\sum_{i\in \N}t_i^\alpha\leq \Big(\sum_{i\in \N}t_i\Big)^\alpha\,.\] 
By \eqref{r_eps_are_negligible}--\eqref{c:13_marzo} we have that $\L^n(\tilde H_{j,\theta,M}^\omega) \to 0$ as $j\to +\infty$. Analogously,
\begin{equation}\label{volume_of_capacitary_region_to_0}
\L^n(\tilde D_{j,\theta,M}^\omega)\leq \theta^n(\eps_j^nN^\omega_{\eps_j}(D))\,.
\end{equation} 
Finally we estimate the energy \eqref{random_functionals} of $(U_{j,\theta,M}^\omega)$. Using \eqref{def_of_recovery_sequence}--\eqref{eq: almost_final_recovery_sequences} and \eqref{growth_of_f}, we estimate
\begin{align}\label{energy_upper_bound}
\int_{D}f(\nabla U_{j,\theta,M}^\omega)\,\mathrm{d}x&=\int_{D\setminus \tilde D_{j,\theta,M}^\omega}f(\nabla U_{j,\theta,M}^\omega)\,\mathrm{d}x+ \int_{\tilde D_{j,\theta,M}^\omega}f(\nabla U_{j,\theta,M}^\omega)\,\mathrm{d}x\nonumber\\
&\leq \int_{D\setminus(\tilde D_{j,\theta,M}^\omega\cup\tilde H_{j,\theta,M}^\omega)}f(\nabla w_{j}^\omega)\,\mathrm{d}x+\sum_{x_{j,i}^\omega\in {(VG)}_{j,\theta,M}^\omega}\int_{B_{\theta\eps_j}
(\eps_jx^\omega_{ j,i})}f(\nabla U_{j,\theta,M}^{\omega })\,\mathrm{d}x\nonumber\\
&\quad +c_2\int_{\tilde H_{j,\theta,M}^\omega}\big(|\phi_{j,\theta,M}^\omega\nabla v_{j,\theta,M}^\omega+\nabla\phi_{j,\theta,M}^\omega\otimes v_{j,\theta,M}^\omega|^q+1|\big)\,.
\end{align}
We deal with each term on the right-hand side of \eqref{energy_upper_bound} separately. For the first term, by \eqref{growth_of_f} and by property (iv) of $(w_j^\o)$,
\begin{align}\label{aux_1_energy_upper_bound}
\int_{D\setminus(\tilde D_{j}^\omega\cup\tilde H_{j,\theta,M}^\omega)}f(\nabla w_{j}^\omega)\,\mathrm{d}x&\leq\int_D f(\nabla w_{j}^\omega)\,\mathrm{d}x+c_1\int_{\tilde D_{j,\theta,M}^\omega\cup\tilde H_{j,\theta,M}^\omega}(1-|\nabla w_{j}^\omega|^q)\,\mathrm{d}x\nonumber\\
&\leq \int_D f(\nabla u_{j})\,\mathrm{d}x+C_k(\theta')^n+\frac{C}{k}+c_1\L^n(\tilde D_{j,\theta,M}^\omega\cup\tilde H_{j,\theta,M}^\omega)\,.
\end{align}	
For the second term on the right-hand side of \eqref{energy_upper_bound},  by a simple change of variables, \eqref{eq:properties_of_correctors}$(ii)$, \eqref{candidate_for_capacity_zeta_j_i} and \eqref{def_of_recovery_sequence}, for every $x^\omega_{ j,i} \in {(VG)}_{j,\theta,M}^{\omega}$ we have 
\begin{align*}
\begin{split}
\int_{B_{\theta\eps_j}
(\eps_jx^\omega_{ j,i})}f(\nabla U_{j,\theta,M}^{\omega }(x))\,\mathrm{d}x&=\int_{B_{\theta\eps_j}(\eps_jx^\omega_{ j,i})}f(\alpha_{\eps_j}^{-1}\nabla \zeta_{j,i}^{\omega}\big(\alpha_{\eps_j}^{-1}(x-\eps_jx^\omega_{ j,i})\big))\,\mathrm{d}x\\
&=\eps_j^n\int_{B_{\theta K_{j}}(0)}\alpha_{\eps_j}^qf(\alpha_{\eps_j}^{-1}\nabla \zeta_{j,i}^{\omega}(y))\,\mathrm{d}y\leq \eps_j^n\varphi^j_{\theta,\rho_{j,i}^{\omega}}( \bar u^\omega_{j,i})+\eta\eps_j^n\,.
\end{split}
\end{align*}
Moreover, by \eqref{eq: very_good_centers} and  \eqref{lower_bound_for_the_truncated_capacities}, 
\begin{align*}
\eps_j^n\sum_{x_{j,i}^\omega\in {(VG)}_{j,\theta,M}^\omega}\varphi^j_{\theta,\rho_{j,i}^{\omega}}( \bar u^\omega_{j,i})&=\eps_j^n\sum_{x_{j,i}^\omega\in {G}_{j,M}^\omega}\varphi^j_{\theta,\rho_{j,i}^{\omega}}( \bar u^\omega_{j,i})-\eps_j^n\sum_{x_{j,i}^\omega\in {(MG)}_{j,\theta,M}^\omega\cap G_{j,M}^\omega}\varphi^j_{\theta,\rho_{j,i}^{\omega}}( \bar u^\omega_{j,i})\nonumber\\
&\leq \eps_j^n\sum_{x_{j,i}^\omega\in {G}_{j,M}^\omega}\varphi^j_{\theta,\rho_{j,i}^{\omega}}( \bar u^\omega_{j,i})+C_2\theta^n(\eps_j^nN_{\eps_j}^\omega(D))\,.
\end{align*}	
Finally, we consider the third term on the right-hand side of \eqref{energy_upper_bound}. By the fact that $0\leq \phi_{j,\theta,M}^\omega\leq 1$, and that $v_{j,\theta,M}^\omega=w_j^\omega$ in $\tilde H_{j,\theta,M}^\omega$ by \eqref{def_of_recovery_sequence} and  \eqref{new_bad_sets}, we have 
\begin{align}\label{aux_3_energy_upper_bound}
\int_{\tilde H_{j,\theta,M}^\omega}\big(|\phi_{j,\theta,M}^\omega\nabla v_{j,\theta,M}^\omega+\nabla\phi_{j,\theta,M}^\omega\otimes v_{j,\theta,M}^\omega|^q+1\big)\,\mathrm{d}x&\lesssim \int_{\tilde H_{j,\theta,M}^\omega}\big(|\nabla w_{j}^\omega|^q+L^q|\nabla\phi_{j,\theta,M}^\omega|^q+1\big)\,\mathrm{d}x\,,
\end{align}	
where we have used that $\|w^\o_j\|_{L^{\infty}(D)}\leq \|u_j\|_{L^{\infty}(D)}\leq L<+\infty$, by property (v) of $(w_j^\o)$ and by the assumption in this step.

The upper bound estimate \eqref{energy_upper_bound}, together with \eqref{aux_1_energy_upper_bound}--\eqref{aux_3_energy_upper_bound}, \eqref{eq:volume_of_bad_region_to_0}, \eqref{volume_of_capacitary_region_to_0},  yields
\begin{align}\label{almost_final_upper_energy}
\int_D f(\nabla U_{j,\theta,M}^\omega)\,\mathrm{d}x	&\leq \int_D f(\nabla u_j)+\eps_j^n\sum_{x_{j,i}^\omega\in G_{j,M}^\omega}\varphi^j_{\theta,\rho_{j,i}^\omega}(\bar u_{j,i}^\omega)\nonumber\\
&\quad +C\int_{\tilde H_{j,\theta,M}^\omega}|\nabla w_{j}^\omega|^q\,\mathrm{d}x+CL^q\int_{D}|\nabla \phi_{j,\theta,M}^\omega|^q\,\mathrm{d}x\nonumber\\
&\quad +C\L^n(\tilde H_{j,\theta,M}^\omega)+C(\eta+\theta^n)(\eps_j^nN^\o_{\eps_j}(D))+C_k(\theta')^n+\frac{C}{k}\,.
\end{align}
By \eqref{sequence_attaining_the_envelope-2}, \eqref{c:c-star}, the growth condition \eqref{growth_of_f}, property (iv) for $(w_j^\o)$, \eqref{eq: correctors_q_energy_to_0}, \eqref{eq:volume_of_bad_region_to_0}, and \eqref{r_eps_are_negligible}--\eqref{c:13_marzo} 
it follows that
\begin{equation*}
\sup_{M,\theta,j}\int_D f(\nabla U_{j,\theta,M}^\omega) \dx  \leq C_{L}<+\infty\,,
\end{equation*}
and hence 
\begin{equation*}
\sup_{M,\theta,j}\|U_{j,\theta,M}^\omega\|_{W^{1,q}(D;\R^m)}\leq C_{L}<+\infty\,,
\end{equation*}
for a constant $C_{L}>0$ independent of $\omega\in \Omega$.  By \eqref{sequence_attaining_the_envelope-1}, proceeding as for the proof of \eqref{v_j_M_L_omega_convergence_is_ok}, we deduce that for $\mathbb{P}$-a.e. $\omega\in \Omega$, 
\begin{equation*}
U_j^{\theta,M,\omega}\rightharpoonup u\ \text{weakly in } W^{1,q}(D;\R^m) \ \text{as } j\to +\infty,\ \theta\to 0^+, \ M\to +\infty\,.
\end{equation*}
Finally, we show that the admissible sequence $(U_j^{\theta,M,\omega})$ is a recovery sequence, up to a diagonal procedure. 
Indeed, taking in \eqref{almost_final_upper_energy} the limsup as $j\to +\infty$, and using the subadditivity of the limsup, \eqref{sequence_attaining_the_envelope-2}, the equi-integrability of $(|\nabla w_j^\o|^q)$, \eqref{eq:volume_of_bad_region_to_0} and  \eqref{1st_law_of_large_numbers_type},
we deduce
\begin{eqnarray*}
\varlimsup_{j\to+\infty}\int_D f(\nabla U_{j,\theta,M}^\omega)\,\mathrm{d}x &\leq & \int_D Qf(\nabla u)\,\mathrm{d}x
+\varlimsup_{j\to+\infty}\ \eps_j^n\sum_{x_{j,i}^\omega\in G_{j,M}^\omega}\varphi^j_{\theta,\rho_{j,i}^\omega}(\bar u_{j,i}^\omega)\nonumber\\
&\quad &+CL^q \varlimsup_{j\to+ \infty}\int_D |\nabla \phi_{j,\theta,M}^\omega|^q\,\mathrm{d}x+C(\eta+\theta^n)\langle N(Q)\rangle\L^n(D)\nonumber\\
&\quad &+C_k(\theta M)^n+\frac{C}{k}\,.
\end{eqnarray*}
Passing further to the limit superior in the above inequality in $\theta\to 0^+$, and $M\to +\infty$, using Proposition \ref{main_prop_for_approximating_the_capacity} and \eqref{eq: correctors_q_energy_to_0} we get 
\begin{equation*}
\varlimsup_{M\to+\infty}\varlimsup_{\theta\to 0^+}\varlimsup_{j\to+\infty}\int_D f(\nabla U_{j,\theta,M}^\omega)\,\mathrm{d}x\leq \int_D Qf(\nabla u)\,\mathrm{d}x+\langle N(Q)\rangle \int_D\varphi(u)\,\mathrm{d}x+C(\eta+{1}/{k})\,.
\end{equation*}
By finally choosing a diagonal sequence $(\tilde u_j^{ \omega})\subset W_0^{1,q}(D;\R^m)$, with $\tilde u_j^{ \omega }:= U_{j,\theta_j,M_j}^\omega$, we conclude that for $\mathbb{P}$-a.e. $\omega \in \Omega$ $(\tilde u_j^\omega)$ is admissible$\ie$ $\tilde u_j^\omega|_{H_{\eps_j}^\omega\cup \partial D}\equiv0$, and 
$$\limsup_{j\to +\infty}\int_{D}f(\nabla \tilde u_j^{\omega})\,\mathrm{d}x\leq \int_{D}Qf(\nabla u)\,\mathrm{d}x+\langle N(Q)\rangle\int_D\varphi(u)\,\mathrm{d}x+ C(\eta+{1}/{k})\,,$$
where we have used that for $\mathbb{P}$-a.e. $\omega \in \Omega$ the sequence $(\tilde u_j^{ \omega})$ converges to $u$ weakly in $W^{1,q}(D;\R^m)$, and hence also in $L^1(D;\R^m)$. 

Hence, \eqref{limsup_inequality_to_show} follows by the arbitrariness of $k\in \N$ and $\eta\in (0,1)$, with $(\tilde u_j^{ \omega})$ as a recovery sequence.

\smallskip

\textit{Step 2: Removing the uniform $L^\infty$-bound on $(u_j)$.} This is done by following verbatim the final argument in \cite{Ansini-Braides}. Indeed, assume first that $u\in W_0^{1,q}(D;\R^m)\cap L^\infty(D;\R^m)$, and let $L:=4\|u\|_{L^\infty(D)}$. Let $\Psi_L\colon \R^m\mapsto \R^m$ be a Lipschitz function with Lipschitz constant at most $1$ such that
\begin{equation*}
\Psi_L(z):=\begin{cases} z \quad \text{if }\ |z|\leq L/2\,,\\
0 \quad \text{if }\ |z|\geq L\,.
\end{cases}
\end{equation*}		
Let again $(u_j)\in W_0^{1,q}(D;\R^m)$ be such that \eqref{sequence_attaining_the_envelope-1}-\eqref{sequence_attaining_the_envelope-2} hold true. Without loss of generality we assume that $u_j\to u$ pointwise $\L^n$-a.e in $D$ as $j\to +\infty$, and that $(|\nabla u_j|^q)$ is equi-integrable on $D$. We now set $u_j^L:=\Psi_L(u_j)$. By the $\mathcal{L}^n$-a.e.  pointwise convergence $u_j$ to $u$ we have that 
$$\L^n\big(\{u_j^L\neq u_j\}\big)\leq \L^n\big(\{|u_j|>2\|u\|_{L^\infty(D)}\geq \|u\|_{L^{\infty}(D)}\}\big)\longrightarrow 0 \ \text{as \ } j\to +\infty\,,$$
hence, we also have that $u_j^L\rightharpoonup u$ weakly in $W^{1,q}(D;\R^m)$ as $j\to +\infty$. Furthermore, by the equi-integrability of $(|\nabla u_j|^q)$ we obtain
$$\lim_{j\to +\infty}\int_D f(\nabla u_j^L)\,\mathrm{d}x=\lim_{j\to +\infty}\int_D f(\nabla u_j)\,\mathrm{d}x=\int_D Qf(\nabla u)\,\mathrm{d}x\,, $$
and we can therefore repeat all the reasonings of Step 1, with $(u_j^L)$ in the place of $(u_j)$ in order to obtain the desired limsup inequality.

Finally, given an arbitrary $u\in W_0^{1,q}(D;\R^m)$, we can approximate it by a sequence $(u^L)\subset W_0^{1,q}(D;\R^m)\cap L^\infty(D;\R^m)$ with respect to the strong $W^{1,q}(D,\R^m)$-topology. Then, for $\mathbb{P}$-a.e. $\omega\in \Omega$ by the lower semicontinuity of the functional 
$$\F_0^{''\omega}:=\Gamma\text{-}\limsup_{j\to +\infty}\F_{\eps_j}^{\omega}$$
with respect to the strong $L^q(D;\R^m)$-topology (cf. \cite[Remark 7.8]{Braides-De Franceschi}), we obtain
$$\F_0^{'' \omega}(u)\leq \liminf_{L\to +\infty}\F_0^{'' \omega}(u^L)=\lim_{L\to +\infty}\F_0(u^L)=\F_0(u)\,,$$
which, by the definition of $\F_0^{'' \omega}$, is just another way of writing the desired limsup inequality in the general case. The proof is now complete. \qedhere

\end{proof}

\section*{Acknowledgements}  L. Scardia was supported by the EPSRC grants EP/V00204X/1 and EP/V008897/1. \EEE
The work of K. Zemas and C. I. Zeppieri was supported by the Deutsche Forschungsgemeinschaft (DFG, German Research Foundation) under Germany's Excellence Strategy EXC 2044 -390685587, Mathematics M\"unster: Dynamics--Geometry--Structure.

\typeout{References}


\begin{thebibliography}{50}
\bibitem{Allaire}
{\sc G.~Allaire}.
\newblock {Homogenization of the Navier-Stokes equations in open sets perforated with tiny holes. II. Noncritical size of the holes for a volume distribution and a surface distribution of holes}. 
\newblock {\em Arch. Rational Mech. Anal.}, {\bf 113(3)} (1990), 261--298.
	
	
\bibitem{Ansini-Babadjian-Zeppieri}
{\sc N.~Ansini, J. F.~Babadjian, C. I. Zeppieri}.
\newblock {The Neumann sieve problem and dimension reduction: a multiscale approach}. 
\newblock {\em Math. Models Methods Appl. Sci.}, {\bf 17(05)} (2007), 681--735.

\bibitem{Ansini-Braides}
{\sc N.~Ansini, A.~Braides}.
\newblock {Asymptotic analysis of periodically-perforated nonlinear media}. 
\newblock {\em J. Math. Pures Appl.}, {\bf 81(5)} (2002), 439--451.
	
\bibitem{Erratum-Ansini-Braides}
{\sc N.~Ansini, A.~Braides}.
\newblock {Erratum to Asymptotic analysis of periodically-perforated nonlinear media}. 
\newblock {\em J. Math. Pures Appl.}, {\bf 81.5} (2002), 439--451.
	
\bibitem{Braides-De Franceschi}
{\sc A.~Braides, A.~Defranceschi}.
\newblock {\em Homogenization of multiple integrals}. 
\newblock Oxford University Press, {\bf Vol. 12}, New York (1998).
	
\bibitem{vitalibraides1996homogenization}
{\sc A.~Braides, A.~Defranceschi, E.~Vitali}.
\newblock {Homogenization of free discontinuity problems}.
\newblock {\em Arch. Rational Mech. Anal.}, {\bf 135(4)} (1996), 297--356.
	
	
\bibitem{Caffarelli_Mellet}
{\sc L. A.~Caffarelli, A.~Mellet}.
\newblock {Random homogenisation of an obstacle problem}.
\newblock {\em Ann. Inst. H. Poincar\'e Anal. Non Lin\'eaire}, 
\newblock {\bf 26(2)} (2009), 375--395.
	
	
\bibitem{Casado-Diaz}
{\sc C.~Calvo-Jurado, J.~Casado-Diaz, M.~Luna-Laynez}.
\newblock {Homogenization of nonlinear Dirichlet problems in random perforated domains}.
\newblock {\em Nonlinear Anal.}, 
\newblock {\bf 133} (2016), 250--274.
	
	
\bibitem{Casado-Diaz2}
{\sc C.~Calvo-Jurado, J.~Casado-Diaz, M.~Luna-Laynez}.
\newblock{Homogenization of the Poisson equation with Dirichlet conditions in random perforated domains}.
\newblock {\em J. Comp. App. Math.}, 
\newblock {\bf 275} (2016), 375--381.			
					

\bibitem{Casado-Diaz-Garroni-2}
{\sc J.~Casado-Diaz, A.~Garroni}.
Asymptotic behaviour of nonlinear elliptic systems on varying domains. \emph{SIAM J. Math. Anal.}, {\bf 31} (2000), 581--624.	
	
\bibitem{Cionarescu_strange}	
{\sc D.~Cioranescu, F.~Murat}.
\newblock {Un terme \'etrange venu d'ailleurs II}.
\newblock {\em Nonlinear partial differential equations and their applications, Coll\'ege de France Seminar, Res. Notes in Math.}
\newblock {\bf Vol. III} (1982).

\bibitem{Daley-Vere-Jones1}
{\sc D. J.~Daley, D.~Vere-Jones}.
\newblock {\em An introduction to the theory of point processes. Vol.I: Elementary theory and
methods}.
\newblock Probability and its applications, Springer-Verlag New York (2003).

	
	
\bibitem{Daley-Vere-Jones}
{\sc D. J.~Daley, D.~Vere-Jones}.
\newblock {\em An introduction to the theory of point processes. Vol.II: General theory and
structures}.
\newblock Probability and its applications, Springer-Verlag New York (2008).

\bibitem{DalMaso-Defranceschi}
{\sc G. Dal Maso, A. Defranceschi}. Limits of nonlinear Dirichlet problems in varying domains. \emph{Manuscripta Math.}, {\bf 61} (1988), 251--278.  

\bibitem{DalMaso-Garroni}
{\sc G. Dal Maso, A. Garroni}. 
\newblock{New results on the asymptotic behaviour of Dirichlet problems in perforated domains.} 
\emph{Math. Models Methods Appl. Sci.}, {\bf{4}(3)} (1994), 373--407.  



\bibitem{Federer-Ziemer}
{\sc H.~Federer, W.~Ziemer}.
\newblock {The Lebesgue set of a function whose distribution derivatives are $p$-th power summable}.
\newblock {\em Indiana Univ. Math. J.}, {\bf 22(2)} (1972), 139--158.	


\bibitem{Focardi}
{\sc M.~Focardi}. Homogenisation of random fractional obstacle problems via $\Gamma$-convergence. \emph{Comm. in PDEs}, {\bf 34(12)} (2009), 1607--1631. 

\bibitem{Fonseca-Mueller-Pedregal}
{\sc I.~Fonseca, S.~M\"uller, P.~Pedregal}.
\newblock {Analysis of concentration and oscillation effects generated by gradients}.
\newblock {\em SIAM J. Math. Anal.}, {\bf 29(3)} (1996), 736--756.	

\bibitem{Giunti}
{\sc A.~Giunti.} Convergence rates for the homogenization of the Poisson problem in randomly perforated domains. \emph{Netw. Heterog. Media}, {\bf 16(3)} (2021), 341--375. 

	
\bibitem{Giunti-Hofer-Velasquez} 
{\sc A.~Giunti, R.~H\"ofer, J.~Vel\'azquez}.
\newblock {Homogenization for the Poisson equation in randomly perforated domains under minimal assumptions on the size of the holes}. 
\newblock {\em Comm. Partial Differential Equations}, {\bf 43(9)} (2018), 1377--1412. 
	
	
\bibitem{Marchenko_Khruslov}
{\sc A. V.~Marchenko, E. Y.~Khruslov}.
\newblock {\em Boundary Value Problems in Domains with Fine-Granulated Boundaries}.
\newblock {Naukova Dumka}, Kiev (1974).
	
	
\bibitem{Marchenko_Khruslov2}
{\sc A. V.~Marchenko, E. Y.~Khruslov}.
\newblock {\em Homogenization of partial differential equations}.
\newblock {Progress in Mathematical Physics}
\newblock {\bf 46}, Birkhäuser Boston, Inc., MA (2006).
	
\bibitem{Papa}
{\sc G. C.~ Papanicolaou, S. R. S.~Varadhan}.
\newblock {Diffusion in regions with many small holes}.
\newblock {\em Springer Berlin Heidelberg} (1980), 190--206.
	

\bibitem{Schneider-Weil}
{\sc R.~Schneider, W.~Weil}.
\newblock {\em Stochastic and Integral Geometry}.
\newblock Probability and its applications, Springer-Verlag New York (2008).


\bibitem{PelScaZep}
{\sc X.~Pellet, L.~Scardia, C. I.~Zeppieri}.
\newblock {Stochastic homogenisation of free-discontinuity functionals in random perforated domains}.
\newblock Accepted in {\em Adv. Calc. Var.} (2023). 
\end{thebibliography}
 \end{document}